\documentclass{article}

\usepackage[T1]{fontenc}
\usepackage[a4paper,hmargin=2.5cm]{geometry}

\usepackage{csquotes}
\usepackage{amsmath,amsfonts,amsthm}
\usepackage{dsfont}

\usepackage{textcomp}
\usepackage{gensymb} 

\usepackage{graphicx} 
\usepackage{subcaption}

\usepackage{orcidlink}

\usepackage{pgfplots}
\pgfplotsset{compat=1.18}
\usepackage{tikz}
\usetikzlibrary{patterns}
\usetikzlibrary{patterns.meta}
\usetikzlibrary{calc}
\usetikzlibrary{arrows.meta,bending}

\usepackage{standalone}

\usepackage{enumitem}
\usepackage{booktabs}   
\usepackage{hyperref}

\usepackage
[
backend=biber,
style=alphabetic,
maxnames=5, 
giveninits=true, 
isbn=false,
backref=true
]
{biblatex}

\DeclareSourcemap{
	\maps[datatype=bibtex]{
		\map[overwrite]{
			\step[fieldsource=doi, final]
			\step[fieldset=url, null]
			\step[fieldset=eprint, null]
		}  
	}
}

\usepackage{todonotes}
\colorlet{cl_mh}{cyan}
\colorlet{cl_ck}{orange!80!black}
\colorlet{cl_js}{green!80!black}

\def\bu{\mathbf{u}}

\def\bv{\mathbf{v}}

\def\bw{\mathbf{w}}

\def\bg{\mathbf{g}}
\def\bn{\mathbf{n}}

\def\bJ{\boldsymbol{J}}

\def\bU{\boldsymbol{U}}

\def\du{\,\mathrm{d}}
\def\D{\,\mathrm{D}}
\newcommand{\dx}{\du x}

\DeclareMathOperator{\dive}{div}

\newcommand{\alphi}[1]{\alpha_{\varphi}^{(\delta)}(#1)}
\newcommand{\alphu}[1]{\alpha_{\bu}^{(\delta)}(#1)}

\theoremstyle{plain}
\newtheorem{theorem}{Theorem}[section]
\newtheorem{proposition}[theorem]{Proposition}
\newtheorem*{proposition*}{Proposition}
\newtheorem{lemma}[theorem]{Lemma}

\newtheorem{corollary}[theorem]{Corollary}

\theoremstyle{definition}
\newtheorem{definition}[theorem]{Definition}
\newtheorem{assumption}[theorem]{Assumption}
\newtheorem{notation}[theorem]{Notation}

\theoremstyle{remark}
\newtheorem{remark}[theorem]{Remark}

\title{A porous media approximation for an obstacle immersed in a binary fluid governed by thermodynamically consistent Cahn--Hilliard/Navier--Stokes equations
}

\author{Michael Hinze%
	\footnote{Mathematical Institute, University of Koblenz, Germany}
	\orcidlink{0000-0001-9688-0150},\,%
	\and Christian Kahle%
	\footnotemark[1]
	\orcidlink{0000-0002-3514-5512},\,%
	\and John Sebastian H. Simon%
    \footnote{Institute of Mathematics, University of the Philippines Diliman, Quezon City, Philippines}
	\orcidlink{0000-0002-7711-9582}
}
\date{\today}

\begin{document}
	
\maketitle

\begin{abstract}
		We propose a diffuse-interface model for two-phase flow in domains with obstacles. Both the fluid phases and the obstacle are represented by phase-field variables. The model combines a thermodynamically consistent Cahn--Hilliard/Navier--Stokes system with a porous medium approximation of the obstacle that extends the two-phase fluid domain into the obstacle.
        We formally justify the latter approximation through matched asymptotic analysis and formally recover the two-phase fluid system outside the obstacle.
        A stable semi-discrete time stepping scheme is introduced and existence of solutions together with global energy estimates are established. Numerical simulations illustrate the consistency of the approximation and its convergence toward the setting without porous medium approximation of the obstacle.
		
		
\end{abstract}

\section{Introduction}
In this work we provide and analyze a mathematical model for two-phase flow over an obstacle, where the flow is governed by a Cahn--Hilliard/Navier--Stokes system with variable densities and viscosities, and the obstacle is modeled/defined by a phase-field.

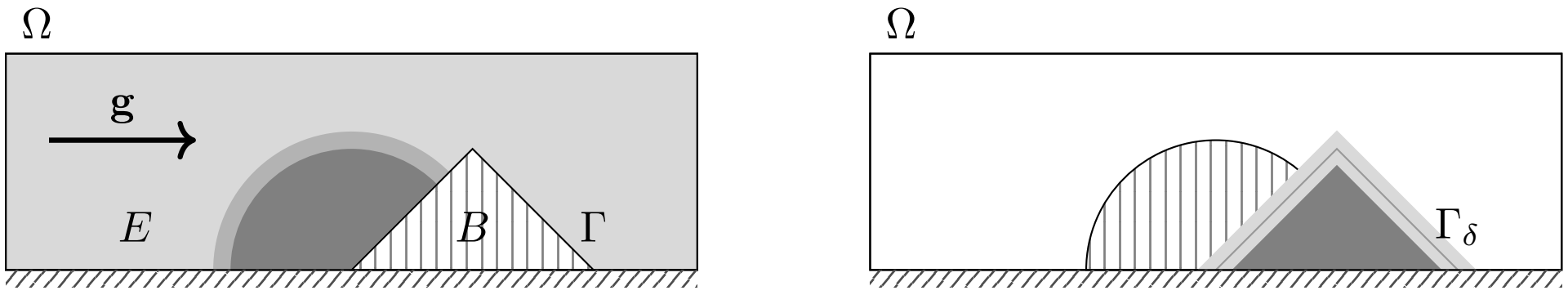
\begin{figure}
    \centering
     {
    \tikzset{every picture/.style={scale=1.4}}
    \begin{tikzpicture}

\colorlet{clFluidA}{black!15}
\colorlet{clFluidI}{black!30}
\colorlet{clFluidB}{black!50}

\begin{scope}[shift={(-5,0)}]

\coordinate (BL) at (0,0);
\coordinate (TR)  at (4,1.25);
\draw[fill=clFluidA] (BL)  rectangle (TR);

\begin{scope}[shift={(2.0,0.0)}]
\pgfmathsetmacro{\rI}{0.7}
\pgfmathsetmacro{\rO}{0.8}
    \fill[color=clFluidI] (\rO,0) arc (0:180:\rO) -- cycle;
\fill[color=clFluidB] (\rI,0) arc (0:180:\rI) -- cycle;

\end{scope}

\begin{scope}[shift={(2.0,0.0)}]
    \coordinate (A) at (0,0);
    \coordinate (B) at (0.7,0.7);
    \coordinate (C) at (1.4,0);
    \pgfmathsetmacro{\r}{0.1}

\draw[fill=white,color=white] (A)--(B)--(C)--cycle;
\draw[pattern={Lines[angle=90, distance=3.5pt,line width=0.5pt]}, pattern color=black!50]  (A)--(B)--(C)--cycle;

\node at (0.7,0.25)  {$B$};

\end{scope}

\draw (BL)  rectangle (TR);

\fill[pattern={Lines[angle=50, distance=2.5pt,line width=0.5pt]}, pattern color=black!70] (0.0,-0.1) rectangle (4.0,0.0);

\draw[very thick,->] (0.25,0.75) -- (1.1,0.75) node[midway,above] {${\mathbf g}$};

\node at (0.75,0.25)  {$E$};

\node at (3.4,0.25) {$\Gamma$};

\node[anchor=south west] at (0.0,1.25)  {$\Omega$};

\end{scope} 


\begin{scope}[shift={(0,0)}]
    
\coordinate (BL) at (0,0);
\coordinate (TR)  at (4,1.25);
\draw[fill=white] (BL)  rectangle (TR);

\begin{scope}[shift={(2.0,0.0)}]

\pgfmathsetmacro{\r}{0.75}
\draw[pattern={Lines[angle=90, distance=3.5pt,line width=0.5pt]}, pattern color=black!50] 
(\r,0) arc (0:180:\r) -- cycle;

\end{scope} 

\begin{scope}[shift={(2.0,0.0)}]
    \coordinate (A) at (0,0);
    \coordinate (B) at (0.7,0.7);
    \coordinate (C) at (1.4,0);
    \pgfmathsetmacro{\r}{0.1}

\draw[color=black!15,fill=black!15] ($(A) - (\r,0)$) -- ($(B) + (0,\r)$) -- ($(C) + (\r,0)$) -- cycle;

\draw[color=black!40] (A)--(B)--(C)--cycle;
\draw[color=black!50,fill=black!50] ($(A) + (\r,0)$) -- ($(B) - (0,\r)$) -- ($(C) - (\r,0)$) -- cycle;



\end{scope}

\draw (BL)  rectangle (TR);

\fill[pattern={Lines[angle=50, distance=2.5pt,line width=0.5pt]}, pattern color=black!70] (0.0,-0.1) rectangle (4.0,0.0);

\node at (3.4,0.25) {$\Gamma_\delta$};

\node[anchor=south west] at (0.0,1.25)  {$\Omega$};

\end{scope}

\end{tikzpicture}
    }
    \caption{Abstract setup of the situation at hand. 
        In this sketch the gravitational force $\mathbf g$ points from left to right and the sketch shows a droplet sliding down a wall over a kink that is attached to the wall.
        On the left, the hold-all domain $\Omega$ is split into a fluid domain $E$ and an obstacle domain $B$. These are separated by an interface $\Gamma$.
        Inside $E$ a two-phase fluid is present, where the distribution of the two fluids is given by a phase field variable $\varphi$. The two fluids are indicated by the very light gray domain and the dark gray circle. The middle gray shell in between indicates the interfacial region between the two fluids and admits a thickness proportional to $\varepsilon >0$.  The obstacle in this sketch is shown as a sharp object and is not part of the fluid domain.
        We approximate this sharp obstacle with a phase field $\psi$ and extend the fluid domain to the complete domain $\Omega$. On the right we show the diffuse obstacle in gray values and the zero level line of the droplet. Especially the interface $\Gamma_\delta$ between fluid domain and obstacle now has a positive thickness of order $\delta >0$.
}
    
    \label{fig:intro:schematicsetup}
  
\end{figure}

For this let us consider a domain $\Omega\subset \mathbb{R}^d$ with $d\in\{2,3 \}$, and a subdomain $B\subset \Omega$, which we shall refer to as an obstacle. 
We suppose that the domain $E := \Omega\setminus \overline{B}$ is occupied by a two-phase fluid with varying density and viscosity. 
In Figure~\ref{fig:intro:schematicsetup} we sketch a physical example, illustrating a bubble sliding along a wall and passing a triangular obstacle.

For both the description of the two-phase flow and the description of the obstacle we exploit a phase-field approach. For the binary fluid we use the thermodynamically consistent model for two-phase flow based from \cite{AGG2012-3MAS} in the variant proposed in \cite{AbelsBreit-NonlinearDensity} that can cope with nonlinear modeling of the phase-field to density relation. 
We combine this with the concepts for topology optimization for one-phase flow from \cite{BorPIJNMF2003,GarckeHecht-TopoOpt,GHHKL2016-IFB} and for Poisson equation, e.g. \cite{GHKKL2023-TopoOpt-Laplace}, to achieve a monolithic formulation of the system. By this approach, the obstacle enters as coefficient into the system.

More precisely, given a fixed phase field $\psi$ that describes the stationary obstacle, we investigate the following model for a velocity field $\bu$ and an associated pressure field $p$ together with a time dependent phase field $\varphi$, that describes the distribution of the two involved fluids, and an associated chemical potential $\mu$
\begin{subequations}
        \begin{align*}
    \frac{1}{2}\big[\rho\partial_t \bu + \partial_t(\rho\bu) + ((\rho\bu + \bJ_\delta)\cdot\nabla)\bu + \dive((\rho\bu + \bJ_\delta) \otimes \bu) \big] &\notag\\
    + \alphu{\psi}\bu- \dive(2 \eta \D\bu) + \nabla p &= -\varphi\nabla\mu + \rho\bg &&\text{ in }I\times \Omega,\\
    \dive \bu & = 0 && \text{ in }I\times \Omega,\\
    \partial_t\varphi + \dive(\varphi\bu - b\alphi{\psi} \nabla\mu ) & = 0 &&\text{ in }I\times \Omega, \\
     - \sigma{\varepsilon}\dive(\alphi{\psi}\nabla \varphi) + \frac{\sigma}{{\varepsilon}}W'(\varphi) & = \mu &&\text{ in }I\times \Omega.
\end{align*}
\end{subequations}
This system is supplemented with the boundary conditions $\bu = 0$, $\partial_\nu \varphi = \partial_\nu \mu = 0$ on $(0,T]\times \partial \Omega$, where $\nu$ denotes the outer unit normal of $\Omega$. Moreover, the initial conditions $\bu(0) = \bu^0$ and $\varphi(0) = \varphi^0$ are satisfied in $\Omega$.
Here the density $\rho(\varphi)$ and the viscosity $\eta(\varphi)$ are given functions, possibly nonlinear, that describe the spatial distribution of the density and viscosity. We abbreviate $\D\bu = \frac{1}{2}(\nabla \bu + (\nabla \bu)^T) $ as the shear stress tensor. 
The gravitational 
acceleration 
is given by $\bg$ and the function $\alphu{\psi}$ is a penalty function, 
that  takes large values inside the obstacle and thus penalizes flow inside the obstacle.
The term $\bJ_\delta$
is a function that render the model (without the penalty terms) thermodynamically consistent and is given in detail in Section~\ref{sec:model}. The term $\varphi\nabla \mu$ models capillary forces at the interface between the two fluids.
The constant $b>0$ describes a mobility of the interface between the two fluids and $\alphi{\psi}$ penalizes the equation inside the obstacle. The constants $\sigma$, $\varepsilon$ and $\delta$ respectively define the surface tension at the fluidic interface between the two fluids, the scaling of the width of the fluidic interface within the phase field $\varphi$, 
and the thickness of the interface between obstacle and fluid domain. Finally $W$ is a double-well potential, that drives the phase field $\varphi$ towards the two pure fluids. 
We note that the time derivative in the momentum equation seems uncommon. However, it appears from a reformulation of the original formulation using the third equation and allows us to derive a thermodynamically consistent time discretization scheme.
For a detailed description we refer to Section~\ref{sec:model}.

\subsubsection*{Contribution}
The main contributions of this work are as follows:
\begin{enumerate}
    \item We provide the above model that allows to simulate two-phase flow over an obstacle which itself is modeled by a phase field function.
    \item We investigate the sharp interface limit of this model with respect to the obstacle approximation using the concept of matched asymptotics to show that, at least formally, the model converges towards two-phase flow in the implicitly given fluid domain if the phase field parameter  of the obstacle tends to zero.
    \item We provide a well-posed time discretization scheme for the above model.
\end{enumerate}
The analytical results are validated by simulations.

\subsubsection*{Motivation}
Our work is motivated by fluid mechanic applications
such as pinning of sliding droplets \cite{BJKR2019CET}, design of ship hulls \cite{MSR2023} and controlling Rayleigh--Taylor instability \cite{YK2023FEAD}, where shape and topology optimization of obstacles would be useful. 
This current work, forms the basis for the mathematical formulation of the underlying physical systems in the context of two-phase flows,  modeled using the phase field approach, with any obstacles that may be present and require optimization.

Shape and topology optimization naturally
accommodates complex topological transitions. 
Representing the obstacles via a phase-field variable  facilitates a mathematically tractable framework for shape and topology optimization, particularly for two-phase flows governed by the Cahn--Hilliard/Navier--Stokes equations. 
In this approach, the sought topology is not explicitly represented and thus topology changes are naturally possible by the model.
Furthermore, this quantification allows the optimization problem to be cast as an optimal control problem constrained by partial differential equations (PDEs), enabling the use of the well developed machinery from PDE-constrained optimization \cite{HPUU2009}.

Accordingly, the current work focuses on the modeling and numerical simulation of two phase flow in the presence of obstacles, serving as a necessary precursor to the development of robust optimization algorithms for systems governed by two-phase flow. By ensuring a stable and consistent forward model, we bridge the gap between abstract theory and the practical numerics required for shape and topology design in physical regimes governed by two-phase flows.

\subsubsection*{Related literature}

The study of two-phase flows is an active field of research due to its relevance to applicable sciences such as engineering, biology and physics, to name a few. To model this phenomena, there are several noteworthy works in the literature such as the work \cite[Section 1.6.7 and Chapter 8]{G1994}.
Gidaspow \cite{G1994} proposes a model that tracks the individual velocity of each fluid phase and considers a transport equation that follows through the concentration of each fluid, also known as the volume-of-fluid (VoF) approach. 
Another standard model for two phase flow of two liquids that are separated by a sharp interface is e.g. stated in \cite{Gross2011}, where numerical methods are also discussed. Here the fluid domain is separated into the two subdomains with the two phases and these two fluids are modeled by Navier--Stokes equation. At the interface between the two phases certain jump conditions define the motion of the interface.
A numerical treatment thus especially needs to track this interface. Several approaches are suggested for this, where we especially mention the tracking by a level-set function as in \cite{Gross2011,}. However a monolithic numerical treatment of the complete two-phase fluid is complicated for this model.

The Cahn--Hilliard equation \cite{CH1959TJCP} models the distribution of two fluids, whose interaction is diffusion dominated. The two phases are described by a phase-field variable, that serves as a smooth indicator function for the binary fluids. An extension to two phase flow is given in \cite{HohenbergHalperin} for equal density flow.
To accommodate the interaction of a binary fluid with different densities, \cite{AGG2012-3MAS} proposed a model that maintains thermodynamic consistency, i.e., the model satisfies local and global dissipation inequalities. This work has been the basis of several works that aim to model, approximately, fluid-to-gas dynamics such as gas bubble dynamics and sliding droplet, see the works of \cite{HTMP2017JCP} and \cite{BJKR2019CET} for illustration. 
Alternative models for two phase flow that describe the distribution of the involved fluids by a phase field are e.g. \cite{Boyer,DingSpeltShu,GPV1996M3AS,LowengrubTruskinowsky}.


Numerically solving the model proposed in \cite{AGG2012-3MAS} can be quite challenging not only due to its nonlinearity of the individual equations, but also from the fact that the equations describing the fluid velocity and the phase distribution are nonlinearly coupled. To circumvent this, \cite{M2012-NMPDE},  \cite{GGT2014-JCM} and  \cite{GGGM2016-CCP} considered splitting schemes that decouples the equations. This is done by plucking out an intermediate velocity from the momentum equation that absorbs the capillary forces. This decoupling method has been since utilized to model flow of dilute polymeric solutions \cite{M2018ESAIM}, non-Newtonian fluids with varying polymer density \cite{S2020JCM},  solid-to-liquid dynamics \cite{RvWGR2026arXiv} and sliding droplets with consideration of contact angle \cite{BK2021-SICON} to name a few.

The sequential coupling of equation for the phase distribution and the velocity field leaves the nonlinearity associated with $W$ within the system, while preserving thermodynamically consistency in the time discrete setting.
The often used convex-concave splitting of $W$ \cite{eyre_CH_semi_implicite,Gruen_Klingbeil_CHNS_AGG_numeric,GHK2016-LinearScheme,RvWGR2026arXiv}
leads to a nonlinear thermodynamically consistent scheme.
As alternative, 
Invariant Energy Quadratization (IEQ) schemes \cite{YANG20171005,YangZhang-Convergence-IEQ} and 
Scalar Auxiliary Variable  (SAV) schemes \cite{SHEN2018407} are proposed that introduce additional equations to the system to fulfill a (modified) energy inequality.


A central question is, how the equations are extended from the fluid domain to the full domain, while inside the domain of interest and at its boundary proper physics are satisfied.
The SIMP approach \cite{Bendsoe1989-SIMP,Bendsoe2004-TopoOpt,Evgrafov2005} 
is a very common approach to extend equations that satisfy Neumann boundary data to a larger domain.
Originating in structural optimization, it assumes an artificial material with low density
to extend the equation of interest to a larger computational domain. 
Mathematically speaking, this leads to a spatially varying coefficient in the second order operator.
For homogeneous Dirichlet boundary data similarly the porous medium approach 
\cite{BorPIJNMF2003} 
assumes an artificial porous medium outside the domain of interest.
This leads to  an additional lower order Darcy term in the equation.

Closely related to these approaches is the diffuse domain approach proposed in \cite{LiLowengrubRaetzVoigt-DiffuseDomainApproach}.
It especially allows to extend an equation from a lower dimensional manifold to a surrounding domain as e.g. in \cite{TeigenSongLowengrubVoigt-Surfactants,GarckeLamStinner-Surfactants,AbelsLamStinner}.

Our formulation for two-phase flow is an  extension of the formulation for one-phase flow that is used for topology optimization in \cite{GarckeHechtHinzeKahle-2015-SIAM,GarckeHecht-TopoOpt,GHHKL2016-IFB,GHKL2018-ACM}.
Extending these works will be investigated in a subsequent paper.

\subsubsection*{Outline}

Section~\ref{sec:model} introduces the two-phase flow model in detail, including the modifications arising from the porous-medium approximation of the obstacle.  A formal matched asymptotic analysis is carried out in Section~\ref{sec:asymptotics}, establishing the consistency of the obstacle approximation. 
Section~\ref{sec:discrete} is devoted to the development of a practical time discretization suitable for numerical implementation, where we provide existence of solutions to the discrete system along with global stability results. Finally, Section~\ref{sec:numerics} presents numerical experiments that validate the proposed obstacle approximation and illustrate the asymptotic behavior predicted by the analysis in Section~\ref{sec:asymptotics}.

\subsection*{Notations and functional analytic tools}

We denote by $L^p(\Omega)$ the space of Lebesgue integrable function of order $p\ge 1$ and the usual Sobolev spaces are denoted as $W^{m,p}(\Omega)$ with $H^m(\Omega) = W^{m,2}(\Omega)$ for $m\in \mathbb{N}$. 
For simplicity, we use $\|\cdot\|$ and $(\cdot,\cdot)_\Omega$ to denote the norm and inner product in $L^2(\Omega)$ or $L^2(\Omega)^d$, for $d=2,3$, $\|\cdot\|_{L^p}$ for the norm in $L^p(\Omega)$ or $L^p(\Omega)^d$, and $\|\cdot\|_{W^{m,p}}$ the norm in $W^{m,p}(\Omega)$ or $W^{m,p}(\Omega)^d$.

To take into account the incompressibility of the fluid we utilize the following solenoidal spaces $H_\sigma = \{{\bu}\in L^2(\Omega)^d : \dive{\bu} = 0\text{ in }\Omega,\, {\bu}\cdot{\bn} = 0\text{ on }\partial\Omega \}$, and $V_\sigma = \{{\bu}\in H^1_0(\Omega)^d: \dive{\bu} = 0 \text{ in } \Omega \}$.
%


Furthermore we define $H = L^2(\Omega)$ and $V = H^1(\Omega)$.

Finally, for $q > 2$ (if $d=2$) or $q\geq d$ (if $d=3$), we introduce the trilinear form 
$a: (L^q(\Omega))^d\times V_\sigma\times V_\sigma\to \mathbb{R}$ defined as 
\begin{align*}
    a(\bv,\bu,\bw) = \frac{1}{2}\left\{ \int_\Omega [(\bv\cdot\nabla)\bu]\cdot \bw \du x - \int_\Omega[(\bv\cdot\nabla)\bw]\cdot \bu \du x \right\}.
\end{align*}
For $\bv\in (L^q(\Omega))^d$ and $\bu \in V_\sigma$ by construction it holds 
$a(\bv,\bu,\bu) = 0$, i,e. $a$ is antisymmetric with respect to the second and third argument.

\section{Modeling}
\label{sec:model}
Let us consider a hold-all Lipschitz domain $\Omega\subset \mathbb{R}^d$ with $d\in\{2,3 \}$ with outer unit normal $\nu_\Omega$. Consider a subdomain $B\subset \Omega$, which we shall refer to as an obstacle. We suppose that the domain $E := \Omega\setminus \overline{B}$ is occupied by a two-phase fluid with varying density and viscosity.

We start with the description of the two-phase model inside the fluid domain $E$ 
and thereafter extend this model to the hold-all domain $\Omega$ including the obstacle.

\subsection{Two-phase flow inside $E$}
 To model the dynamics inside $E$ we consider a thermodynamically consistent phase-field model proposed in \cite{AGG2012-3MAS}, with the additional assumption that the density of the fluid may be non-linearly dependent on the phase-field \cite{AbelsBreit-NonlinearDensity}. 
 
Let $T>0$ and $I:= (0,T]$ be a time interval. 
Within $E$ there is a two-phase fluid that consists of  two immiscible fluids.
The model assumes, that there is a small domain with width proportional to $\varepsilon$ between the two fluids, where both fluids coexists.
With this assumption, the distribution of the two fluids is described by  the phase field $\varphi: I\times E \to \mathbb{R}$. It serves as a smooth indicator function such that $\varphi(x)=1$ indicates the first phase and $\varphi(x)=-1$ indicates the second phase, while $|\varphi(x)|<1$ indicates the domain where both fluids coexist -- the so-called diffuse interface.
Associated to $\varphi$ is a chemical potential $\mu : I\times E \to \mathbb{R}$. 
Moreover, the fluid velocity and pressure of the two-phase fluid are given by  $\bu : I\times E\to \mathbb{R}^d$ and  $p:I\times E\to \mathbb{R}$. 
Such quantities satisfy the following system
    \begin{subequations}
    \label{system:prime}
        \begin{align}
            \rho\partial_t\bu + (( \rho\bu + \bJ )\cdot\nabla)\bu + {R\frac{\bu}{2}} - \dive(2\eta \D\bu) + \nabla p & = -\varphi\nabla\mu + \rho\bg&& \text{ in }I\times E, \label{eq:sysE:NS1}\\
            \dive\bu & = 0&& \text{ in }I\times E,\label{eq:sysE:NS2}\\
              \partial_t\varphi + \dive(\varphi\bu - b\nabla \mu) &= 0 &&\text{ in }I\times E, \label{eq:sysE:CH1}\\
             - \sigma{\varepsilon}\Delta\varphi + \frac{\sigma}{{\varepsilon}}W'(\varphi) &= \mu &&\text{ in }I\times E, \label{eq:sysE:CH2}
        \end{align}
    \end{subequations}
    with boundary conditions $\bu = 0$, $\partial_\nu \varphi = \partial_\nu \mu = 0$ on $I\times \partial E$, where $\partial_\nu$ denotes the derivative in outer normal direction.
    Furthermore, the initial conditions are $\bu(0) = \bu^0$ and $\varphi(0) = \varphi^0$ for given sufficiently smooth data $\bu^0$ and $\varphi^0$. 
    Here, $\rho \equiv \rho(\varphi)$ and $\eta \equiv \eta(\varphi)$ are the spatially varying density and viscosity of the binary fluid which are both defined via the phase-field $\varphi$.
    Additionally, we abbreviate $\bJ = -b\rho^\prime\nabla\mu$ 
    and $R = -b\nabla\rho^\prime\cdot\nabla\mu$.
    The function $W$ is called the free energy potential that attains exactly two minima at $-1$ and $1$, i.e. at the values that indicate the two pure fluids.
    The constant $b>0$ is called the mobility, the constant ${\varepsilon}>0$  is proportional to the diffuse interface between the two fluids. The constant $\sigma = c_W\sigma_s>0$ includes the surface tension at the interface, where $\sigma_s$ denotes the physical surface tension and $c_W >0$ is a constant that appears from the phase field modeling of the system. It is defined as 
    \begin{align*}
        \frac{1}{c_W} = \int_{-\infty}^{\infty} 2 W(\widetilde{\Phi}_0(s))\du s,
    \end{align*}
     where $\widetilde{\Phi}_0$ satisfies $\partial_{zz} \widetilde{\Phi}_0(z) = W'(\widetilde{\Phi}_0(z))$. 
     Such objects are derived from formal matched asymptotics in \cite{AGG2012-3MAS} and especially, $\widetilde{\Phi}_0$ represents a first order approximation of $\varphi$.
    Further specifications of the free energy potential $W$, density $\rho$ and viscosity $\eta$ 
    are discussed in Section~\ref{sec:discrete}.

\bigskip

We use  a reformulation of the momentum equation \eqref{eq:sysE:NS1}, which has been proven convenient for defining the weak formulation and obtaining an energy dissipation equality. 
To begin, we multiply $\rho^\prime$ to \eqref{eq:sysE:CH1}, and utilize $\dive{\bu}=0$ which leads to
\begin{align}
    \partial_t\rho + \dive( \rho \bu + \bJ ) = R.     \label{rhoeq}
\end{align}    
Consequently, multiplying $\frac{\bu}{2}$ to \eqref{rhoeq} and adding the result to \eqref{eq:sysE:NS1} yield
\begin{align}
        \frac{1}{2}\big[\rho\partial_t \bu + \partial_t(\rho\bu) + ((\rho\bu + \bJ)\cdot\nabla)\bu + \dive((\rho\bu + \bJ) \otimes \bu) \big]
        - \dive(2 \eta \D\bu) + \nabla p = -\varphi\nabla\mu + \rho\bg.        
    \label{unifNS-rhoeq}
\end{align}
We stress, that this reformulation is independent of the actual boundary condition of the involved functions. It will  be useful for the numerical analysis of the time-stepping scheme, that we derive in Section~\ref{sec:discrete}.

\subsection{Phase-field approximation of the obstacle $B$}
To incorporate the obstacle into the system, we follow both \cite{GHKKL2023-TopoOpt-Laplace} and the porous medium approach from \cite{BorPIJNMF2003,GarckeHecht-TopoOpt,GHHKL2016-IFB}.

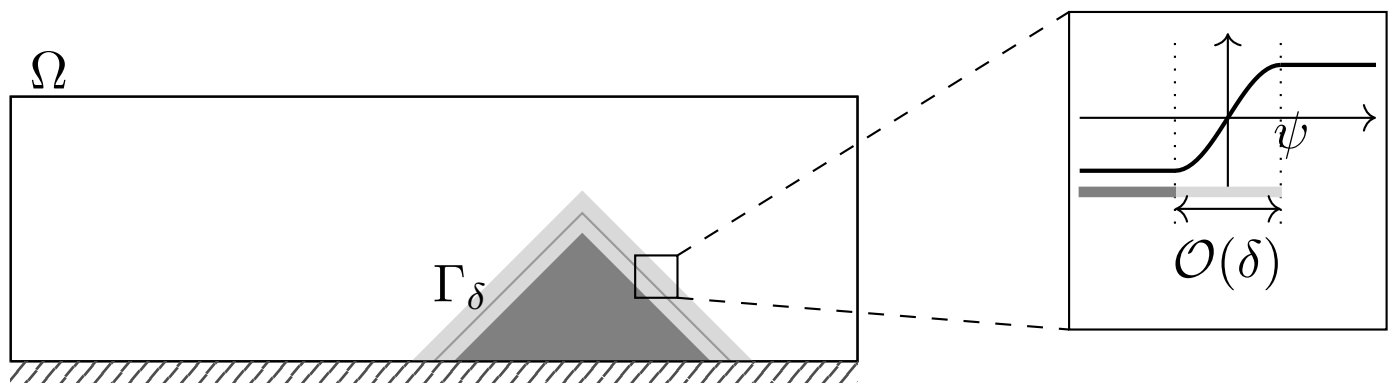
\begin{figure}
\centering
{
\tikzset{every picture/.style={scale=1.4}}
\begin{tikzpicture}

\colorlet{clFluidA}{black!15}
\colorlet{clFluidI}{black!30}
\colorlet{clFluidB}{black!50}
















\coordinate (BL) at (0,0);
\coordinate (TR)  at (4,1.25);
\draw[fill=white] (BL)  rectangle (TR);

\begin{scope}[shift={(2.0,0.0)}] 
    \coordinate (A) at (0,0);
    \coordinate (B) at (0.7,0.7);
    \coordinate (C) at (1.4,0);
    \pgfmathsetmacro{\r}{0.1}

\draw[color=black!15,fill=black!15] ($(A) - (\r,0)$) -- ($(B) + (0,\r)$) -- ($(C) + (\r,0)$) -- cycle; 
\draw[color=black!40] (A)--(B)--(C)--cycle;
\draw[color=black!50,fill=black!50] ($(A) + (\r,0)$) -- ($(B) - (0,\r)$) -- ($(C) - (\r,0)$) -- cycle;
 

 \node[anchor=east] at ($0.5*(B) + 0.5*(A)$) {$\Gamma_\delta$};
 
\end{scope}


\node[anchor=south west] at (0.0,1.2)  {$\Omega$};

\draw (BL)  rectangle (TR);

\fill[pattern={Lines[angle=50, distance=2.5pt,line width=0.5pt]}, pattern color=black!70] (0.0,-0.1) rectangle (4.0,0.0);

  \begin{scope}[shift={(5.0,0.4)}]

    
    \draw (-2.05,-0.1) rectangle ++(0.20,0.20);

    \draw[dashed] (-1.85,-0.1) -- (0,-0.25);
    \draw[dashed] (-1.85,0.1) -- (0,1.25);

    \draw (0,-0.25) rectangle(1.5,1.25);

    \draw[->] (0.05,0.75) -- (1.45,0.75); 
    \draw[->] (0.75,0.40) -- (0.75,1.15); 

    \draw[dotted] (0.5,0.25) -- (0.5,1.15);
    \draw[dotted] (1.0,0.25) -- (1.0,1.15);


    \colorlet{clphi}{black!100}
    \draw[thick,color=clphi,domain=-1:1] plot (0.75 - \x/4,{0.75-0.25*(sin(pi/2*\x r))});
    \draw[thick,color=clphi,solid] (0.05,0.5) -- (0.5,0.5);
    \draw[thick,color=clphi,solid] (1.0,1.0) -- (1.45,1.0);

    \draw (0.875,0.87) node[anchor=north west]{\color{clphi}$\psi$};

    \draw[fill,color=black!00] (1.0,0.38) rectangle ++(0.45,0.04);
    \draw[fill,color=black!15] (0.5,0.38) rectangle ++(0.5,0.04);
    \draw[fill,color=black!50](0.05,0.38) rectangle ++(0.45,0.04);
    \draw[<->] (0.5,0.32)--(1.0,0.32);
    \draw (0.75,0.275) node[anchor=north] {$\mathcal O (\delta)$};

  \end{scope}

\end{tikzpicture}
}
\caption{
Sketch of the phase-field $\psi$, that describes the diffuse obstacle. We show the phase-field function $\psi$ in gray scale and plot the function $\psi$ across the  diffuse interface. Here dark gray indicates $\psi=-1$, i.e. the obstacle, and white indicates $\psi=+1$, i.e. the fluid domain.
}
\label{fig:intro:schematicobstaclesetup}
\end{figure}

We  assume that the interface between obstacle and flow domain is diffuse and of width proportional to $\delta>0$. 
To describe the location of the obstacle we consider a smooth
stationary 
indicator function $\psi : \Omega \to \mathbb{R}$ in the form of a phase field, i.e.  $\psi(x) = -1$ indicates the obstacle, while $\psi(x) = 1$ indicates the fluid domain. 
Again $|\psi|<1$ indicates a  diffuse region between the two phases, i.e. between fluid domain and obstacle. 
In Figure~\ref{fig:intro:schematicobstaclesetup} we sketch this situation.
We show the obstacle in gray scale and plot the function $\psi$ across the  diffuse interface.
To incorporate such description into the governing system,
we introduce interpolation functions 
\begin{itemize}
        \item $\alpha_\bu^{(\delta)} : \mathbb{R}\to \mathbb{R}$, 
        such that $\alphu{1} = 0 $ and  $\alphu{-1} = \delta^{-1}\overline{\alpha}_\bu$, 
        where  $\overline{\alpha}_\bu > 0$,
        \item $\alpha_\varphi^{(\delta)} : \mathbb{R}\to \mathbb{R}$ satisfying $\alphi{1} = 1 $ and  $\alphi{-1} = \delta\underline{\alpha}_\varphi$, where  $\underline{\alpha}_\varphi>0$.
\end{itemize}
Especially, we consider the obstacle as porous medium and  extend \eqref{system:prime} and \eqref{unifNS-rhoeq} from $E$ to $\Omega$ as
\begin{subequations}
\label{system:pores}
    \begin{align}
        \frac{1}{2}\big[\rho\partial_t \bu + \partial_t(\rho\bu) + ((\rho\bu + \bJ_\delta)\cdot\nabla)\bu + \dive((\rho\bu + \bJ_\delta) \otimes \bu) \big] &\notag\\
        + \alphu{\psi}\bu- \dive(2 \eta \D\bu) + \nabla p &= -\varphi\nabla\mu + \rho\bg &&\text{ in }I\times \Omega, \label{NS}\\
    \dive \bu & = 0 && \text{ in }I\times \Omega, \label{divfree}\\
    \partial_t\varphi + \dive(\varphi\bu - b\alphi{\psi} \nabla\mu ) & = 0 &&\text{ in }I\times \Omega,  \label{phieq}\\
         - \sigma{\varepsilon}\dive(\alphi{\psi}\nabla \varphi) + \frac{\sigma}{{\varepsilon}}W'(\varphi) & = \mu &&\text{ in }I\times \Omega,\label{mueq}
    \end{align}
\end{subequations}
where $\bJ_\delta = -b\alphi{\psi}\rho^\prime\nabla\mu$,
with boundary conditions $\bu = 0$, $\partial_\nu \varphi = \partial_\nu \mu = 0$ on $I\times \partial \Omega$,  and initial conditions $\bu(0) = \bu^0$ and $\varphi(0) = \varphi^0$ in $\Omega$. 
We notice that the boundary conditions of the original system \eqref{system:prime} on the boundaries of the subdomain $B$ are dropped and are also relaxed due the inclusion of the  interpolants $\alpha_\varphi^{(\delta)}$ and $\alpha_{\bu}^{(\delta)}$. 
Nevertheless, the boundary data on $B$ are given implicitly, and can be recovered, as we show in the formal matched asymptotics in Section~\ref{sec:asymptotics}.

\begin{remark}
We note that, while both $\varphi$ and $\psi$ are phase-fields, the phase-field $\varphi$ is obtained by solving the instationary fluid equations \eqref{system:pores} and represents the distribution of the two fluids, while $\psi$ is a chosen data and represents the stationary distribution of obstacle and fluid domain.
Especially, we choose $\psi$ to precisely satisfy $\psi = \pm1$  inside the fluid domain and obstacle domain and consequently $\nabla \psi = 0$, while for $\varphi$  inside the pure fluid phases $\varphi=\pm1$ generally  holds only with equality in the limit $\epsilon\to 0$.
\end{remark}

\begin{remark}
    We note that we also modify the flux $\bJ$ by including $\alphi{\psi}$ to achieve $\bJ_\delta$. This is for consistency reasons. In this way, we could equivalently first extend \eqref{system:prime} to $\Omega$ and then perform the reformulation that leads to \eqref{unifNS-rhoeq}.
\end{remark}

   \begin{remark}
   \label{rm:conservation-phiavg}
        We note that by integrating \eqref{phieq} with respect to space and using integration by parts and the boundary conditions we recover conservation for the phase field $\varphi$ in $\Omega$, i.e.,
        \begin{align*}
            \frac{d}{dt}\int_{\Omega} \varphi \du x = 0.
        \end{align*}
        On the other hand, multiplying \eqref{phieq} by $\frac{\psi + 1}{2}$, which can be interpreted as indicator function of the fluid domain, and integrating the resulting product gives us, after integration by parts and using the boundary conditions,
        \begin{align}
            \label{mass}
            \begin{aligned}
                \frac{d}{dt}\int_{\Omega}\frac{\psi + 1}{2}\varphi \du x = \frac{1}{2}(\varphi\bu, \nabla \psi)_\Omega - \frac{b}{2}( \alphi{\psi}\nabla\mu, \nabla \psi )_\Omega\\ 
            = \frac{1}{2}\int_{\Gamma_\delta}\varphi\bu\cdot \nabla \psi \du x 
            - \frac{b}{2}\int_{\Gamma_\delta} \alphi{\psi}\nabla\mu\cdot \nabla \psi \du x
            \end{aligned}
        \end{align}
        where $\Gamma_\delta := \{x\in\Omega: -1 < \psi(x) < 1   \}$ and we use that $\nabla \psi = 0$ if $\psi=\pm1$.
        
        As we show in Section~\ref{sec:conservation-phiavg} by formal matched asymptotics the integrals involving $\Gamma_\delta$ formally 
        go to zero as $\delta \to 0$.
        This result is consistent with the fact that as $\delta\to 0$ the integral on the left-hand side goes to the integral in $E$, which is the expected conservation for the phase field in system \eqref{system:prime}.
    \end{remark}


\section{Sharp interface asymptotics for the obstacle}
\label{sec:asymptotics}
 As a way of validation, we conduct a formal asymptotic analysis to establish that \eqref{system:prime} is the asymptotic limit of \eqref{system:pores} as the parameter $\delta>0$ goes to zero. For this purpose let us recall the governing equation in $(0,T]\times \Omega$, with highlight on the dependence on $\delta$ of the involved variables 
    \begin{subequations}
    \label{system:delpores}
        \begin{align}
                \frac{1}{2}\big[\rho_\delta\partial_t \bu_\delta + \partial_t(\rho_\delta\bu_\delta) + ((\rho_\delta\bu_\delta + \bJ_\delta)\cdot\nabla)\bu_\delta + \dive((\rho_\delta\bu_\delta + \bJ_\delta) \otimes \bu_\delta) \big] &\notag \\
                + \alphu{\psi}\bu_\delta - \dive(2 \eta_\delta\D\bu_\delta) + \nabla p_\delta 
                &= -\varphi_\delta\nabla\mu_\delta + \rho_\delta\bg, \label{delNS}\\
            -\dive{\bu_\delta} &= 0,\label{delsolNS}\\
             \partial_t\varphi_\delta + \dive(\varphi_\delta\bu_\delta - b\alphi{\psi_\delta} \nabla\mu_\delta ) & = 0, \label{delphieq}\\
             - \sigma{\varepsilon}\dive(\alphi{\psi_\delta}\nabla \varphi_\delta) + \frac{\sigma}{{\varepsilon}}W'(\varphi_\delta) &= \mu_\delta,\label{delmueq}
        \end{align}
    \end{subequations}
    where $\rho_\delta := \rho(\varphi_\delta)$ and $\eta_\delta := \eta(\varphi_\delta)$, together with the initial conditions $\bu_\delta(0) = \bu^0$ and $\varphi_\delta(0) = \varphi^0$ and the boundary conditions
    \begin{align*}
        \bu_\delta = 0,\quad \partial_\nu \varphi_\delta = \partial_\nu\mu_\delta = 0  \text{ on }\partial\Omega.
    \end{align*}

Exclusively for the sake of the asymptotic investigation we state the following assumptions on system~\eqref{system:delpores}.
\begin{assumption}
\mbox{}
     \begin{itemize}
        \item We restrict to $d=2$, i.e., the two-dimensional setting.  The analysis for three dimensions follows the same procedure, but is more technically complex.
         \item  We assume,  that system~\eqref{system:delpores} admits unique sufficiently smooth solution $(\bu_\delta,p_\delta,\varphi_\delta,\mu_\delta)$.
         \item 
   We assume sufficient regularity for the free-energy potential $W$, the density $\rho$ and the viscosity $\eta$ to allow Taylor expansion.
   \item The interpolants $\alpha_{\varphi}^{(\delta)}$ and $\alpha_{\bu}^{(\delta)}$ take the following forms:
    \begin{align*}
        \alphi{s} = \frac{1+s}{2} + \delta\underline{\alpha}_\varphi\frac{1-s}{2}, 
        \qquad \alphu{s} = \frac{\overline{\alpha}_\bu}{\delta}\frac{1-s}{2},
    \end{align*}
    where $\overline{\alpha}_\varphi,\underline{\alpha}_\bu > 0$. 
     \end{itemize}  
\end{assumption}


\paragraph{Expansions of the phase-field $\psi_\delta$} \mbox{}\\
The phase-field $\psi_\delta$ describes the distribution of the bulk regions inside $\Omega$, i.e. fluid domain and obstacle, by taking the values $\pm 1$.
We use the notation for the bulk regions $E_\delta:= \{x\in\Omega: \psi_\delta(x) = 1 \}$ and 
$B_\delta:= \{x\in\Omega: \psi_\delta(x) = -1 \}$,
and denote by $\Gamma_\delta := \{x\in\Omega: -1 < \psi_\delta(x) <1  \}$ the interface. Furthermore, we impose the following assumptions for the phase-field $\psi_\delta$ and its expansions:
\begin{assumption}\mbox{}
\label{expansion:assumption}
    \begin{enumerate}[label = (\Roman*)]
        \item For $\delta \to 0$, the phase-field $\psi_\delta$ converges towards a function $\psi_0\in\{-1,1\}$ a.e. in $\Omega$. Consequently, we can expand $\psi_\delta$ as $\psi_\delta = \psi_0 + \delta\psi_1 + \ldots$ in $\Omega$. 
        Furthermore, it holds $E = \{ x\in\Omega : \psi_0(x) = 1 \}$ and $B = \{x\in\Omega: \psi_0(x)=-1  \}$ and the sets $\{x \in \Omega: \psi_\delta(x) = 0 \}$ and $\Gamma = \partial E\cap\partial B$ coincide.
        \item In the bulk regions $E_\delta$ and $B_\delta$ we have $\psi_1 \equiv \psi_2 \equiv \ldots \equiv 0$.
        \item After the introduction of new coordinate system $(s,z)$, see \eqref{coordchange}, we have $\psi_\delta(x) = \Psi_\delta(s,z)$ with the expansion $\Psi_\delta = \Psi_0 + \delta \Psi_1 + \ldots$, where $\Psi_0$ is independent of the variable $s$ with $\Psi_0 (s,0) = 0$ and $ \lim_{z\to \pm \infty}\Psi_0 = \pm 1$. Furthermore, in $\Gamma_\delta$ we have $\Psi_1 \equiv \Psi_2 \equiv \ldots \equiv 0$.\label{interior:assumption}
    \end{enumerate}
\end{assumption}



\subsection{Outer expansions}
We assume that a variable of interest $y_\delta \in \{\bu_\delta, p_\delta, \varphi_\delta, \mu_\delta \}$ can be expanded with respect to $\delta$ as 
\begin{align}\label{genexpansion}
    y_\delta = y_0 + \delta y_1 + \ldots,
\end{align}
with certain functions $y_i$, $i=0,1,\ldots$.

Moreover, by using Taylor's theorem we obtain the following expansions with respect to $\delta$ for the fluid density,  viscosity and the potential $W$
\begin{align*}
    \rho_\delta = \rho(\varphi_\delta) &= \rho(\varphi_0) + \delta \rho'(\varphi_0)\varphi_1 + \ldots =: \rho_0 + \delta \rho_1 + \ldots,\\
    \eta_\delta = \eta(\varphi_\delta) &= \eta(\varphi_0) + \delta \eta'(\varphi_0)\varphi_1 + \ldots =: \eta_0 + \delta \eta_1 + \ldots,\\
    W'(\varphi_\delta)&= W'(\varphi_0) + \delta W''(\varphi_0)\varphi_1 + \ldots\, .
\end{align*}
%

\paragraph{Expansion at order $\mathcal O(\delta^{-1})$}\mbox{}\\
    Plugging these expansions into \eqref{system:delpores}, the leading order $\mathcal{O}(\delta^{-1})$ in \eqref{delNS} gives us
    \begin{align*}
        \frac{\overline{\alpha}_\bu}{2}(1-\psi_0)\bu_0 = 0
    \end{align*}
    which implies 
    \begin{align}\label{UinB}
        \bu_0 = 0\text{ whenever }\psi_0\ne 1.
    \end{align}


\paragraph{Expansion at order $\mathcal O(\delta^0)$}\mbox{}\\
    Taking the $\mathcal{O}(\delta^0)$-order terms in \eqref{delNS} leads us to
    \begin{align*}
				&\begin{aligned}
					&\frac{1}{2}\big(\rho_0\partial_t \bu_0 + \partial_t(\rho_0\bu_0)\big)  + \frac{1}{2}\left[\left(\rho_0\bu_0 - b\frac{\psi_0+1}{2}\rho'(\varphi_0)\nabla \mu_0 \right)\cdot \nabla \right]\bu_0\\
                    & + \frac{1}{2} \dive \left( \left(\rho_0\bu_0 
                    - b\frac{\psi_0+1}{2}\rho'(\varphi_0)\nabla \mu_0 \right) \otimes \bu_0 \right) 
                    - \frac{\overline{\alpha}_\bu}{2}( \psi_1\bu_0 + (\psi_0-1)\bu_1 )  \\
                    & - 2\dive( \eta_0 \D\bu_0)  
					 +  \varphi_0\nabla\mu_0 + \nabla p_0 = \rho_0 \bg.
				\end{aligned}
    \end{align*}
    In the bulk region $E_\delta$, where $\psi_0 = 1$ and $\psi_1=0$, we get
    \begin{align*}			
        \frac{1}{2}\big(\rho_0\partial_t \bu_0 + \partial_t(\rho_0\bu_0)\big)  + \frac{1}{2}\left[\left(\left(\rho_0\bu_0 + \bJ_0   \right)\cdot \nabla \right)\bu_0 + \dive \left( \left(\rho_0\bu_0 + \bJ_0  \right) \otimes \bu_0 \right) \right]&\\
         - 2\dive( \eta_0 \D\bu_0)  
         + \nabla p_0 &= -\varphi_0\nabla\mu_0  + \rho_0 \bg.
    \end{align*}
    where $\bJ_0 = -b\rho'(\varphi_0)\nabla \mu_0$. This is the reformulated momentum equation \eqref{unifNS-rhoeq} from system~\ref{system:prime}.
    
    The incompressibility condition $\dive \bu_0 = 0$ from \eqref{eq:sysE:NS2} is  recovered by taking the expansion of $\bu_\delta$ and its divergence.

    The terms of order $\mathcal{O}(\delta^0)$ in \eqref{delphieq} and \eqref{delmueq} give us
    \begin{align*}
         \partial_t\varphi_0  + \dive\left(\varphi_0\bu_0 - b\frac{\psi_0+1}{2} \nabla\mu_0 \right)  &= 0,\\
         -\sigma\varepsilon\dive\left(\frac{\psi_0+1}{2} \nabla\varphi_0  \right) + \frac{\sigma}{\varepsilon} W'(\varphi_0) &= \mu_0.
    \end{align*}
    In the bulk region $E_\delta$, where $\psi_0=1$, the equations above  reduce to
    \begin{align*}
         \partial_t\varphi_0  + \dive (\varphi_0\bu_0-b \nabla\mu_0 )  &= 0,\\
        -\sigma\varepsilon\Delta\varphi_0 + \frac{\sigma}{\varepsilon} W'(\varphi_0) &= \mu_0,
    \end{align*}
    and thus to \eqref{eq:sysE:CH1}--\eqref{eq:sysE:CH2}.
    
    On the other hand, in the bulk region $B_\delta$,  where $\psi_0=-1$, these equations reduce to
    \begin{align}
         \partial_t \varphi_0&=0, \label{eq:asy:phi-leading-B}\\
         \mu_0 &= \frac{\sigma}{\varepsilon} W'(\varphi_0).
    \end{align}
    We also see that the boundary conditions 
    \begin{align*}
        \bu_0 = 0,\quad \partial_\nu \varphi_0 = \partial_\nu\mu_0 = 0  \text{ on }\partial\Omega
     \end{align*}
    are satisfied by considering the leading order of the expansions at the boundary.

\subsection{Inner expansions}

In the interface $\Gamma_\delta$ we introduce a parametrization in a neighborhood of the boundary $\Gamma$. We consider the signed distance function $d(x)$ to $\Gamma$, which satisfies $d(x) > 0 $ if $x\in E$ and $d(x) < 0 $ if $x \in B$, and set $z = \frac{d}{\delta}$ as the $\frac{1}{\delta}$-scaled distance function from $\Gamma$. We consider the parametrization $\gamma(s)$ of $\Gamma$ by arc-length $s$, and mark by $\nu$ the unit normal to $\Gamma$ pointing from $B$ into $E$.

For a sufficiently smooth function $y$, we can thus write for $x\in \Gamma(\delta)$
\begin{align}\label{coordchange}
    y(x) = y(\gamma(s) + \delta z\nu(\gamma(s)) ) =: Y(s,z).
\end{align}

The introduction of the new $(s,z)$-coordinates gives us the following formula for the derivatives \cite[Section 5.2]{GHHKL2016-IFB}: 
\begin{itemize}[label=$\bullet$]
    \item $\nabla y = \displaystyle \frac{1}{\delta}(\partial_z Y)\nu + \nabla_\Gamma Y + \ldots $;
    \item $\dive \boldsymbol{y} = \displaystyle \frac{1}{\delta}\partial_z\boldsymbol{Y}\cdot\nu + \dive_\Gamma \boldsymbol{Y} + \ldots$;
    \item $\D\boldsymbol{y} = \displaystyle \frac{1}{2\delta}( \partial_z \boldsymbol{Y} \otimes \nu + (\partial_z \boldsymbol{Y} \otimes \nu)^T ) + D_\Gamma \boldsymbol{Y} + \ldots$;
    \item $\dive\boldsymbol{a} = \displaystyle \frac{1}{\delta}(\partial_z \boldsymbol{A}) \nu + \dive_\Gamma \boldsymbol{A} + \ldots$;
    \item $\Delta y = \displaystyle \frac{1}{\delta^2} \partial_{zz}Y + \frac{1}{\delta}\dive_\Gamma((\partial_z Y)\nu) + \ldots =  \frac{1}{\delta^2} \partial_{zz}Y - \frac{1}{\delta}\kappa \partial_zY + \ldots$ ,
\end{itemize}
where $\boldsymbol{y}(x) = \boldsymbol{Y}(s,z)$ is a vector-valued function, $\boldsymbol{a}(x) = \boldsymbol{A}(s,z)$ is a tensor-valued function, $\nabla_\Gamma$ denotes the surface gradient on $\Gamma$, $D_\Gamma \boldsymbol{y} = \frac{1}{2}(\nabla_\Gamma \boldsymbol{y} + (\nabla_\Gamma \boldsymbol{y})^T)$,  $\dive_\Gamma$ denotes the surface divergence, and $\kappa$ denotes the mean curvature of $\Gamma$ defined as $\partial_s \nu = -\kappa\tau$, where $\tau$ is the unit tangent vector such that $(\nu,\tau)$ is positively oriented.

For the variables of interest $\bu_\delta, p_\delta, \varphi_\delta$ and $\mu_\delta$ we use $\bU_\delta, P_\delta, \Phi_\delta$ and $M_\delta$ in the new coordinate system. We assume that for $Y_\delta\in \{ \bU_\delta, P_\delta, \Phi_\delta, M_\delta \}$ the following expansion holds
\begin{align*}
    Y_{\delta} = Y_0 + \delta Y_1 + \ldots,
\end{align*}
where all of the terms, of course, depend on the variables $(s,z)$.

In the interface $\Gamma_\delta$ we aim for the representations (i.e., the outer and inner expansions) of $y_\delta \in \{\bu_\delta, p_\delta, \varphi_\delta, \mu_\delta \}$ to match. To do so, we have, as $z\to \pm \infty$, the following matching conditions 
\begin{align*}
    Y_0(s,z) & =  y_0(0^\pm),\\
    \partial_z Y_0(s,z) & = 0,\\
    Y_1(s,z) & = y_1^\pm + (\nabla y_0(0^\pm)\cdot\nu)z,\\
    \partial_z Y_1(s,z) & = (\nabla y_0(0^\pm)\cdot\nu) ,
\end{align*}
where for a smooth function $f$
\begin{align*}
    f(0^\pm)(\gamma) := \lim_{r\to 0} f(\gamma\mp r \nu(\gamma))\quad\text{for }\gamma\in\Gamma.
\end{align*}
We refer to \cite{GS2006-IFB} for the detailed derivation.

We can follow the arguments used in \cite[Section~5.2.1]{GHHKL2016-IFB} to show that the Dirichlet condition for the Navier--Stokes equations on $\Gamma$ holds. Indeed, from the leading order $\mathcal{O}(\delta^{-1})$ in the expansion of \eqref{delsolNS}, we have 
\begin{align}\label{ord1div}
    \partial_z \bU_0 \cdot \nu = 0.
\end{align}
On the other hand, the leading order $\mathcal{O}(\delta^{-2})$ in \eqref{delNS} is
\begin{align}
    \frac{1}{2}\partial_z \left( \eta_0 (\partial_z \bU_0 + (\partial_z \bU_0 \cdot \nu)\nu) \right) = 0,
\end{align}
where $\eta_0$ is evaluated as $\eta(\Phi_0)$. Taking the integral over $z$ and invoking \eqref{ord1div} to the equation above leads to
\begin{align}
    \eta_0\partial_z \bU_0 = A(s),
\end{align}
where $A(s)$ is a function independent on $z$. Form the matching conditions we get
\begin{align*}
    A(z) = \lim_{z\to \pm\infty} \eta_0\partial_z \bU_0 = \eta(\phi_0(0^{\pm})) \cdot 0 = 0.
\end{align*}
From the computation above and the fact that $\eta_0 >0$, we infer that $\partial_z \bU_0 = 0$. Integrating over $z$ gives us $\bU_0 = B(s)$ for some function $B(s)$ independent on $z$. From the matching conditions and \eqref{UinB}, we get
\begin{align}
    B(s) = \lim_{z\to -\infty} \bU_0 = \bu_0(0^-) = 0.
\end{align}
On the other hand, we have 
\begin{align}
    \bu_0(0^+) = \lim_{z\to +\infty} \bU_0 = B(s) = 0.
\end{align}

To proceed with recovering the Neumann condition for the limiting concentration $\varphi_0$, we draw a particular attention to the interpolant $\alphi{\cdot}$ which is given as
\begin{align}  
    \alphi{\Psi_\delta(s,z)} &= \frac{\Psi_0(s,z)+1}{2} 
    + \delta \left(\frac{\underline{\alpha}_\varphi(1-  \Psi_0(s,z))}{2} \right) =: \widehat{\alpha}_\varphi^{(0)}(s,z) + \delta \widehat{\alpha}_\varphi^{(1)}(s,z).
\end{align}
We infer from Assumption~\ref{expansion:assumption}\ref{interior:assumption} that $\widehat{\alpha}_\varphi^{(0)}$ is independent on $s$, and
\begin{align}\label{alphalim}
    \lim_{z\to +\infty} \widehat{\alpha}_\varphi^{(0)}(z) = 1,\quad
    \lim_{z\to -\infty} \widehat{\alpha}_\varphi^{(0)}(z) = 0,\quad \lim_{z\to +\infty} \widehat{\alpha}_\varphi^{(1)}(z) = 0,
    \text{ and }
    \lim_{z\to -\infty} \widehat{\alpha}_\varphi^{(1)}(z) = \underline{\alpha}_\varphi.
\end{align}
Plugging the inner expansions into \eqref{delmueq}, the leading order $\mathcal{O}(\delta^{-2})$ gives us
\begin{align}
    \partial_z\left(\widehat{\alpha}_{\varphi}^{(0)} \partial_z \Phi_0 \right) = 0.
\end{align}
Integrating over $z$, we get $C(s) = \widehat{\alpha}_{\varphi}^{(0)} \partial_z \Phi_0 $, where $C(s)$ is a function dependent only on $s$ but independent on $z$. 
Thus, by using the limits of $\widehat{\alpha}_{\varphi}^{(0)}$ and $\partial_z \Phi_0$, which we infer from the matching condition 
for $\partial_z \Phi_0$, we get 
\begin{align}
    C(s) & = \lim_{z \to \pm \infty} \widehat{\alpha}_{\varphi}^{(0)} \partial_z \Phi_0  = \left\{ 
    \begin{aligned}
        1\cdot 0 &&\text{as }z\to +\infty\\ 
        0\cdot 0 &&\text{as }z\to -\infty
    \end{aligned} 
    \right\}  = 0.
\end{align}
From this, we infer that $\widehat{\alpha}_{\varphi}^{(0)}\partial_z\Phi_0 = 0 $ everywhere.
Moreover, it holds $\widehat{\alpha}_{\varphi}^{(0)} >0$ for any finite $z$ and thus $\partial_z\Phi_0 = 0$ for any finite $z$.

The next singular term $\mathcal{O}(\delta^{-1})$ in the expanded \eqref{delmueq} gives us
\begin{align}\label{ord1sing}
    0 = \partial_z\left( \widehat{\alpha}^{(0)}_\varphi \partial_z \Phi_1 + \widehat{\alpha}^{(1)}_\varphi \partial_z \Phi_0  \right) + \partial_z\left(\widehat{\alpha}^{(0)}_\varphi \nabla_\Gamma \Phi_0  \right)\cdot\nu + \dive_{\Gamma}\left( \widehat{\alpha}^{(0)}_\varphi \partial_z \Phi_0 \nu \right).
\end{align}
Since the normal vector $\nu$ is independent on $z$, we can pull $\nu$ inside the derivative in the second term above, and since the tangential gradient is perpendicular to the normal vector the second term vanishes, i.e. $\nabla_\Gamma \Phi_0 \cdot \nu = 0$. Furthermore, since $\widehat{\alpha}_{\varphi}^{(0)}\partial_z \Phi_0 = 0$, \eqref{ord1sing} simplifies to
\begin{align}
     0 = \partial_z\left( \widehat{\alpha}^{(0)}_\varphi \partial_z \Phi_1 + \widehat{\alpha}^{(1)}_\varphi \partial_z \Phi_0  \right)
\end{align}
Integrating over $z$ leads to $\widehat{\alpha}^{(0)}_\varphi \partial_z \Phi_1 + \widehat{\alpha}^{(1)}_\varphi \partial_z \Phi_0 = D(s)$, where $D(s)$ is a function dependent only on $s$ but independent on $z$. From the matching conditions we have
\begin{align*}
    D(s) = \lim_{z\to -\infty} \widehat{\alpha}^{(0)}_\varphi \partial_z \Phi_1 + \widehat{\alpha}^{(1)}_\varphi \partial_z \Phi_0 = 0\cdot [\nabla\varphi_0(0^-)\cdot \nu] + \underline{\alpha}_\varphi\cdot 0 = 0.
\end{align*}
On the other hand, we have
\begin{align*}
    \nabla\varphi_0(0^+)\cdot \nu = \nabla\varphi_0(0^+)\cdot \nu + 0\cdot 0 = \lim_{z\to +\infty} \widehat{\alpha}^{(0)}_\varphi \partial_z \Phi_1 + \widehat{\alpha}^{(1)}_\varphi \partial_z \Phi_0= D(s) = 0.
\end{align*}
Using the same arguments on \eqref{delphieq}, we can show that $\nabla\mu_0(0^+)\cdot \nu = 0$, in fact we also get that
\begin{align}\label{Muprop}
    \widehat{\alpha}_{\varphi}^{(0)} \partial_z M_0 = 0, \text{ and } \widehat{\alpha}^{(0)}_\varphi \partial_z M_1 + \widehat{\alpha}^{(1)}_\varphi \partial_z M_0  = G(s),
\end{align}
where, through \eqref{alphalim} and the matching conditions, one can establish $G(s) = 0$.
From these we get the Neumann conditions for the concentration $\varphi_0$ and the chemical potential $\mu_0$ on $\Gamma$.


\begin{remark}
     In summary, we see that the element $(\bu_0,p_0,\varphi_0,\mu_0)$, for which $(\bu_\delta,p_\delta,\varphi_\delta,\mu_\delta)$ converges to, satisfies the following system
     \begin{subequations}\label{system:asymplimit}
        \begin{align}
                \frac{1}{2}\big[\rho(\varphi_0)\partial_t \bu_0 + \partial_t(\rho(\varphi_0)\bu_0) + ((\rho(\varphi_0)\bu_0 + \bJ_0)\cdot\nabla)\bu_0&\notag\\  
                + \dive((\rho(\varphi_0)\bu_0 + \bJ_0) \otimes \bu_0) \big] - \dive(2 \eta(\varphi_0)\D\bu_0) + \nabla p_0  
                &= -\varphi_0\nabla\mu_0 + \rho(\varphi_0)\bg && 
            \text{ in }(0,T]\times E, \label{asympNS} \\
             \dive \bu_0 &= 0 &&\text{ in }(0,T]\times E,\label{asympdiv}\\
             \partial_t\varphi_0 + \dive(\varphi_0\bu_0 - b \nabla\mu_0) &= 0 &&\text{ in }(0,T]\times E,\label{asympphieq} \\
             - \sigma{\varepsilon}\Delta \varphi_0 + \frac{\sigma}{{\varepsilon}}W'(\varphi_0) &= \mu_0 && \text{ in }(0,T]\times E,\label{asympmueq}
        \end{align}
    \end{subequations}
    with boundary conditions $\bu_0 = 0$ and $\partial_\nu\varphi_0 = \partial_\nu \mu_0 = 0$ on $(0,T]\times \partial E$ and initial conditions $\bu_0(0) = \bu^0|_{E}$ and $\varphi_0(0) = \varphi^0|_{E}$ in $E$.

    Moreover, by defining $R_0 = -b\nabla\frac{\partial \rho}{\partial\varphi}(\varphi_0)\cdot \nabla\mu_0$, \eqref{asympphieq} when multiplied to $\frac{\partial\rho}{\partial\varphi}(\varphi_0)$ gives us 
    \begin{align*}
        \partial_t \rho(\varphi_0) + \dive(\rho(\varphi_0)\bu_0 + \bJ_0) = R_0.
    \end{align*}
    Applying this identity to \eqref{asympNS} leads to
    \begin{align}
        \rho(\varphi_0)\partial_t\bu_0 + ((\rho(\varphi_0)\bu_0 + \bJ_0)\cdot\nabla)\bu_0 + R_0\frac{\bu_0}{2} - \dive(2\eta(\varphi_0)\D \bu_0) + \nabla p_0 = - \varphi_0\nabla\mu_0 + \rho(\varphi_0)\bg.
    \end{align}
    This shows that $(\bu_0,p_0,\varphi_0,\mu_0)$ solves \eqref{system:prime}.
\end{remark}

\subsubsection{On the conservation of the average of $\varphi_\delta$.}
\label{sec:conservation-phiavg}

    In Remark~\ref{rm:conservation-phiavg}, we
    calculated
    \begin{align*}
        \frac{d}{dt}\int_{\Omega}\frac{\psi + 1}{2}\varphi \du x 
        = \frac{1}{2}\int_{\Gamma_\delta}\varphi\bu\cdot \nabla \psi \du x 
        - \frac{b}{2}\int_{\Gamma_\delta} \alphi{\psi}\nabla\mu\cdot \nabla \psi \du x
    \end{align*}
    and we conjectured that as $\delta \to 0 $ the integrals involving $\Gamma_\delta$ vanishes. In the same spirit as the formal computations as above, we show that indeed these integrals go to zero. 

    With the introduced $(s,z)$-coordinates, we can write $\Gamma_\delta$ as 
    \begin{align}
        \Gamma_\delta = \{ (\gamma(s) + \delta z\nu(\gamma(s)))\in \Omega : 0\le s\le \ell,\, z\in \mathbb{R}\}
    \end{align}
    for some fixed constant $\ell>0$. 
    
    Letting $T:\mathbb{R}^2\to \mathbb{R}^d$ be the transformation defined as 
        $$x = T(s,z) = \gamma(s) + \delta z\nu(\gamma(s)),$$
    we compute the first integral as 
    \begin{align*}
        \int_{\Gamma_\delta} \varphi_\delta\bu_\delta\cdot \nabla\psi_\delta \du x 
        = \lim_{r\to \infty}\int_0^\ell \int_{-r}^{r} \Phi_\delta\bU_\delta \cdot \left( \frac{1}{\delta}(\partial_z \Psi_\delta )\nu + \nabla_\Gamma \Psi_\delta 
        + \ldots \right) \sqrt{ \mathrm{det}(J_T^{\top}J_T) } \du z \du s,
    \end{align*}
    where $J_T$ is the Jacobian of the transformation $T$. By noting that, for sufficiently small $\delta>0$, we have
    \begin{align}
        \sqrt{ \mathrm{det}(J_T^{\top}J_T) } = \delta(|\gamma'(s)| - \delta z \kappa ),
    \end{align}
    and by using the expansions for $\Phi_\delta$, $\bU_\delta$ and $\Psi_\delta$, we can further compute that
    \begin{align*}
       \int_{\Gamma(\delta)} \varphi_\delta\bu_\delta\cdot \nabla\psi_\delta \du x 
       &= \lim_{r\to \infty}\int_0^\ell \int_{-r}^{r} \Phi_\delta\bU_\delta \cdot \left( \frac{1}{\delta}(\partial_z \Psi_\delta )\nu + \nabla_\Gamma \Psi_\delta + \ldots \right) \sqrt{ \mathrm{det}(J_T^{\top}J_T) } \du z \du s\\
         &= \lim_{r\to \infty}\int_0^\ell \int_{-r}^{r} \Phi_0\bU_0 \cdot (\partial_z \Psi_0 )\nu\, |\gamma'(s)| \du z \du s + \mathcal{O}(\delta)
    \end{align*}
    We recall from the previous part that $\partial_z\bU_0 = \boldsymbol{0}$ and $\partial_z\Phi_0 = 0$, i.e., $\bU_0$ and $\Phi_0$ are independent on $z$. 
    We hence get at order $\delta^0$
    \begin{align*}
        \lim_{r\to \infty}\int_0^\ell \int_{-r}^{r} \Phi_0\bU_0 \cdot (\partial_z \Psi_0 )\nu\, |\gamma'(s)| \du z \du s 
        &= \lim_{r\to \infty}\int_0^\ell \Phi_0(\bU_0 \cdot \nu) |\gamma'(s)|\du s \int_{-r}^{r} (\partial_z \Psi_0 ) \du z\\
        & = \lim_{r\to \infty}[\Psi_0(r) - \Psi_0(-r)]\int_\Gamma \varphi_0(\bu_0 \cdot \nu) \du \sigma = 0,
    \end{align*}
    where the last line is achieved via the transformation $\sigma = \gamma(s)$ for elements $\sigma\in\Gamma$, the boundary condition $\bu_0 = 0$ on $\Gamma$, and $\Psi_0(r) \to 1$ and $\Psi_0(-r) \to -1$ as $r\to\infty$.
    Therefore, as $\delta \to 0$ we have $\int_{\Gamma(\delta)} \varphi_\delta\bu_\delta\cdot \nabla\psi_\delta \du x \to 0$.

    Similarly, we can show,
    that
    \begin{align*}
        &\int_{\Gamma_\delta} \alphi{\psi}\nabla\mu_\delta\cdot \nabla\psi_\delta \du x \to 0\quad\text{as }\delta\to 0.
    \end{align*}
    Indeed, we have
    \begin{align*}
        &\int_{\Gamma_\delta} \alphi{\psi}\nabla\mu_\delta\cdot \nabla\psi_\delta \du x\\ & = \lim_{r\to\infty}\int_0^\ell \int_{-r}^r \left[ \frac{1}{\delta^2} \widehat{\alpha}^{(0)}_\varphi \partial_z M_0 \partial_z \Psi_0 + \frac{1}{\delta}\left(\widehat{\alpha}^{(0)}_\varphi\partial_z M_1 \partial_z + \widehat{\alpha}^{(1)}_\varphi\partial_z M_0 \right)\partial_z  \Psi_0  \right]\delta(|\gamma'(s)| - \delta z\kappa) \du z\du s + \mathcal{O}(\delta).
    \end{align*}
    Finally, from \eqref{Muprop}, we infer that the integrals with $\widehat{\alpha}^{(0)}_\varphi \partial_z M_0$ and $\widehat{\alpha}^{(0)}_\varphi\partial_z M_1 + \widehat{\alpha}^{(1)}_\varphi\partial_z M_0$ equate to zero.

    
\section{The time-discrete system}
\label{sec:discrete}

To transition from a continuous theoretical framework to a format suitable for practical numerical implementation, we establish a well-posed discrete scheme.
We begin with stating a formal weak formulation for the continuous system \eqref{system:pores} 
and derive an energy equality that solutions to this system satisfy. In Section~\ref{sec:discrete} we thereafter propose a time discretization scheme that mimics the energy equality in the time discrete setting.

\begin{remark}
    In the following we restrict to the solenoidal setting and restrict our argumentation to solenoidal velocity fields. Therefore especially the pressure $p$ and the solenoidality condition on $\bu$ drop from the system. By using standard techniques the pressure can be recovered, see e.g. \cite[Lem.~I~2.1]{GiraultRaviart_FEM_for_NavierStokes}.
\end{remark}
\begin{definition}
    We formally define a weak formulation of the continuous system.
       We say that the sufficiently smooth elements $\bu:[0,T]\times\Omega\to \mathbb{R}^d$, $\dive{\bu}=0$,
       $\varphi:[0,T]\times\Omega\to \mathbb{R}$ and $\mu:[0,T]\times\Omega\to \mathbb{R}$ make up a weak solution to \eqref{system:pores} if they satisfy the system
    \begin{subequations}
    \label{system:weakpores}
        \begin{align}
         \frac{1}{2} (( \rho\partial_t \bu + \partial_t(\rho\bu)  ), \bw)_\Omega+ a(\rho\bu + \bJ_\delta, \bu, \bw) & \notag \\
           +  (\alphu{\psi}\bu, \bw)_\Omega + 2( \eta \D\bu, \D\bw)_\Omega &= - (\varphi\nabla\mu,\bw)_\Omega   + (\rho \bg,\bw)_\Omega, \\
            (\partial_t \varphi, \zeta) - (\varphi\bu,\nabla\zeta)_\Omega + b(\alphi{\psi}\nabla\mu,\nabla\zeta)_\Omega &= 0, 
            \label{eq:weak-CH1}\\
             \sigma\varepsilon(\alphi{\psi}\nabla\varphi,\nabla\vartheta)_\Omega + \frac{\sigma}{\varepsilon}(W'(\varphi),\vartheta)_\Omega &= (\mu,\vartheta)_\Omega, 
            \label{eq:weak-CH2}
        \end{align}
    \end{subequations}
    for all $\bw\in V_\sigma$, $\zeta,\vartheta\in V$ and for a.e. $t\in (0,T]$, and the initial conditions $\bu(0) = \bu^0 \in H_\sigma$ and $\varphi(0) = \varphi^0\in V $.
\end{definition}

\begin{remark}
   We note that the trilinear form is achieved by invoking the property
    \begin{align*}
        a(\bv,\bu,\bw) &= \frac{1}{2}\left\{ 
        \int_\Omega \left[(\bv \cdot \nabla)\bu\right]  \cdot \bw \dx - 
          \int_\Omega \left[(\bv \cdot \nabla) \bw\right] \cdot \bu \dx
        \right\}\\
        &= \frac{1}{2}\left\{ \int_\Omega [(\bv\cdot\nabla)\bu]\cdot \bw \du x + \int_\Omega \dive(\bv\otimes\bu)\cdot \bw \du x \right\}
    \end{align*}
	as long as $\bu$ and $\bw$ satisfy $\bu\cdot\nu_\Omega = \bw\cdot\nu_\Omega = 0$ on $\partial\Omega$.
\end{remark}

System \eqref{system:weakpores} is thermodynamically consistent in the sense, that it satisfies the following energy equality.     
    \begin{lemma}
          Let  $\bu:[0,T]\times\Omega\to \mathbb{R}^d$, $\varphi:[0,T]\times\Omega\to \mathbb{R}$ and $\mu:[0,T]\times\Omega\to \mathbb{R}$ satisfy system \eqref{system:weakpores}. 
          Then they satisfy the following energy equality.
        \begin{align}
            \begin{aligned}
            \frac{d}{dt}\left[ 
            \int_\Omega \frac{1}{2}\rho|\bu|^2 
            + \sigma 
            \left( \frac{\varepsilon}{2}\alphi{\psi}|\nabla\varphi|^2 
            + \frac{1}{\varepsilon}W(\varphi)\right)
            \du x 
            \right]& \\ 
              + \int_\Omega \alphu{\psi}|\bu|^2 \du x 
              +\int_\Omega 2\eta|D\bu|^2 \du x 
            &+ \int_\Omega b\alphi{\psi}|\nabla\mu|^2 \du x 
            = \int_\Omega \rho\bg\cdot\bu \du x.
            \end{aligned}
            \label{energy:total}
        \end{align}
    \end{lemma}
    \begin{proof}
        Indeed, take $\bw = \bu$, $\zeta = \mu$ and $\vartheta = \partial_t\varphi$ in \eqref{system:weakpores} and taking the sum give us the energy equality.
    \end{proof}

    \begin{remark}
    \label{rm:energyEquality}
        By integrating with respect to time, the energy equality \eqref{energy:total} can equivalently be written as
        \begin{align}
            \mathcal E(t)
             &+ \int_0^t \int_\Omega \alphu{\psi}|\bu|^2 \du x \du t
             + \int_0^t \int_\Omega 2\eta|D\bu|^2 \du x \du t
             +  \int_0^t \int_\Omega  b\alphi{\psi}|\nabla\mu|^2 \du x \du t \\
           &  = \mathcal E(0) +  \int_0^t \int_\Omega \rho\bg\cdot\bu \du x \du t, 
             \label{energy:total:integrated}
        \end{align}
        where
        \begin{align}
          \mathcal E(t) = \int_\Omega \frac{1}{2}\rho(\varphi(t))|\bu(t)|^2\du x 
           + \sigma \int_\Omega \left( \frac{\epsilon}{2}\alphi{\psi}| \nabla \varphi(t)|^2 + \frac{1}{\epsilon}W(\varphi(t)) \right) \du x
           \label{energy:def}
        \end{align}
        denotes the energy of the system at some time $t\in I$. 
        Here, the first integral denotes the kinetic energy of the system and the second integral denotes the interfacial energy, which is approximately the product of surface tension and interfacial length.
        The time discrete scheme that we propose in Definition~\ref{defn:weak} mimics this equality in a time discrete setting, see Corollary~\ref{cor:energy:global}.
    \end{remark}

\subsection{Statement of the time-discrete system}  

To discretize the time domain $(0,T]$ we introduce time points $0=t_0 < \ldots<t_{m-1}< t_m < \ldots < t_M=T$, that are for simplicity of presentation equidistant with distance $\tau = t_m-t_{m-1}$.
For a function $y \in \{\bu,\varphi,\mu\}$ we approximate $y^m \approx y(t_m)$ and thus approximate the time depending function $y$ by a sequence $(y^m)_{m=1}^M$ of $M$ time independent functions.

Before we state the time stepping scheme to approximate system \eqref{system:weakpores} we impose some assumptions on the domain $\Omega$, the potential $W$, the density $\rho$ and the viscosity $\eta$.

\begin{assumption}\mbox{}
\label{assumption:timediscrete}
    \begin{enumerate}[label = \Roman*.]
    \item For $d=2$ the domain $\Omega$ is a bounded polygonal domain. The initial data $\varphi_0$ satisfies $\varphi_0 \in L^{\infty}(\Omega)$ and also $\alphi{\psi} \in L^\infty(\Omega)$.
    
    \item For $d=3$ the boundary $\partial \Omega$ of $\Omega$ is of class $C^{1,1}$ and the initial data $\varphi_0$ satisfies $\varphi_0 \in C^{0,1}(\overline \Omega)$. Moreover $\alphi{\psi} \in C^{0,1}(\overline \Omega)$ holds.
    
        \item The potential $W$ is radially unbounded and  can be decomposed into $W = W_+ + W_{-}$, where $W_+$ is a convex function
            and $W_-$ is a concave function. 
            Furthermore, we assume that $W$, $W_+$ and $W_-$ satisfy
            \begin{align}
                W(s) \le C(1 + |s|^4),\quad 
                |W_+'(s)| \le C(1 + |s|^3),\quad 
               |W_-'(s)| \le C(1 + |s|^3), \quad\forall s\in\mathbb{R}.\label{est:potential}
            \end{align}
        \item 
            The density and  viscosity satisfy $\rho \in W^{1,\infty}(\mathbb{R})$ and $\eta \in L^{\infty}(\mathbb{R})$. We suppose further that there exist constants $\rho_{min},\rho_{max},\eta_{min},\eta_{max}>0$ such that
            \begin{align*}
                \rho_{min} \le \rho(s) \le \rho_{max},\quad \eta_{min} \le \eta(s) \le \eta_{max}, \quad \forall s\in\mathbb{R}.
            \end{align*}
    \end{enumerate}
\end{assumption}

\begin{notation}
\label{not:PHI}
For $d=2$ we define the space $\Phi = L^{\infty}(\Omega)$. 
For $d=3$ we define the space $\Phi = C^{0,1}(\overline \Omega)$.
\end{notation}
        
We are now in the position to state the finite difference time-stepping scheme for the numerical approximation of system \eqref{system:weakpores}.
      \begin{definition}
      \label{defn:weak}
		Given $(\bg^m)_{m=1}^M \in (L^2(\Omega)^d)^M$, 
        $\bu^0 \in H_\sigma$ 
        and 
        $\varphi^0 \in V\cap \Phi$, 
        we say that 
        $(\bu^m)_{m=1}^M \in V_\sigma^M$, 
        $(\varphi^m)_{m=1}^M \in (V\cap \Phi)^M$ 
        and $(\mu^m)_{m=1}^M \in W^{1,3}(\Omega)^M$ 
        constitute a semi-discrete weak solution to \eqref{system:weakpores} if they satisfy the following difference equation for all $\bw\in V_\sigma$, $\zeta,\vartheta\in V$ and  for all $m=1,\ldots,M$
		\begin{subequations}
        \label{system:weakdiscretepores}
			\begin{align}
					\frac{1}{2\tau} \left( (\rho^m + \rho^{m-1}) \bu^{m} - 2\rho^{m-1}\bu^{m-1}  , \bw \right)_\Omega+ a(\rho^{m-1}\bu^{m-1} + \bJ_\delta^{m}, \bu^{m}, \bw) 
                    & \notag\\
					+  (\alphu{\psi}\bu^{m}, \bw)_\Omega +\, 2( \eta^{m} \D\bu^{m}, \D\bw)_\Omega +  (\varphi^{m-1}\nabla\mu^{m},\bw)_\Omega & = (\rho^{m} \bg^{m},\bw)_\Omega,
				     \label{NSdisc}\\
                    \frac{1}{\tau}(\varphi^{m} - \varphi^{m-1}, \zeta) - (\varphi^{m-1}\bu^{m-1},\nabla\zeta)_\Omega + {\tau}\left( \frac{|\varphi^{m-1}|^2}{\rho^{m-1}}\nabla\mu^m, \nabla \zeta \right)& \notag \\ 
                    +b(\alphi{\psi}\nabla\mu^{m},\nabla\zeta)_\Omega & = 0,
               \label{phieqdisc}\\
				\sigma\varepsilon(\alphi{\psi}\nabla\varphi^{m},\nabla\vartheta)_\Omega + \frac{\sigma}{\varepsilon}(W_+'(\varphi^{m}) + W_-'(\varphi^{m-1}),\vartheta)_\Omega &= (\mu^{m},\vartheta)_\Omega.\label{mueqdisc}
			\end{align}
		\end{subequations}        
        Here, we use the abbreviation $\rho^m:= \rho(\varphi^m)$, $\eta^m = \eta(\varphi^m)$ and $\bJ^m_\delta = -b\alphi{\psi}\rho^\prime(\varphi^m)\nabla\mu^m$ for $m = 0,1,\ldots, M$.
	\end{definition}

   \begin{remark}
       Note that scheme \eqref{NSdisc}--\eqref{mueqdisc} is a sequentially coupled scheme that allows us to solve for the variables $\varphi^m$ and $\mu^m$ from \eqref{phieqdisc}--\eqref{mueqdisc} first before resolving the fluid equations \eqref{NSdisc}. 
       The difference equation \eqref{phieqdisc} is in fact obtained by estimating the term $(\varphi\bu,\nabla\zeta)_\Omega$ by $(\varphi^{m-1}\bu^\star,\nabla\zeta)_\Omega$ where $\bu^\star$ is an intermediate velocity field---not necessarily a solenoidal field---that is assumed to take care of the capillary forces in the continuity equation \eqref{NS}, and satisfies
       \begin{align}
           \frac{1}{\tau}\rho^{m-1}(\bu^\star - \bu^{m-1}) = -\varphi^{m-1}\nabla\mu^m.\label{intermediate}
       \end{align}
       We note that such decoupling method has been used in several works, these include the following references \cite{M2012-NMPDE,GGT2014-JCM,GGGM2016-CCP,BKR2019-JCP,BK2021-SICON}.
   \end{remark}

\begin{remark}
\label{rem:higherRegularityNecessary}
   The regularity $\varphi^{m-1} \in L^\infty(\Omega)$ is necessary for $d=2$ and the regularity $\varphi^{m-1} \in C^{0,1}(\overline \Omega)$ is sufficient in $d=3$ because $|\varphi^{m-1}|^2$ appears as coefficient in the diffusion operator in \eqref{phieqdisc}.
   
   The regularity $\mu^m \in W^{1,3}(\Omega)$ is necessary such that $J_\delta \bu $ has sufficient regularity to make the trilinear form $a(\cdot,\cdot,\cdot)$ in \eqref{NSdisc} well defined in three spatial dimensions.
\end{remark}

\begin{remark}
    By using $\zeta = 1$ in \eqref{phieqdisc} we obtain conservation of the mean value of $\varphi$, i.e. $\int_\Omega \varphi^m\du x = \int_\Omega \varphi^{m-1}\du x$.
\end{remark}

\subsection{Well-posedness of the time-discrete system}

We show existence of unique solutions to system \eqref{system:weakdiscretepores} and that these solutions satisfy a discrete version of the energy equality \eqref{energy:total}.

We start with the existence result for one time step.
\begin{proposition}
\label{prop:pointwise}
    Let $\bg^m\in L^2(\Omega)^d$, $(\bu^{m-1},\varphi^{m-1})\in H_\sigma\times V\cap \Phi$.
    Then there exist unique 
    $\bu^{m} \in V_\sigma$, 
    $\varphi^{m} \in V\cap \Phi$, 
    $\mu^{m} \in W^{1,3}(\Omega)$ 
    satisfying \eqref{system:weakdiscretepores}.
\end{proposition}
\begin{proof}
The proof directly follows the corresponding proof in \cite[Lem.~3.5]{BK2021-SICON}.
In comparison to \cite{BK2021-SICON} for \eqref{phieqdisc}--\eqref{mueqdisc} we now have a given uniformly bounded and strictly positive diffusion coefficient $\alphi{\psi}$ 
and for \eqref{NSdisc} we have an additional non-negative lower order term $\alphu{\psi}$. 
Both changes do not affect the proof of \cite[Lem.~3.5]{BK2021-SICON}.
We summarize the steps.

System \eqref{system:weakdiscretepores} is sequentially coupled. 
We therefore first discuss \eqref{phieqdisc}--\eqref{mueqdisc}. The existence of unique solutions $\varphi^m,\mu^m \in V$ follows as in \cite{HintermuellerHinzeTber} by considering a suitable minimization problem. 
The additional regularity $\varphi^m \in C(\overline \Omega)$ follows from \cite[Thm.~4.8]{Troeltzsch} and requires the monotonicity of $W^\prime_+$.
For $d=2$ the additional regularity $\mu^m \in W^{1,3}(\Omega)$ follows from \cite[Thm.~1.9]{2012_Dhamo_Diss_BoundaryControl}. 
For $d=3$, using the higher boundary regularity, we observe $\varphi^m \in W^{2,p}(\Omega)\hookrightarrow C^{0,1}(\overline \Omega)$. This in turn gives $\mu^m \in H^2(\Omega) \hookrightarrow W^{1,3}(\Omega)$.

The existence of unique $\bu^m \in V_\sigma$ directly follows from Lax--Milgram's theorem.
\end{proof}

\begin{remark}
    As an alternative, to get $\bu^m\in V_\sigma,\varphi^m \in V, \mu^m \in V$, one might use
a standard Galerkin approach and Brouwer's fixed point theorem, see e.g. \cite{GHK2016-LinearScheme,Gruen_Klingbeil_CHNS_AGG_numeric}. The additional regularity then follows as before.
\end{remark}

Because the additional regularity $\varphi^m \in \Phi$ and $\mu^m \in W^{1,3}(\Omega)$ exactly fits the necessary regularity for the existence of unique solutions in Proposition~\ref{prop:pointwise}, we argue by induction to get the existence of solutions to \eqref{system:weakdiscretepores}.

\begin{corollary}
    Let $(\bg^m)_{m=1}^M \in (L^2(\Omega)^d)^M$, $\bu^0\in H_\sigma$,  $\varphi^0\in V\cap \Phi$, and let Assumption~\ref{assumption:timediscrete} hold.
    Then there exists a unique weak solution 
    $(\bu^m)_{m=1}^M \in V_\sigma^M$, 
    $(\varphi^m)_{m=1}^M \in (V \cap \Phi)^M$, 
    $(\mu^m)_{m=1}^M \in (W^{1,3}(\Omega))^M$
    to the time discrete system \eqref{system:weakdiscretepores}.
\end{corollary}

    The continuous system \eqref{system:weakpores} satisfies the energy equality \eqref{energy:total:integrated}.
    We next show that the discrete scheme satisfies a discrete variant of \eqref{energy:total:integrated}.
    For convenience, we restate the energy of the system as defined in  \eqref{energy:def} for some $\bu^m,\varphi^m$
    \begin{align*}
        \mathcal{E}(\bu^m,\varphi^m) := 
        \left[ \int_\Omega \frac{1}{2}\rho(\varphi^m)|\bu^m|^2 \du x 
        + \sigma\int_\Omega \alphi{\psi} \frac{\varepsilon}{2}|\nabla \varphi^m|^2 + \frac{1}{\varepsilon} W(\varphi^m) \du x \right].
    \end{align*}

     \begin{lemma} 
     \label{lemma:energy:discrete}
        Given $\bu^{m-1} \in V_\sigma$ and $\varphi^{m-1} \in V\cap \Phi$.
        Let $\bu^m \in V_\sigma$, $\varphi^m \in V\cap \Phi$, $\mu^m \in W^{1,3}(\Omega)$ denote the unique solution to \eqref{system:weakdiscretepores}.

        Then the following energy inequality holds true
        \begin{align}
        \label{energy:totaldisc}
            \begin{aligned}
            \mathcal{E}(\bu^m,\varphi^m) &
            + \tau \left(2\int_\Omega \eta^m|\D \bu^m|^2 \du x 
            + b \int_\Omega \alphi{\psi}|\nabla\mu^m|^2 \du x 
            + \int_\Omega \alpha_{\bu}^{(\delta)}(\psi)|\bu^m|^2\du x
            \right)\\
            & +\left[ 
            \int_\Omega \frac{\rho^{m-1}}{2}|\bu^m-\bu^{m-1}|^2\du x 
           + \frac{\sigma\epsilon}{2}\int_\Omega \alphi{\psi}|\nabla(\varphi^m-\varphi^{m-1})|^2 \du x 
           \right]\\
            & \le\, \mathcal{E}(\bu^{m-1},\varphi^{m-1}) 
            + \tau\left( \rho^{m}\bg^m, \bu^m \right)_\Omega .                
            \end{aligned}
        \end{align}
    \end{lemma}
    \begin{proof}
        The proof follows the lines of \cite[Thm.~3.1]{GGT2014-JCM} or \cite[Thm.~3.4]{BK2021-SICON}. 
        To be precise, use $\bw = \bu^m$, $\zeta = \mu^m$ and $\vartheta = \frac{1}{\tau}(\varphi^m - \varphi^{m-1} )$ in \eqref{NSdisc}, \eqref{phieqdisc} and \eqref{mueqdisc}, respectively, and invoke the convexity of $W_+$ and concavity of $W_-$. 
        Finally, a carefully scaled Young's inequality is needed to treat the terms coming out of $\bu^\star$.
    \end{proof}

By iterating Lemma~\ref{lemma:energy:discrete}, we come to the following discrete version of the energy equality \eqref{energy:total:integrated}.
\begin{corollary}
    \label{cor:energy:global}
    Let $(\bu^m)_{m=1}^M \in V_\sigma^M$, $(\varphi^m)_{m=1}^M \in (V \times \Phi)^M$, $(\mu^m)_{m=1}^M \in (W^{1,3}(\Omega))^M$ denote a solution to \eqref{system:weakdiscretepores}.
    Then it satisfies
    \begin{align}
    \begin{aligned}
       \mathcal{E}(\bu^m,\varphi^m) &
            + \tau \sum_{k=1}^m\left(2\int_\Omega \eta^k|\D \bu^k|^2 \du x 
            + b \int_\Omega \alphi{\psi}|\nabla\mu^k|^2 \du x 
            + \int_\Omega \alpha_{\bu}^{(\delta)}(\psi)|\bu^k|^2\du x
            \right)\\
            & + \sum_{k=1}^m\left[ 
            \int_\Omega \frac{\rho^{k-1}}{2}|\bu^k-\bu^{k-1}|^2\du x 
           + \frac{\sigma\epsilon}{2}\int_\Omega \alphi{\psi}|\nabla(\varphi^k-\varphi^{k-1})|^2 \du x 
           \right]\\
        & \le\, \mathcal{E}(\bu^{0},\varphi^{0}) 
        + \tau\sum_{k=1}^m\int_\Omega \rho^{k}\bg^k\cdot \bu^k\du x.
    \end{aligned}
    \end{align}
\end{corollary}

\begin{remark}
  We point out that the matched asymptotics done in Section~\ref{sec:asymptotics} for the continuous system can be analogously applied to the strong form of the time-discrete system \eqref{system:weakdiscretepores}. 
  Same as in the continuous case, the solutions $(\bu^{m},p^m,\varphi^{m},\mu^{m})$ should be assumed to possess the asymptotic expansion as in \eqref{genexpansion}. 
  One can also utilize the intermediate velocity field $\bu^\star$, instead of the fully-expanded, in \eqref{phieqdisc}. 
  This calls for the expansion
\begin{align*}
    \bu^\star = \left(\bu^{m-1}_0 -\tau \beta_0 \nabla \mu^m_0 \right) + \delta \left(\bu^{m-1}_1 -\tau[ \beta_1\nabla \mu^m_0 +  \beta_0\nabla \mu^m_1] \right) + \ldots =: \bu^\star_0 + \delta \bu^\star_1 +\ldots,
\end{align*}
where $\beta_0 = \frac{\varphi^{m-1}_0}{\rho^{m-1}_0}$, $\beta_1 = \frac{\rho^{m-1}_0\varphi^{m-1}_1-\rho^{m-1}_1\varphi^{m-1}_0}{(\rho^{m-1}_0)^2}$. Following the same arguments as in the asymptotic analysis as above, one gets 
    \begin{subequations}
    \label{system:asymplimitdisc}
    \begin{align}
        \frac{1}{2\tau}\big[(\rho(\varphi_0^m) + \rho(\varphi_0^{m-1}))\bu_0^m - 2 \rho(\varphi_0^{m-1})\bu_0^{m-1} \big]&  \notag \\ 
        + \frac{1}{2}\big[((\rho(\varphi_0^{m-1})\bu_0^{m-1} + \bJ_0^{m})\cdot\nabla)\bu_0^{m}  + \dive((\rho(&\varphi_0^{m-1})\bu_0^{m-1} + \bJ_0^{m}) \otimes \bu_0^{m}) \big] \notag\\ 
        - \dive(2 \eta(\varphi_0^{m})\D\bu_0^{m}) + \nabla p_0^{m} 
        &= -\varphi_0^{m-1}\nabla\mu_0^{m} + \rho(\varphi_0^m)\bg
    && \text{ in }(0,T]\times E, \label{asympNSdisc} \\ 
         \dive \bu_0^{m} &= 0 &&\text{ in }(0,T]\times E,\label{asympdivdisc}\\
         \frac{1}{\tau}(\varphi_0^{m}-\varphi_0^{m-1}) + \dive(\varphi_0^{m-1}\bu_0^\star - b \nabla\mu_0^{m}) &= 0 &&\text{ in }(0,T]\times E,\label{asympphieqdisc} \\
          - \sigma{\varepsilon}\Delta \varphi_0^{m} + \frac{\sigma}{{\varepsilon}}\left(W_+'(\varphi_0^{m}) + W_-'(\varphi_0^{m-1})\right) &=\mu_0^{m} &&\text{ in }(0,T]\times E,\label{asympmueqdisc}
    \end{align}
\end{subequations}
 with boundary conditions $\bu_0^m = 0$ and $\partial_\nu\varphi_0^m = \partial_\nu \mu_0^m = 0$ on $(0,T]\times \partial E$ and initial conditions $\bu_0^0 = \bu^0|_{E}$ and $\varphi_0^0 = \varphi^0|_{E}$ in $E$. 
The resulting time-discrete system \eqref{system:asymplimitdisc} is a time-discretization of \eqref{system:prime}.
\end{remark}


\section{Numerical investigation}
\label{sec:numerics}
As a numerical example we simulate a droplet that slides over an obstacle.
We validate the  energy inequality of the scheme and
 investigate how the overflowing liquid penetrates the porous obstacle.
In particular, we investigate the behavior with respect to variations in $\delta \to 0$.

For this we consider the physical setup as outlined in the introduction in Figure~\ref{fig:intro:schematicsetup}, which is adapted from \cite{BJKR2019CET}.
In a domain $\Omega$ we consider the left boundary to be a wall to which a triangular obstacle $B$ is attached. A liquid droplet is attached to the wall and slides down the wall under the force of gravity, passing over the obstacle. The second fluid is given by gas that surrounds the droplet.

\begin{figure}
     \centering
          \begin{subfigure}{0.45\textwidth}
         \centering
         \includegraphics[width=\textwidth]{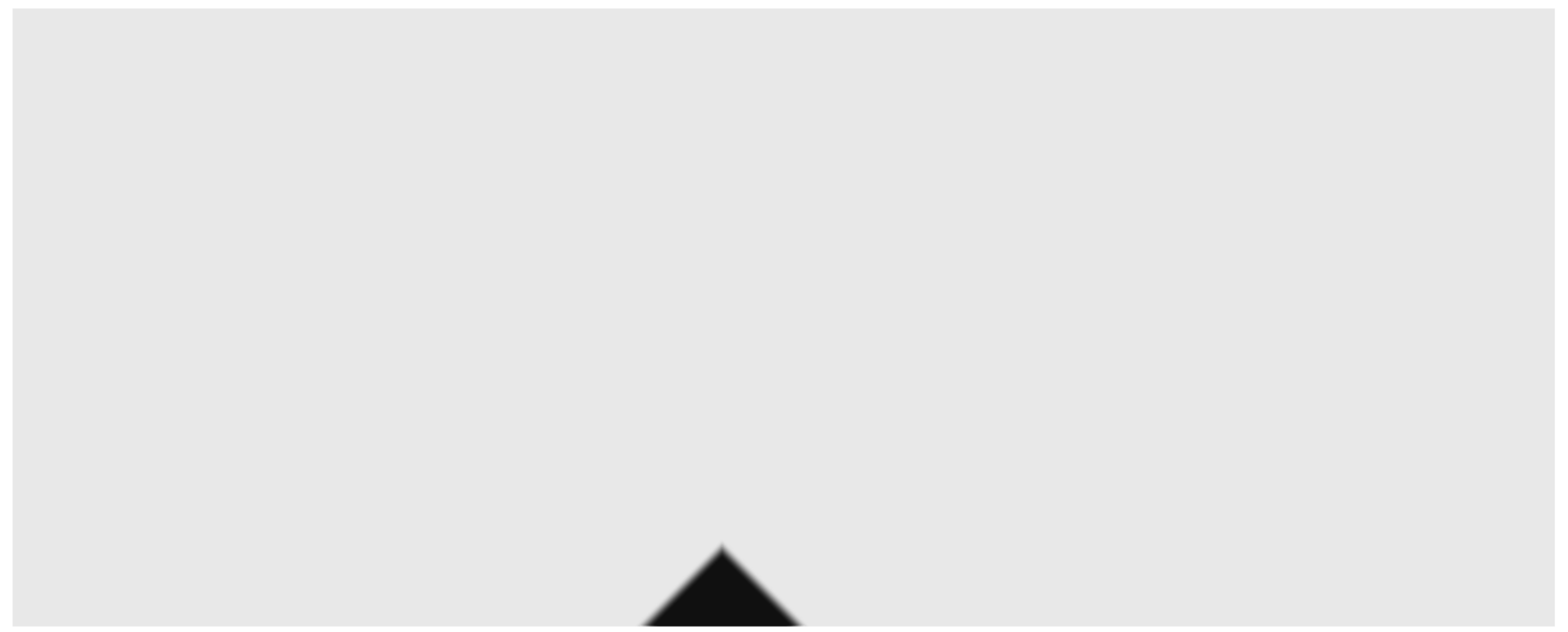}
         \caption*{(a)}
     \end{subfigure}
     \begin{subfigure}{0.45\textwidth}
         \centering
         \includegraphics[width=\textwidth]{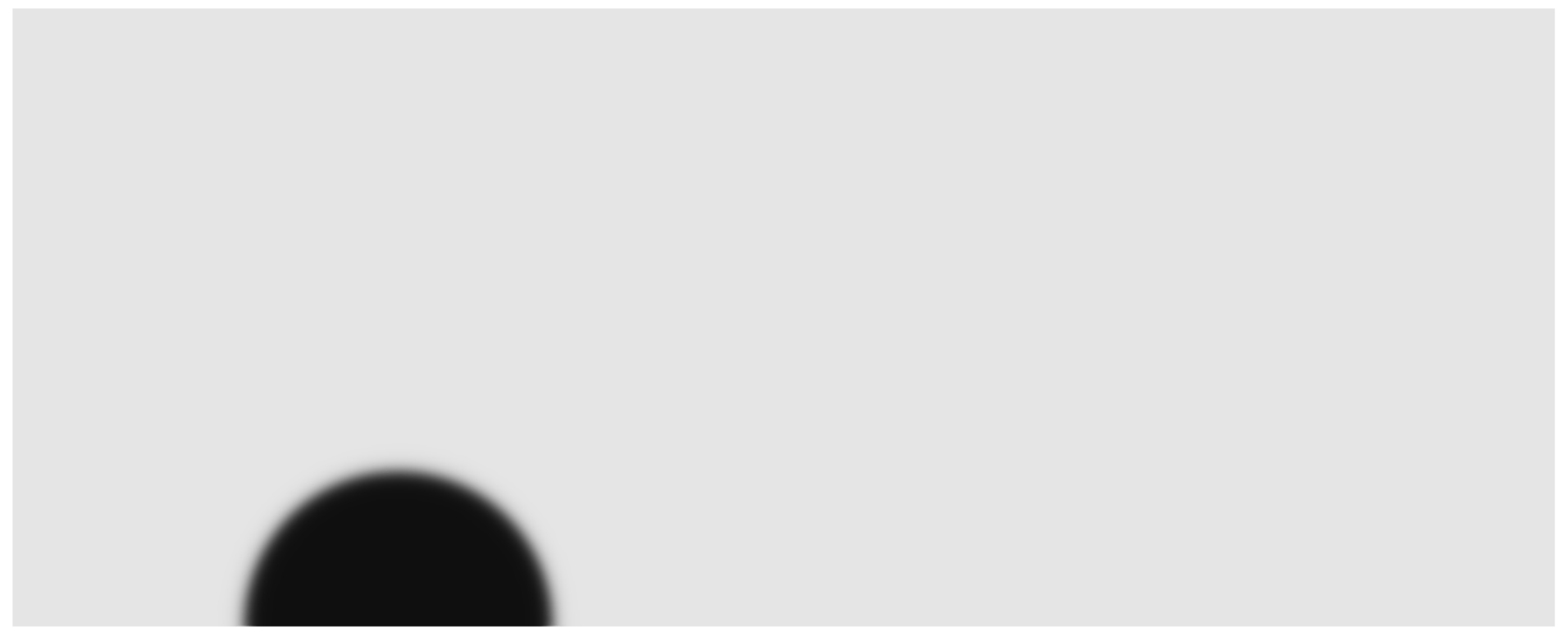}
         \caption*{(b)}
     \end{subfigure}

   \caption{This figure shows the initial set-up of the droplet example. The images are rotated by $90\degree$ and the $x$-axis is pointing upwards, while gravity points from left to right.
    Image (a) visualizes the obstacle parametrized by $\psi$, where $\psi=-1$ inside the black region, that corresponds to the impermeable region, i.e. the obstacle, and $\psi=1$ in the lighter region, that is the fluid region.
    Image (b) shows the initial $\varphi_0$, where $\varphi_0 = 1$ in the darker region---indicating the liquid--- and $\varphi_0 = -1$ in the lighter region---which corresponds to the gas.}
     \label{fig:daminit}
\end{figure}

\paragraph*{Setup - geometry}\mbox{}\\
We consider the domain $\Omega = (0,1) \times (0,2.5)$.
The obstacle $B$ is a triangular structure that is centered at $(0,1.35)$ and of height $0.125$ and given by $B = \{ x \in \Omega \mid \|x-(0,1.35)\|_1 \leq 0.125\}$. 
We describe this obstacle with the phase field $\psi$ given by 
\begin{align}
\label{psigiven}
    \psi(x) = \widetilde{\Psi}_0\left( \frac{1}{\delta}\big(\|x - (0,1.35)\|_1-0.125 \big) \right)
    \quad \text{ where }\quad
    \widetilde{\Psi}_0(z) = 
    \left\{
        \begin{aligned}
            &\sin({z}) &&\text{if }|z| \le \frac{\pi}{2},\\
            &\mathrm{sign}(z)&&\text{otherwise}.
        \end{aligned}    
    \right.
\end{align}
We show $\psi$ in Figure~\ref{fig:daminit}(a).

\begin{remark}
We note that the formal inner expansion in Section~\ref{sec:asymptotics} is derived for a sufficiently smooth obstacle boundary.
    This triangular obstacle does not satisfy this regularity assumption at its vertices.
    We consider this as a numerical robustness test.
   
\end{remark}


\paragraph*{Setup - two-phase fluid}\mbox{}\\
For the fluid density $\rho$, let $\rho_1$ and $\rho_2$ denote the densities of the involved fluids. Let $\rho(-1) = \rho_1$ and $\rho(1) = \rho_2$ where we assume that $0<\rho_1<\rho_2$. Furthermore, we use the following smooth, nonlinear, monotone and strictly positive function
\begin{align*}
    \rho(\varphi) = \left\{ 
        \begin{aligned}
        & \frac{1}{4}\rho_1 &&\text{if } \varphi \le -A_\rho^{-1},\\
        &\frac{1}{\rho_1}\left( \frac{\rho_2-\rho_1}{2}\varphi + \frac{\rho_2+\rho_1}{2} \right)^2 + \frac{1}{4}\rho_1 &&\text{if } -A_\rho^{-1} < \varphi < -1-\frac{\rho_1}{\rho_2-\rho_1},\\
        &\frac{\rho_2-\rho_1}{2}\varphi + \frac{\rho_2+\rho_1}{2} &&\text{if }-1-\frac{\rho_1}{\rho_2-\rho_1}\le \varphi \le 1+\frac{\rho_1}{\rho_2-\rho_1},\\
        &-\frac{1}{\rho_1}\left( \frac{\rho_2-\rho_1}{2}\varphi - \frac{\rho_2+\rho_1}{2}\right)^2 + \rho_2+\frac{3}{4}\rho_1 &&\text{if }1+\frac{\rho_1}{\rho_2-\rho_1} < \rho < A_\rho^{-1},\\
        &\rho_2 + \frac{3}{4}\rho_1 &&\text{if }A_\rho^{-1} \le \varphi,
        \end{aligned}
    \right.
\end{align*}
where $A_\rho = \frac{\rho_2-\rho_1}{\rho_2+\rho_1}$ denotes the Atwood number of the two-phase fluid.
 For the  viscosity we denote the viscosities of the two fluids by $\eta_1$ and $\eta_2$ and use $\eta(-1) = \eta_1$ and $\eta(1) = \eta_2$ with $0 < \eta_1 < \eta_2$, and use an analogous function $\eta(\varphi)$ as above.
Because we can in general not guarantee, that $|\varphi|\leq C$ for some suitable and a-priori known $C>0$, we use these non-linear functions for $\rho$ and $\eta$ to guarantee, that both physical values stay positive during the simulation. Note that as long as $|\varphi|\leq 1 + \frac{\rho_1}{\rho_2-\rho_1}$ holds, we have a linear relation between $\varphi$ and $\rho$ and $\eta$.

For the free-energy potential $W$ we choose the double-well potential with a quadratic extension for $|\varphi|>1$, namely
\begin{align*}
    W(\varphi) = \left\{
        \begin{aligned}
            &\frac{1}{4}(1-\varphi^2)^2 &&\text{if } |\varphi|\le 1,\\
            &(|\varphi|-1)^2 &&\text{if }|\varphi| > 1,
        \end{aligned}
    \right.
\end{align*}    
and to use the following convex-concave splitting
\begin{align*}
    W_+(\varphi) = 
    \begin{cases}        
            \frac{1}{4}(\varphi^4 - 1) &\text{if } |\varphi|\le 1,\\
            \frac{1}{2}(3\varphi^2 - 4|\varphi| + 1) &\text{if }|\varphi| > 1,
    \end{cases}
    \qquad W_-(\varphi) = \frac{1}{2}(1-\varphi^2).
\end{align*}  
By extending $W$ quadratically for $|\varphi|>1$, an $L^\infty$-bound on $\varphi$ holds, see \cite{Caffarelli1995129}.
The choice of the free-energy potential also leads to $\widetilde{\Phi}_0(z) = \tanh\left(\frac{z}{\sqrt{2}} \right)$ and consequently to  $c_W = \frac{3}{2\sqrt{2}}$.

Furthermore, we use the following interpolation functions
\begin{align*}
    \alphu{\psi} =\overline{\alpha}_{\bu}  \delta^{-1}\frac{1-\psi}{2},
    \qquad
    \alphi{\psi} =  \frac{1+\psi}{2} +  \underline{\alpha}_\varphi \delta\frac{1-\psi}{2}.
\end{align*}
Finally, the velocity field is initialized with $\bu^0 = 0$ and 
the liquid-gas distribution is initialized with
\begin{align*}
     \varphi^0(x) = \widetilde{\Phi}_0\left(\frac{1}{\varepsilon} \big( 0.25 - \|x - (0,1.875)\|_2 \big) \right),
\end{align*}
which defines a spherical droplet centered at $(0,1.875)$ with radius $0.25$.
We visualize $\varphi_0$ in Figure~\ref{fig:daminit}(b).


\paragraph*{Setup - numerical parameters}\mbox{}\\
We approximate system \eqref{system:weakdiscretepores} using finite-elements on a triangulation of $\Omega$ provided by the software package \texttt{FreeFem++} \cite{H2012JNM}. 
 Our implementation in \texttt{FreeFEM++} is documented in the GitHub repository
 \mbox{\href{https://github.com/jhsimon/ThermoConsistentCHNS}{jhsimon/ThermoConsistentCHNS}.}

We use piecewise linear and globally continuous finite elements for  $\varphi$ and $\mu$ and the stable MINI element for $\bu$ and $p$, i.e. $\bu$ and $p$ are piecewise linear and globally continuous and additionally $\bu$ is enriched by triangle-wise bubble functions.
Due to the Dirichlet boundary data, it is well known that the pressure is only defined up to a constant, i.e. its mean value.
To treat the mean value of $p$ in \texttt{FreeFem++}  we modify the solenoidality condition $\dive{\bu} =0 $ to $\dive{\bu} - \theta p = 0$, with $\theta=10^{-6}$,
i.e. we use a pressure penalty approximation 
as e.g. in \cite[(4.8)]{PS2021JMFM}. This formulation introduces an error to the system, that converges linearly with respect to $\theta$, see \cite[Prop.~4.3.3]{BBF2013}. 

We notice that the choice of the free-energy potential $W$, and its convex-concave splitting leads to a nonlinear system \eqref{system:weakdiscretepores} on every time step.
To resolve this we use Newton's method.

\begin{figure}
     \centering
         \centering
         \includegraphics[trim=0 600 300 150,clip,angle=90,width=0.5\textwidth]{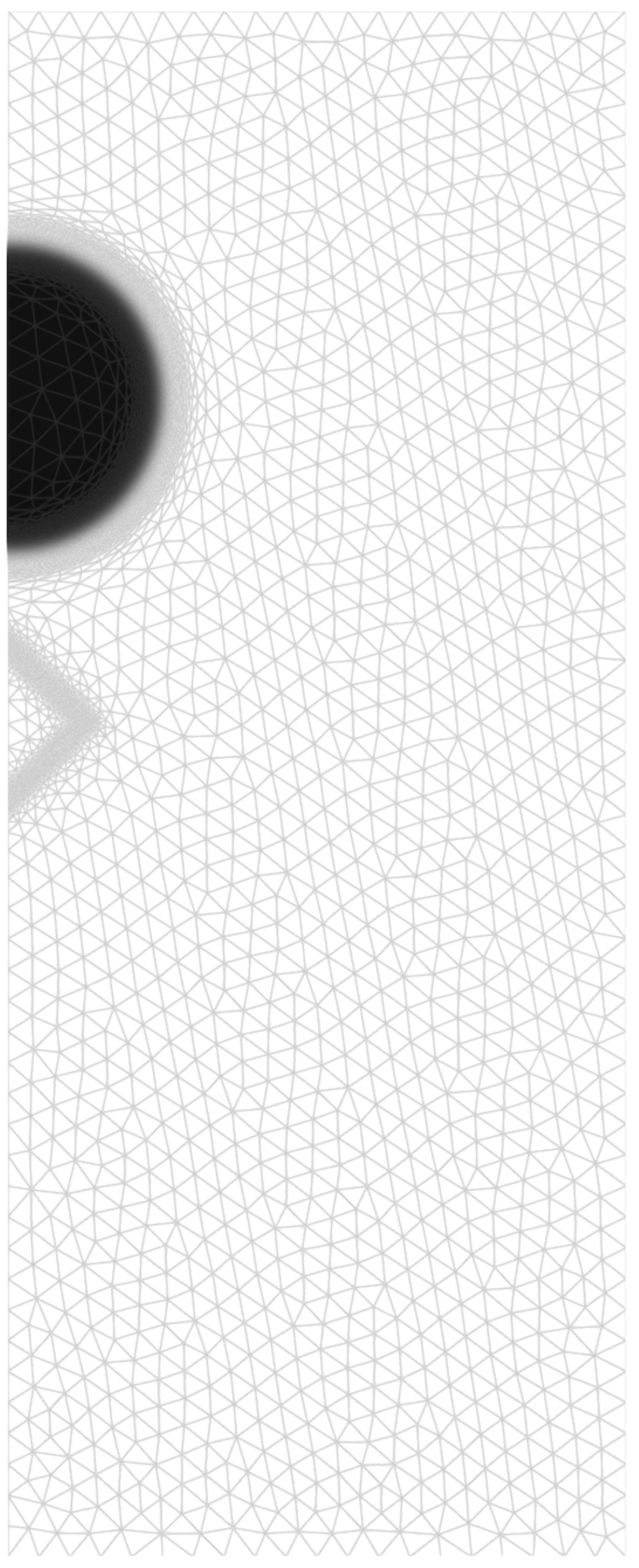}
   \caption{This figure  illustrates the refinement of the triangulation around the interface of the two fluids and around the interface that approximates the boundary of the obstacle. We show only a subset of $\Omega$ for illustration reasons and we rotate the domain by $90^{\degree}$ such that gravity points to the right.
   }
     \label{fig:adapt}
\end{figure}

Both phase-fields $\psi$ and $\varphi$ admit a very localized structure, namely they are (approximately) constant in large regions of $\Omega$ and rapidly change their value over short length scales $\delta$ and $\epsilon$, respectively.
Therefore it is meaningful to use local mesh adaptation throughout the simulation. We use the default adaptivity routine from \texttt{FreeFEM++}. This refines a mesh based on the \textit{approximate} curvature of a given function. 
We use the function $\psi\varphi^m$ and use the adaptation algorithm after every step, with minimum edge size set to $h_{\min} = 0.25\delta$.
The choice of the function and minimum edge size allows for the mesh to be refined on the fluid and obstacle interfaces while taking into account that the simulations we consider may have thinner obstacle interface. This also guarantees sufficient resolution for the interface associated with $\varphi$.
Figure \ref{fig:adapt} illustrates the refined mesh for a simulation.

For the fluid-parameters we take the set-up from the benchmark paper \cite{Hysing_Turek_quantitative_benchmark_computations_of_two_dimensional_bubble_dynamics}, see Table~\ref{tab:params} (left). The further parameters that define the numerical setup are stated in Table~\ref{tab:params} (right).
We simulate the droplet sliding over the obstacle for a sequence of $\delta\to 0$ and also compare this with the result for $\delta=0$, i.e. for the corresponding sharp setting, system~\eqref{system:prime}/\eqref{unifNS-rhoeq}.

\begin{remark}
    As we do not provide a convergence analysis with respect to $\tau$ and $h$ we do not provide parameter studies with respect to $\tau$ and $h$. 
    We assume that the taken values are sufficiently small to not significantly spoil the numerical results with respect to $\delta$.
\end{remark}

\begin{table}[t]
    \centering 
    \begin{tabular}{cccccc}
        \toprule
         $\rho_1$ & $\rho_2$ & $\eta_1$ & $\eta_2$ & $\bg$ & $\sigma_s$ \\
        \midrule
         $1$ &  $1000$ & $0.1$ & $10$ & $[0,-9.8]^\top$ & $24.5$\\
        \bottomrule
    \end{tabular}
    \hspace{1ex}
        \begin{tabular}{cccccccccc}
        \toprule
        $\overline{\alpha}_{\bu}$& $\underline{\alpha}_\varphi$ &$b$ & $\varepsilon$ & $h_{\min}$ & $T$ & $\tau := T/M$ \\
        \midrule
          $10^5$& $10^{-2}$&  $2\times10^{-5}$ & $10^{-2}$ & $0.25\delta$ &  $25$ & $0.01$\\
        \bottomrule
    \end{tabular}
    
    \caption{Fixed parameter values for the numerical setup. The parameters that define the fluid are taken from the benchmark paper
    \cite{Hysing_Turek_quantitative_benchmark_computations_of_two_dimensional_bubble_dynamics}.
    The parameter $\delta$ is varying for the individual tests and is stated in the relating paragraphs.
    }
    \label{tab:params}
\end{table}

\subsection{Numerical tests}

\begin{figure}
     \centering
     \begin{subfigure}{0.2\textwidth}
         \centering
         \includegraphics[width=\textwidth]{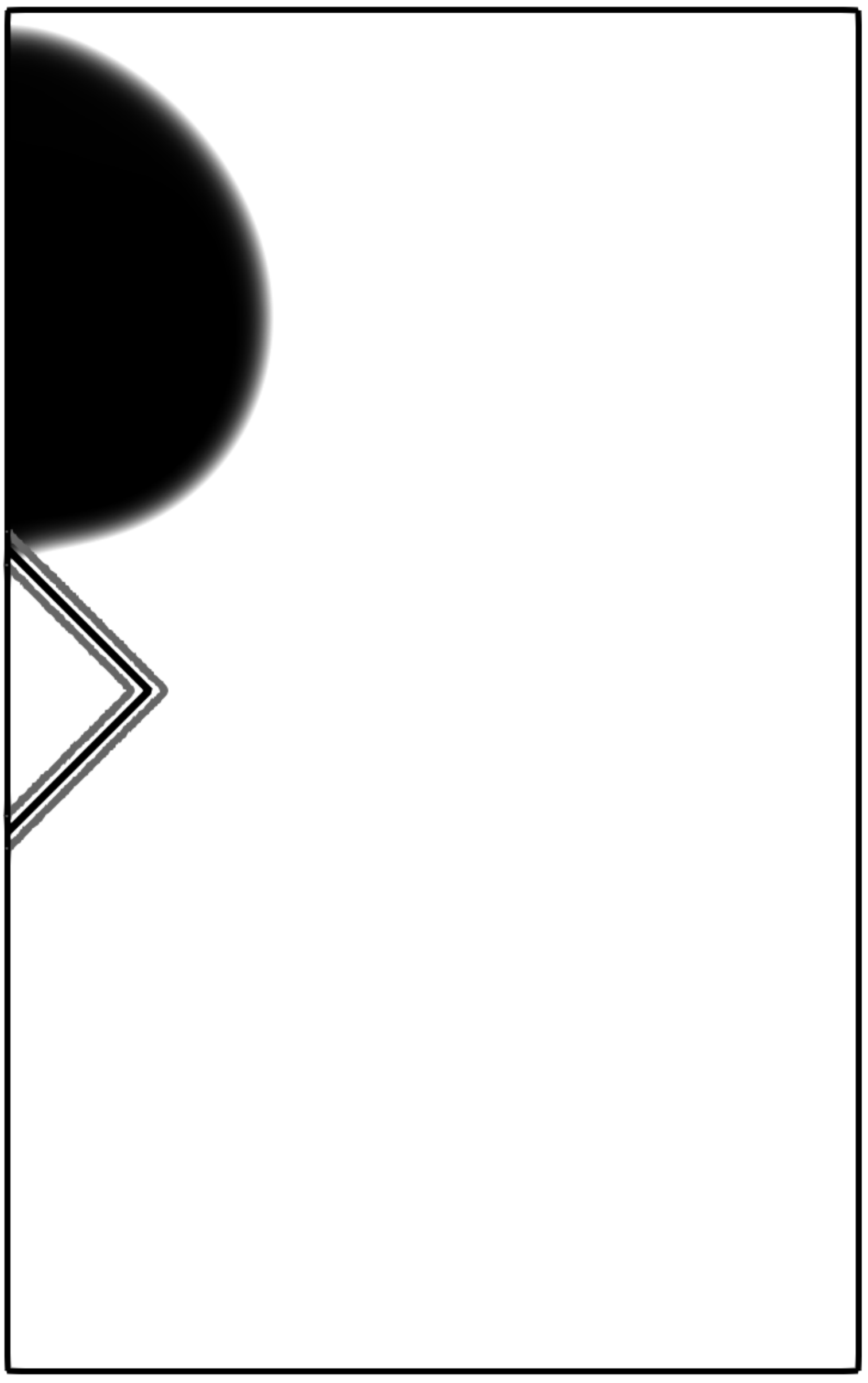}
         \caption*{(a)}
     \end{subfigure}
     \begin{subfigure}{0.2\textwidth}
         \centering
         \includegraphics[width=\textwidth]{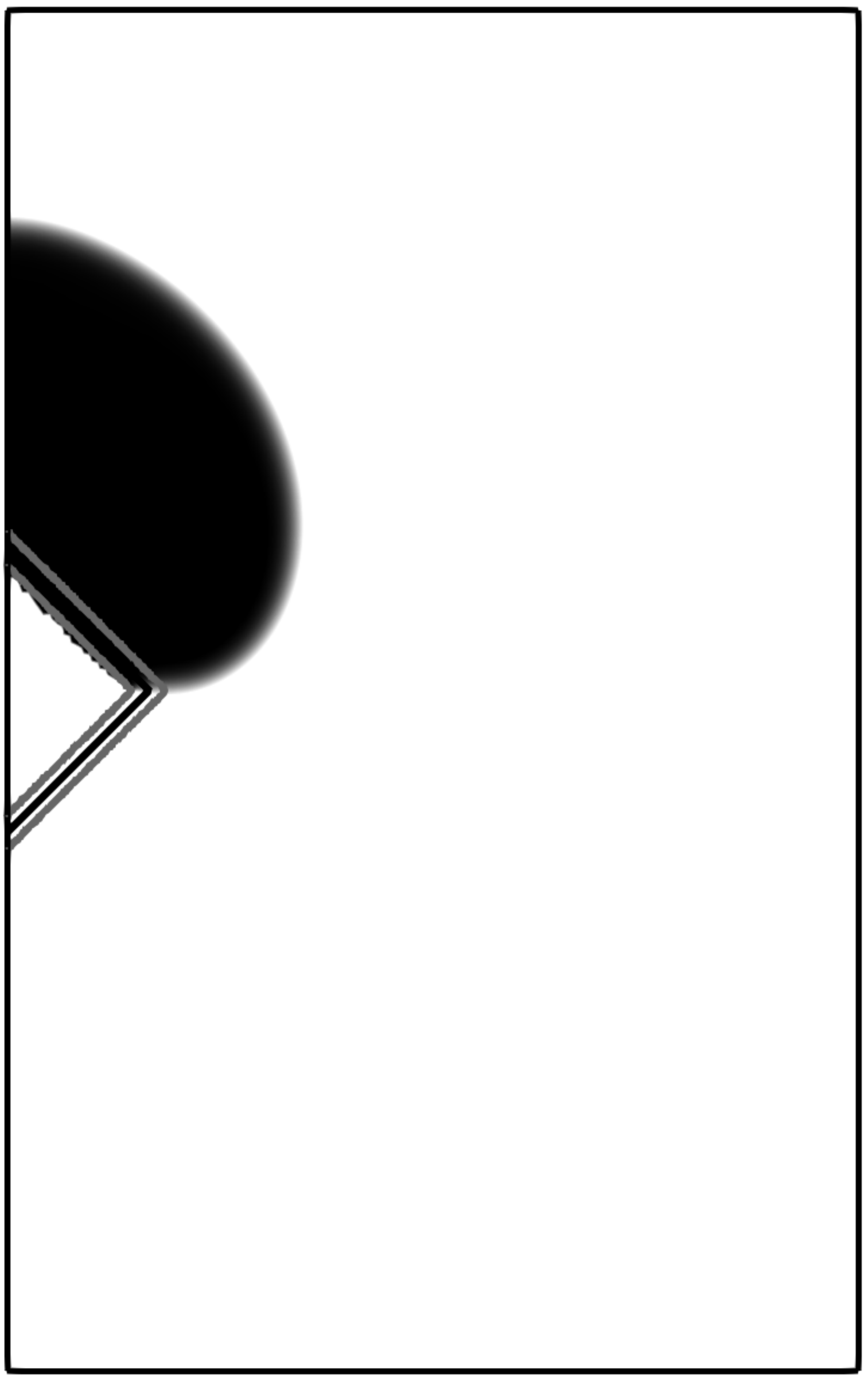}
         \caption*{(b)}
     \end{subfigure}
     \begin{subfigure}{0.2\textwidth}
         \centering
         \includegraphics[width=\textwidth]{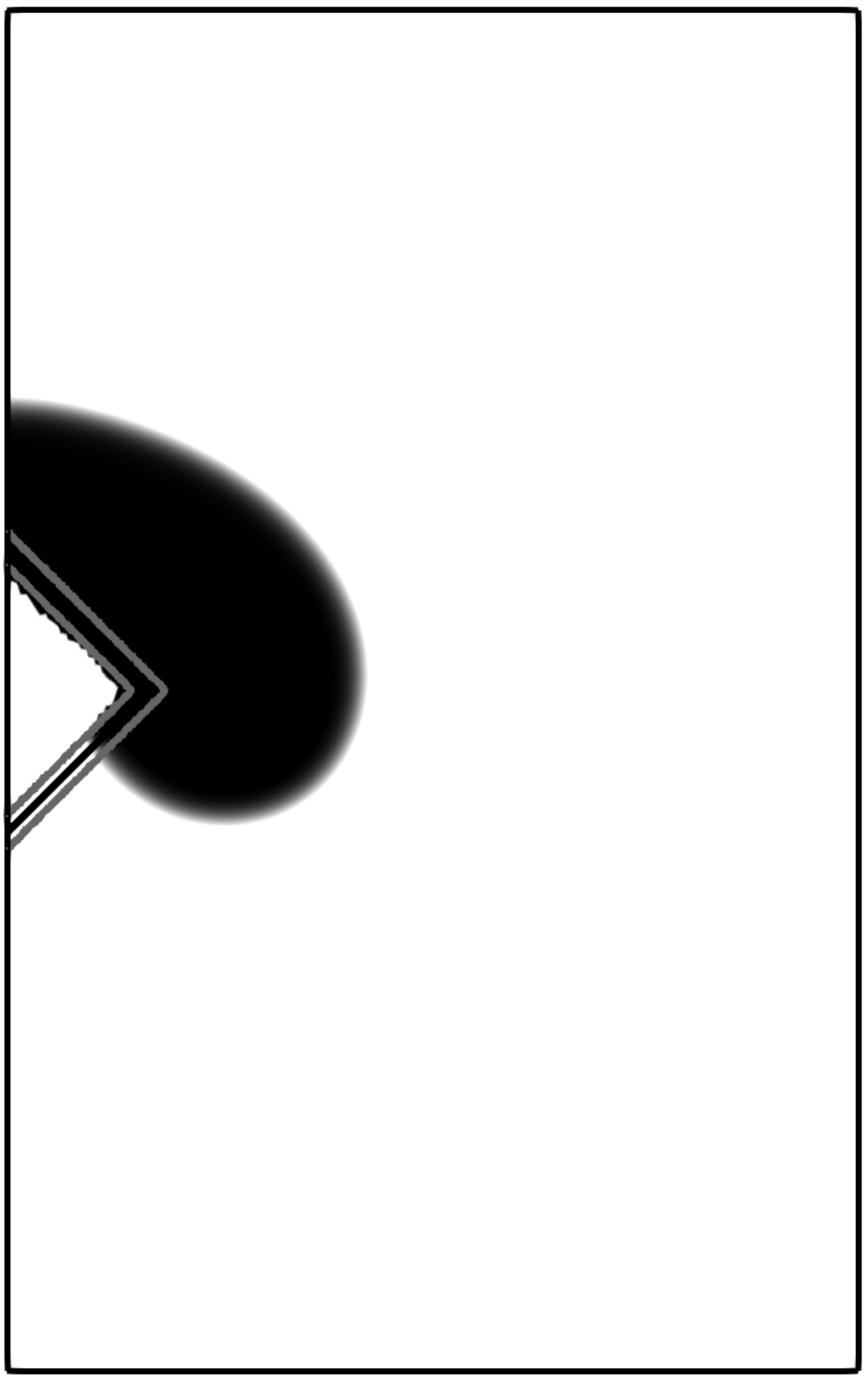}
         \caption*{(c)}
     \end{subfigure}
     \begin{subfigure}{0.2\textwidth}
         \centering
         \includegraphics[width=\textwidth]{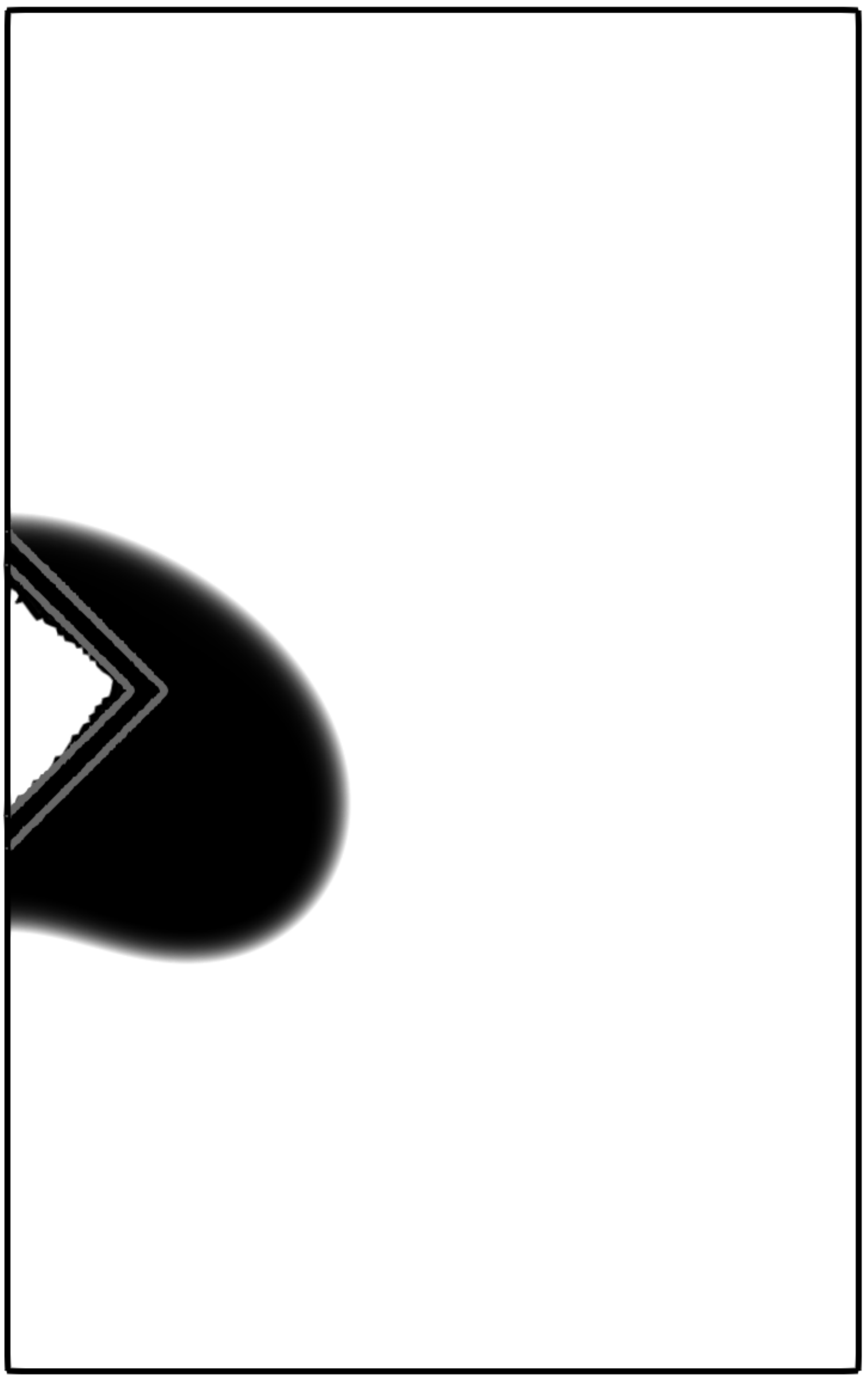}
         \caption*{(d)}
     \end{subfigure}
     \begin{subfigure}{0.2\textwidth}
         \centering
         \includegraphics[width=\textwidth]{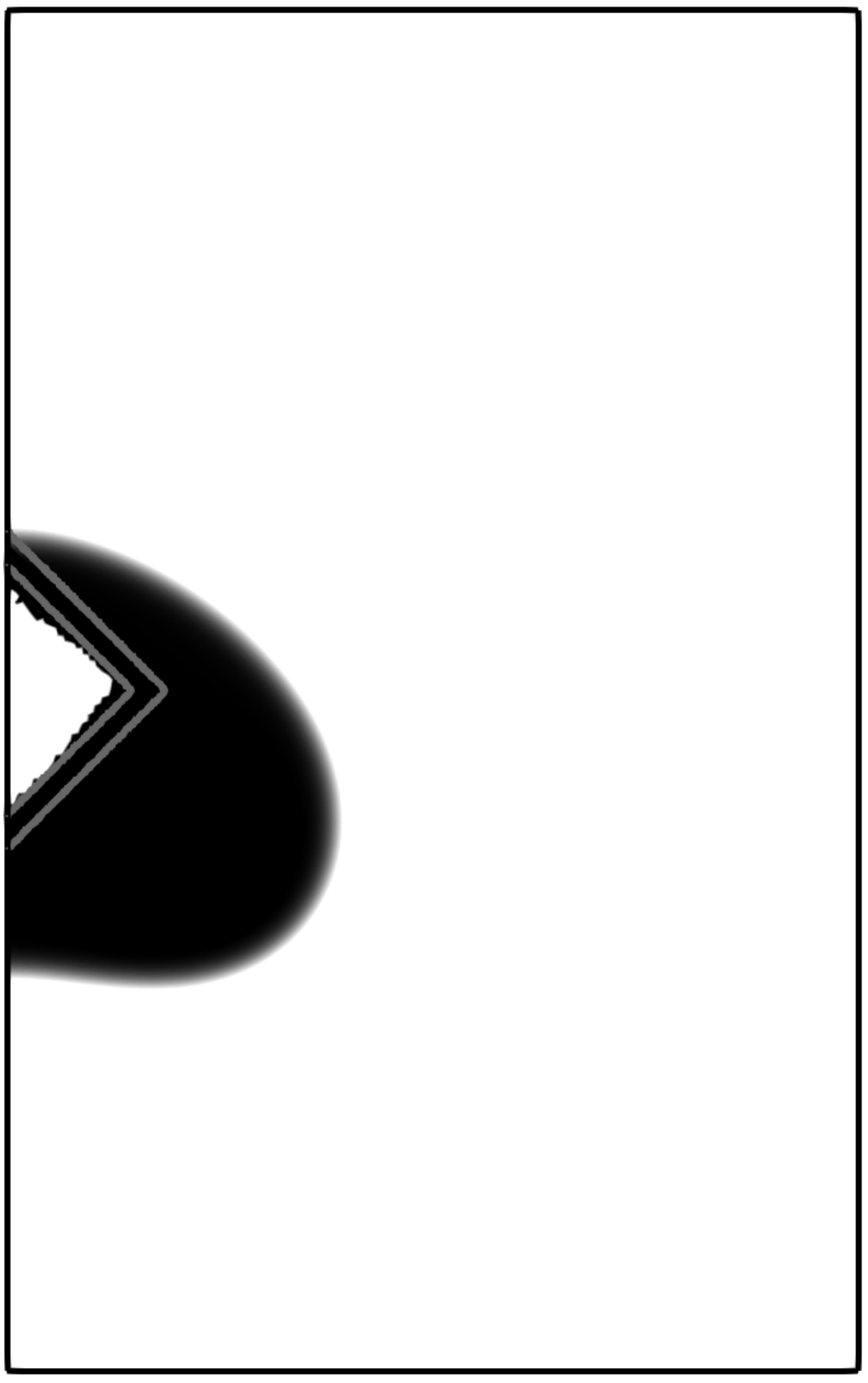}
         \caption*{(e)}
     \end{subfigure}
     \begin{subfigure}{0.2\textwidth}
         \centering
         \includegraphics[width=\textwidth]{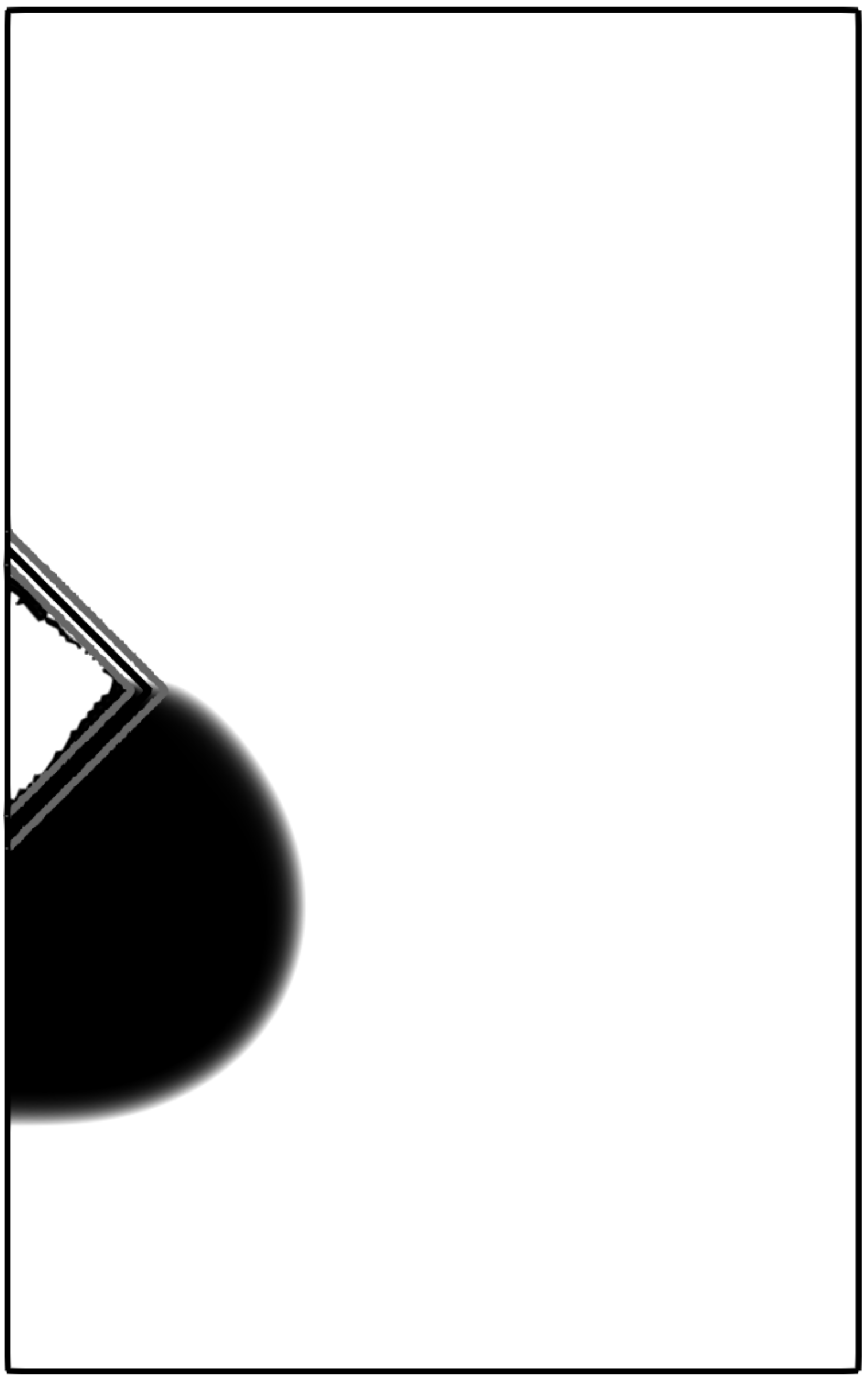}
         \caption*{(f)}
     \end{subfigure}
     \begin{subfigure}{0.2\textwidth}
         \centering
         \includegraphics[width=\textwidth]{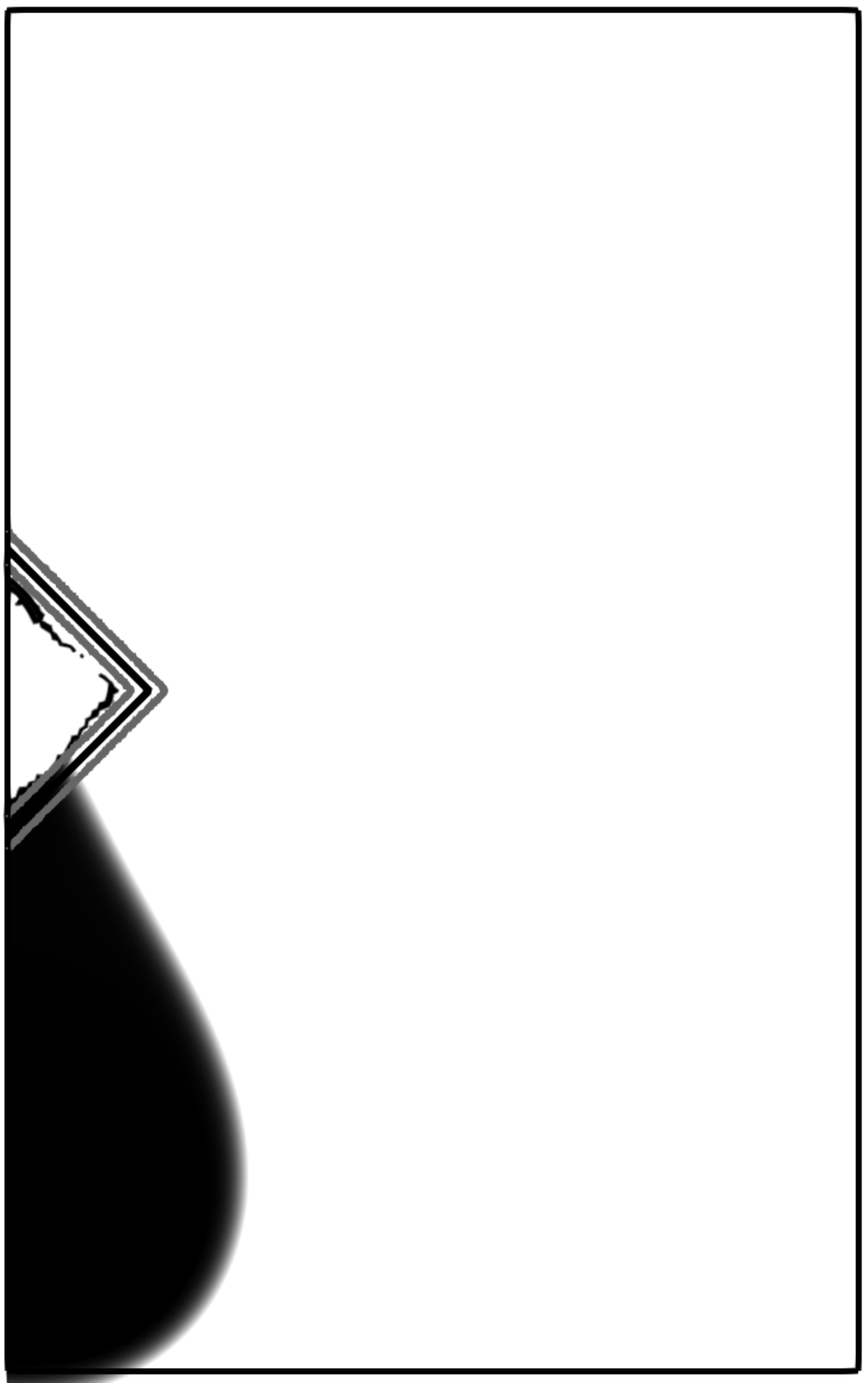}
         \caption*{(g)}
     \end{subfigure}
     \begin{subfigure}{0.2\textwidth}
         \centering
         \includegraphics[width=\textwidth]{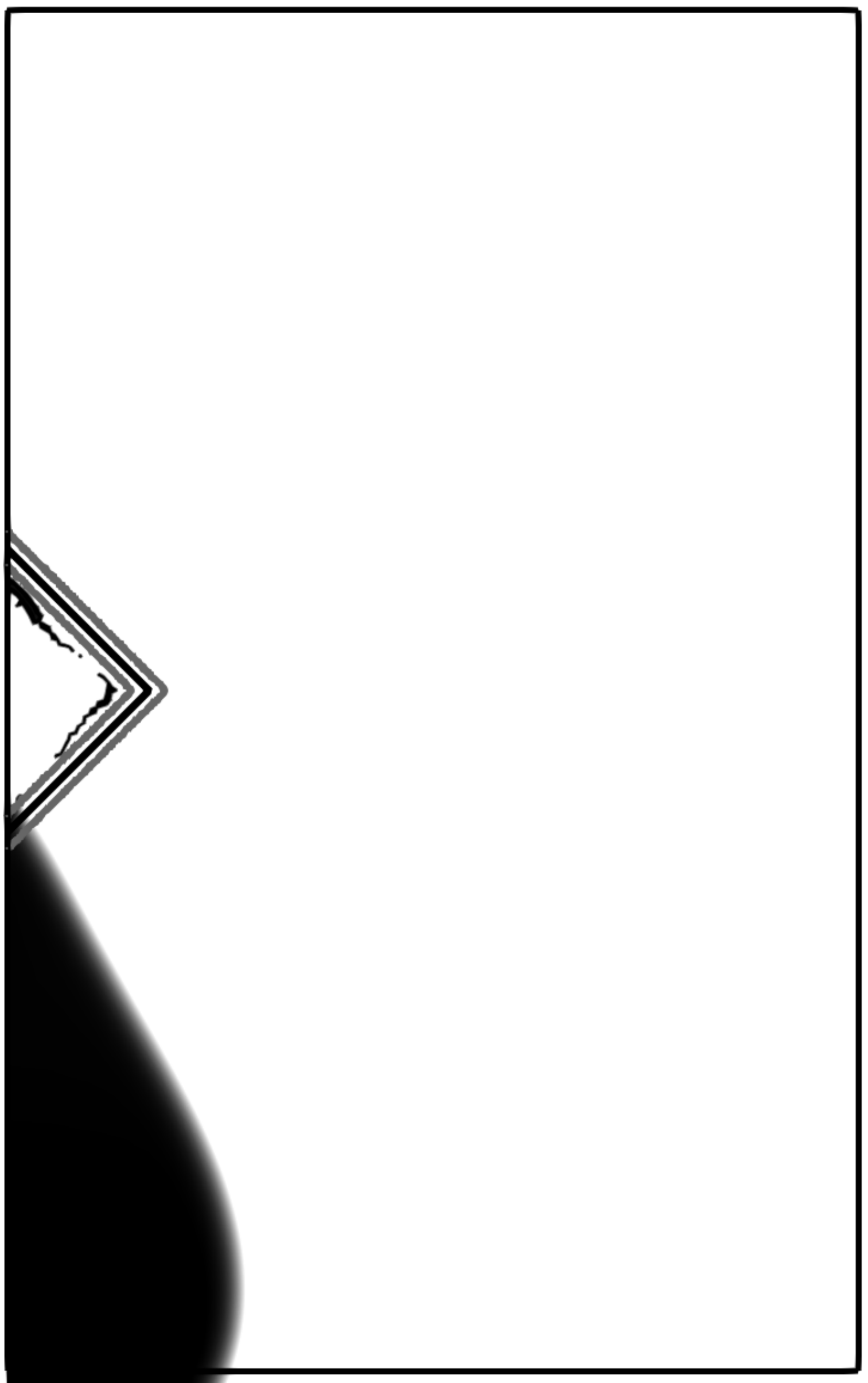}
         \caption*{(h)}
     \end{subfigure}
   \caption{This figure shows the evolution of the fluid concentration $\varphi$ at times (a) $t = 2.65$, (b) $t = 5.15$, (c) $t = 9.50$, (d) $t = 12.00$, (e) $t = 12.80$, (f) $t = 14.40 $, (g) $t = 18.50$ and (h) $t = 20.26$. Together with the level-sets $\{x\in\Omega: \psi(x) =0 \}$, $\{x\in\Omega: \psi(x) =-0.99 \}$ and $\{x\in\Omega: \psi(x) = 0.99 \}$, the images show a zoomed-in view---over the subdomain $(0,0.75)\times(0.75,1.95)$---of the evolution of the droplet and its behavior as it slides down the triangular obstacle. In particular we observe some spill-over, i.e. some part of the liquid droplet penetrates the obstacle.}
     \label{fig:damflow}
\end{figure}

\paragraph{Droplet sliding over a triangular structure.} \mbox{}\\
We fix $\delta = 10^{-2}$ and simulate the provided setup. In Figure~\ref{fig:damflow} we show snapshots of the simulation.
We represent the obstacle $\psi$ by the three isolines $\psi \equiv \pm 0.99$ and $\psi\equiv 0$ and show the phase field $\varphi$ in grayscale.

The droplet slides as expected over the obstacle. Especially, we observe the expected protrusion of the leading part of the droplet when it crosses the kink of the obstacle, see Figure~\ref{fig:damflow} (c).  Such phenomena is physically relevant as illustrated in \cite[Fig.~9]{BJKR2019CET}, which they call \textit{nearly pinning}.

We also observe that, since the obstacle is approximated by a porous medium, the fluid is allowed to flow through the interface $\{x\in\Omega : -1 < \psi(x)< 1\}$ and penetrate the obstacle. Therefore, we observe a \textit{spill-over} in such region, see Figures \ref{fig:damflow} (c-h). 

\begin{figure}
		\centering

        {
        \pgfplotsset{
            every axis/.style={
            width=10cm,
            height=2.5cm
            }
            }
        \definecolor{mycolor1}{rgb}{0.0, 0.0, 0.0}%
\definecolor{mycolor2}{rgb}{0.0039216, 0.2352941, 0.8666667}%
\definecolor{mycolor3}{rgb}{0.5843137, 0.0039216, 0.5607843}%
\definecolor{mycolor4}{rgb}{0.6431373, 0.0000000, 0.3333333}%
\definecolor{mycolor5}{rgb}{0.5176471, 0.2470588, 0.0000000}%
\definecolor{mycolor6}{rgb}{0.3882353, 0.3176471, 0.0000000}%
\definecolor{mycolor7}{rgb}{0.0039216, 0.3764706, 0.1764706}
\begin{tikzpicture}

\pgfmathsetmacro{\DataSpacing}{10}

\begin{axis}[%
scale only axis,   
xmin=-0.1,
xmax=25.2,
xlabel = {$t$},
ymin=-2.3045,
ymax=-2.2835,
ylabel = {$\mathcal{M}$},
ytick={-2.312, -2.308, -2.304, -2.300, -2.296, -2.292, -2.288, -2.284},
yticklabels ={$-2.312$,$-2.308$,$-2.304$,$-2.300$,$-2.296$,$-2.292$,$-2.288$, $-2.284$},
axis background/.style={fill=white},
axis x line*=bottom,
axis y line*=left,
grid = major,
axis lines=box,
]
\addplot [color=mycolor1, line width=1pt]
table[row sep=crcr]{%
0.01	-2.28749\\
0.51	-2.28759\\
1.01	-2.28765\\
1.51	-2.28767\\
2.01	-2.28769\\
2.51	-2.28773\\
3.01	-2.28902\\
3.51	-2.2905\\
4.01	-2.2918\\
4.51	-2.29276\\
5.01	-2.29308\\
5.51	-2.29327\\
6.01	-2.29331\\
6.51	-2.29339\\
7.01	-2.29346\\
7.51	-2.29377\\
8.01	-2.29402\\
8.51	-2.29442\\
9.01	-2.29508\\
9.51	-2.29566\\
10.01	-2.29643\\
10.51	-2.29733\\
11.01	-2.2981\\
11.51	-2.29887\\
12.01	-2.29961\\
12.51	-2.30013\\
13.01	-2.30036\\
13.51	-2.30009\\
14.01	-2.29894\\
14.51	-2.29697\\
15.01	-2.29642\\
15.51	-2.2961\\
16.01	-2.29577\\
16.51	-2.29543\\
17.01	-2.29497\\
17.51	-2.29453\\
18.01	-2.29396\\
18.51	-2.29333\\
19.01	-2.29278\\
19.51	-2.29209\\
20.01	-2.29153\\
20.51	-2.29105\\
21.01	-2.29059\\
21.51	-2.2906\\
22.01	-2.2906\\
22.51	-2.2906\\
23.01	-2.29062\\
23.51	-2.29061\\
24.01	-2.29061\\
24.51	-2.29065\\
25.01	-2.29063\\
};

\addplot[only marks, mark=square*, mark options={fill=mycolor1, draw=mycolor1}, mark size=3pt, draw=mycolor2, forget plot] table[row sep=crcr]{%
x	y\\
2.51	-2.28773\\
};
\addplot[only marks, mark=triangle*, mark options={fill=mycolor1, draw=mycolor1}, mark size=4pt, draw=mycolor3, forget plot] table[row sep=crcr]{%
x	y\\
13.01	-2.30036\\
};
\addplot[only marks, mark=*, mark options={fill=mycolor1, draw=mycolor1}, mark size=3pt, draw=mycolor4, forget plot] table[row sep=crcr]{%
x	y\\
21.01	-2.29059\\
};

\draw[dashed, mycolor1] (0, -2.28773) -- (25, -2.28773);

\draw[dashed, mycolor1] (0, -2.30036) -- (25, -2.30036);

\draw[dashed, mycolor1] (0, -2.29059) -- (25, -2.29059);

\draw[dashed, mycolor1] (2.51, -2) -- (2.51, -3);

\draw[dashed, mycolor1] (13.01, -2) -- (13.01, -3);

\draw[dashed, mycolor1] (21.01, -2) -- (21.01, -3);

\end{axis}

\end{tikzpicture}%
        }
		\caption{ This figure shows the evolution of $\mathcal M = \int_\Omega \frac{1+\psi}{2}\varphi \du x$ over the simulation time. The square corresponds to the situation depicted in Figure~\ref{fig:damflow}(a), the triangle marks approximately situation Figure~\ref{fig:damflow}(e) and the circle marks Figure~\ref{fig:damflow} (h).
        We observe that $\mathcal M$ stays constant while the droplet approaches the obstacle and decreases while the droplet overflows it. After leaving the obstacle the initial value of $\mathcal M$ is not fully recovered, i.e. some liquid is left inside the obstacle.
        }
		\label{fig:mass-trend}
	\end{figure}
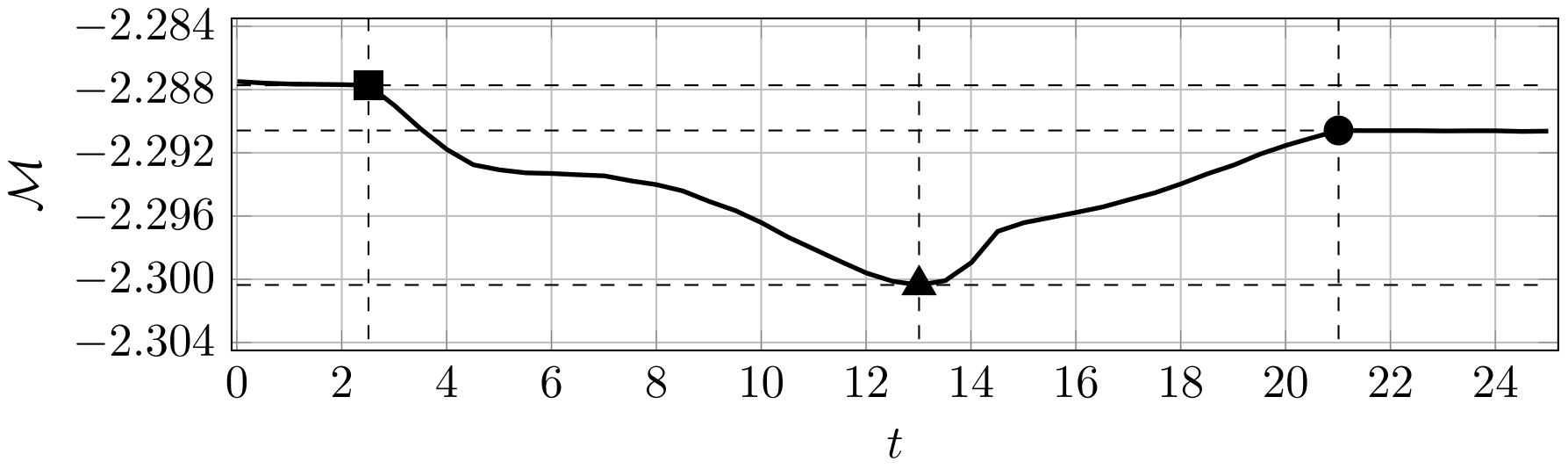

For a further investigation of this spill-over effect we turn our attention to the quantity 
\begin{align*}
    \mathcal{M} = \int_\Omega \frac{\psi+1}{2}\varphi \du x,
\end{align*}
which we have introduced in \eqref{mass} and which measures the amount of $\varphi$ inside the fluid domain, given by $\psi=1$.
Since the droplet is given by $\varphi=+1$ and inside the obstacle we initially have $\varphi^0 = -1$, i.e. gas, the spill-over effect should decrease this value.
In Figure~\ref{fig:mass-trend} we show the evolution of $\mathcal M$ over the simulation time. We observe conservation of $\mathcal M$ while the droplet approaches the obstacle, then decreases 
while the droplet overflows the obstacle and gradually increases again, but does not reach its starting value, because some portion of the droplet is left inside the obstacle.

On the other hand, this spill-over effect is localized close to the interface of $\psi$.
This fits the formal matched asymptotics, where we show that to leading order, $\partial_t \varphi_0=0$, i.e. to first order $\varphi$ is stationary 
inside the obstacle, see \eqref{eq:asy:phi-leading-B}.
On the other hand we observe that this is not fully reversed as the droplet leaves the obstacle.
We investigate this in the next paragraph in more detail.


\paragraph{An  illustration of the fluid inside the obstacle.} \mbox{}\\
From a conceptual point of view, the inner part of the obstacle is not of interest for the simulation while the modeling also assign values to $\varphi$ inside the obstacle. Especially and as seen in the previous paragraph, we might observe a spill-over effect, when liquid sliding over an obstacle partially penetrates the obstacle. 

\begin{figure}
     \centering
     \begin{subfigure}{0.2\textwidth}
         \centering
         \includegraphics[width=\textwidth]{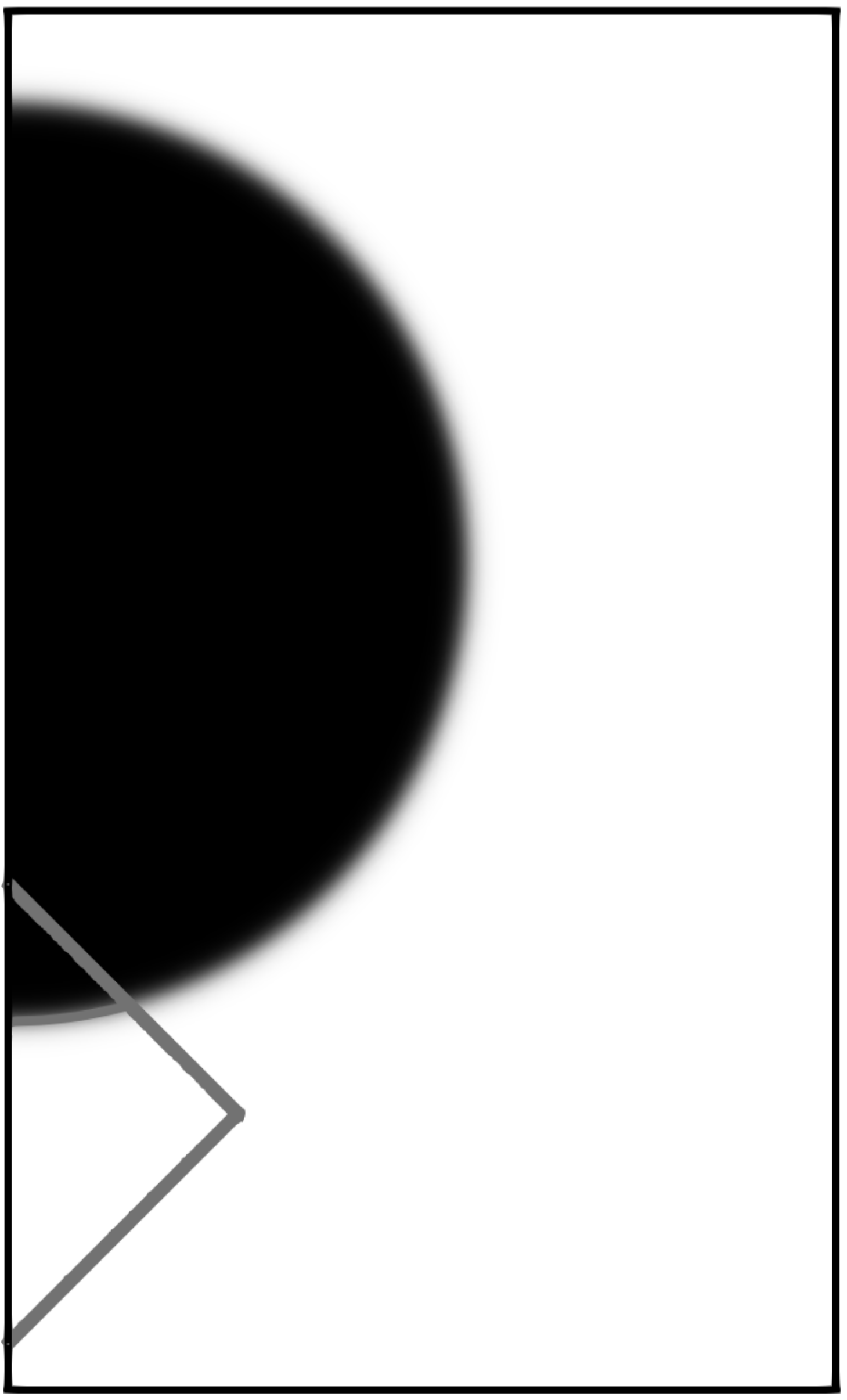}
         \caption*{(a)}
     \end{subfigure}
     \begin{subfigure}{0.2\textwidth}
         \centering
         \includegraphics[width=\textwidth]{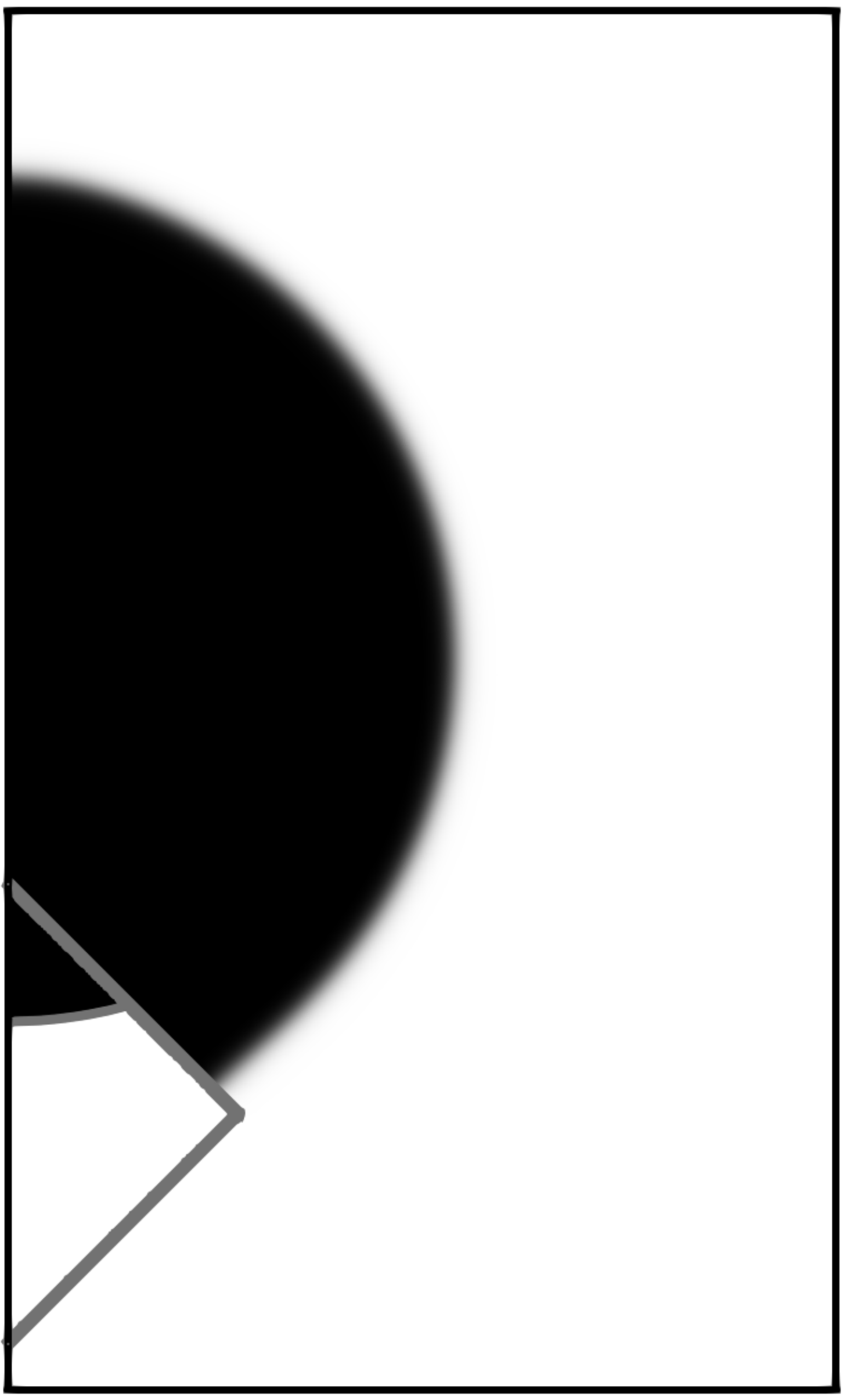}
         \caption*{(b)}
     \end{subfigure}
     \begin{subfigure}{0.2\textwidth}
         \centering
         \includegraphics[width=\textwidth]{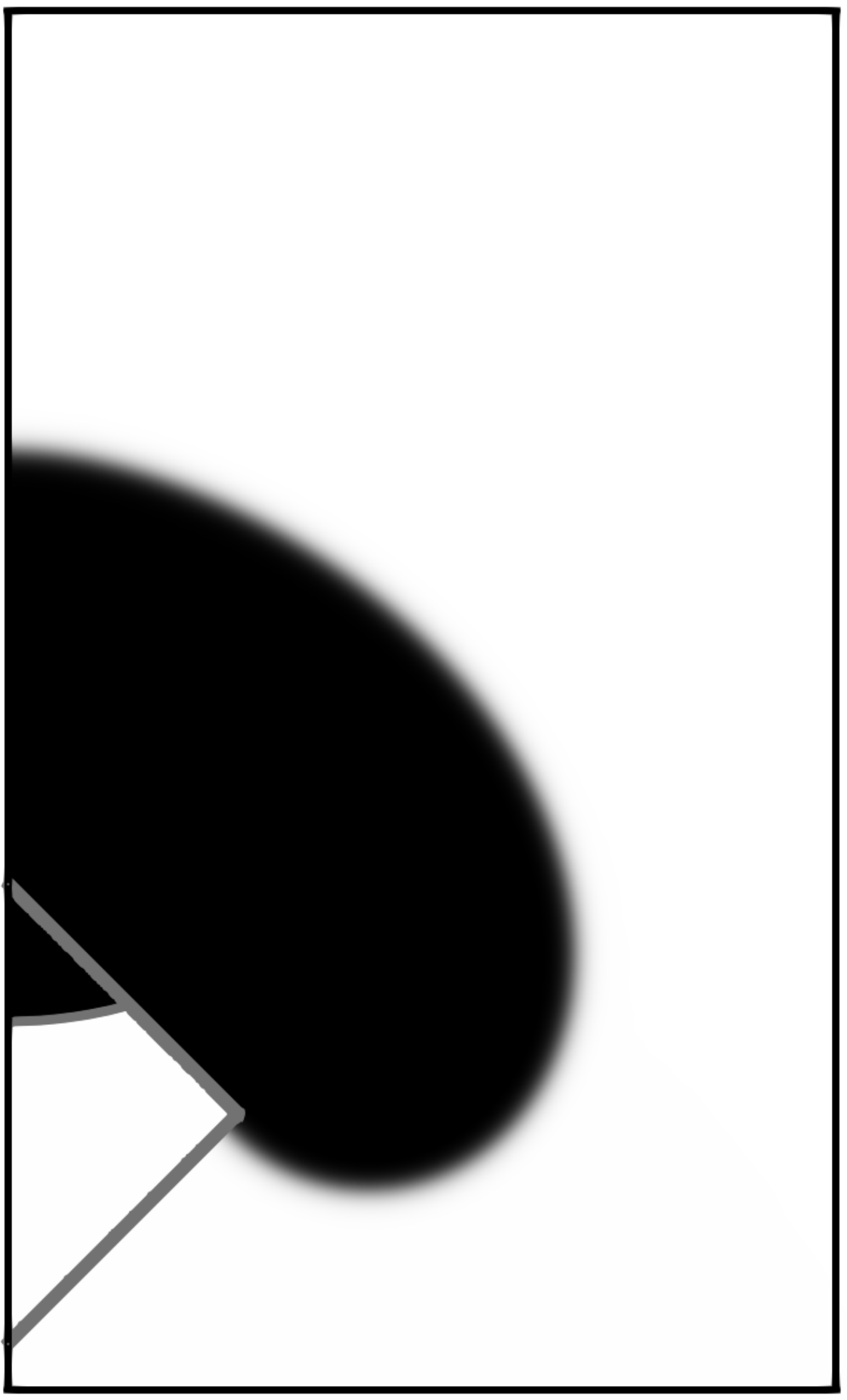}
         \caption*{(c)}
     \end{subfigure}
     \begin{subfigure}{0.2\textwidth}
         \centering
         \includegraphics[width=\textwidth]{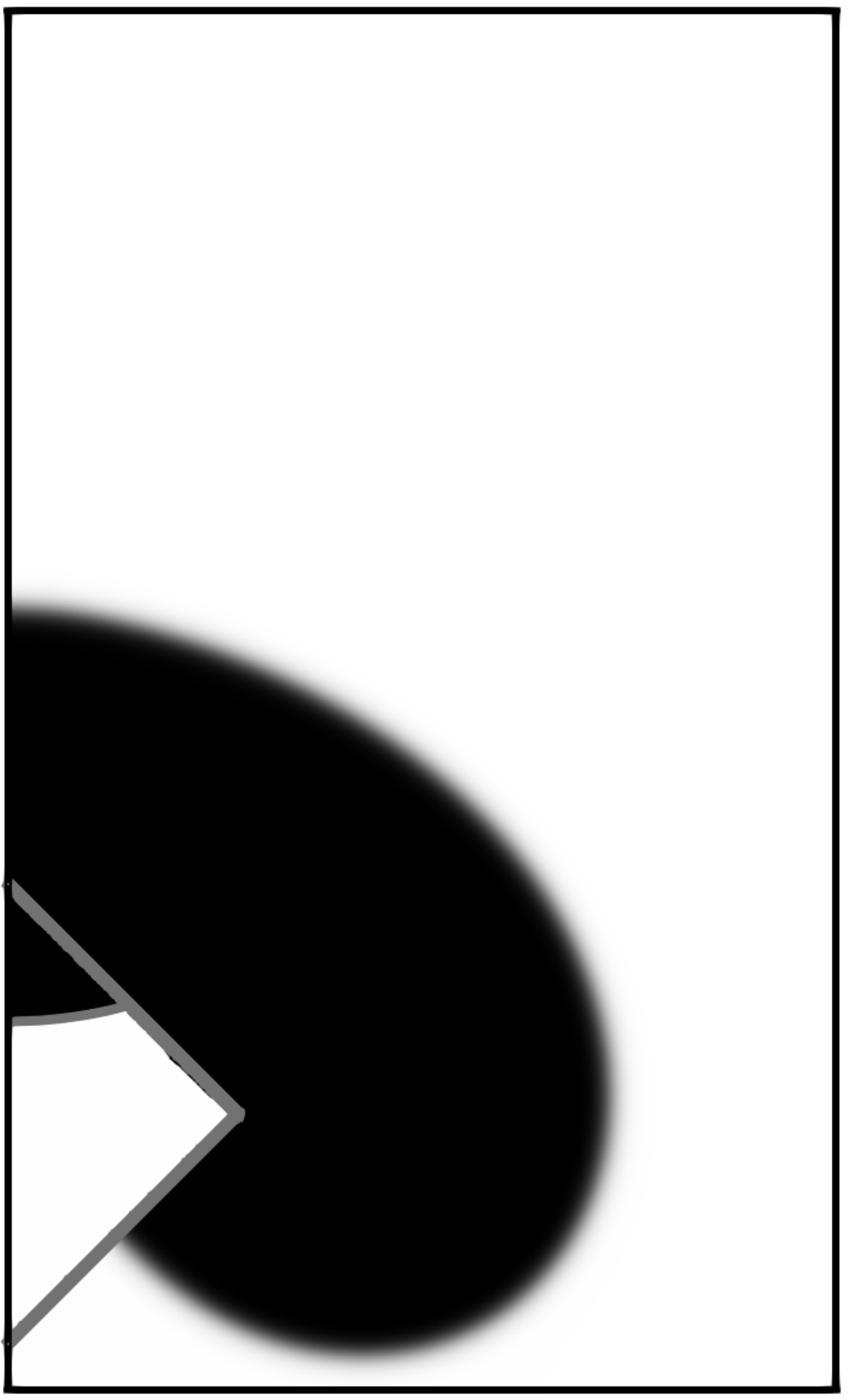}
         \caption*{(d)}
     \end{subfigure}
   \caption{This figure illustrates the example where we initialize the droplet and the obstacle to initially intersect. In particular, it shows the evolution of $\varphi$ at times (a) $t = 0.00$, (b) $t = 2.50$, (c) $t = 15.00$ and (d) $t = 25.00$. 
   Together with the level-set $\{x\in\Omega: \psi(x) =0 \}$ and the sector of a circle that coincides with the initial zero-level set of $\varphi$ inside the obstacle, the images show a zoomed-in view of the evolution of the droplet and its behavior as it slides down the triangular obstacle. In agreement with \eqref{UinB} and \eqref{eq:asy:phi-leading-B}, we observe no qualitative change in the phase field $\varphi$ inside the obstacle as the remainder of the droplet moves over the triangular obstacle. }
     \label{fig:intersected}
\end{figure}

We note that formal matched asymptotics, \eqref{UinB} and \eqref{eq:asy:phi-leading-B}, imply that  to first order there should be no motion and change in $\varphi$ inside the obstacle. 
To illustrate this behavior we modify the set-up in such a way that the obstacle and the initial droplet have an intersection, i.e., the set $\{x\in \Omega: \psi(x) = -1\text{ and }\varphi(x) = 1 \}$ is not empty. For the sake of this illustration we use the following initial set-up for the fluid
\begin{align*}
    \varphi^0(x) = \widetilde{\Phi}_0\left(\frac{1}{\varepsilon} \big( 0.25 - \|x - (0,1.7)\|_2 \big) \right),
\end{align*}
and the following set-up for the obstacle 
\begin{align*}
    \psi(x) = \widetilde{\Psi}_0\left( \frac{1}{\delta}\big(\|x - (0,1.4)\|_1-0.125 \big) \right).
\end{align*}
We  use an obstacle interface thickness $\delta = 1.2\times 10^{-3}$, while the other parameters remain the same. In Figure~\ref{fig:intersected}, we illustrate the evolution of the droplet.  
As a point of reference we also plot the zero-level set of the obstacle and a sector of a circle that coincides with the initial zero-level set of $\varphi$ inside the obstacle in gray. 
Figure \ref{fig:intersected} (a) shows the the initial set-up of the droplet and the zero-level set of the obstacle, while Figure \ref{fig:intersected} (b)--(d) shows the evolution of the droplet as time passes. We observe that the part of the phase-field $\varphi$ inside the obstacle remains fixed while the rest is allowed to flow, which is in total agreement with the equations \eqref{UinB} and \eqref{eq:asy:phi-leading-B}, that state that to first order $\varphi$ is stationary 
inside the obstacle.

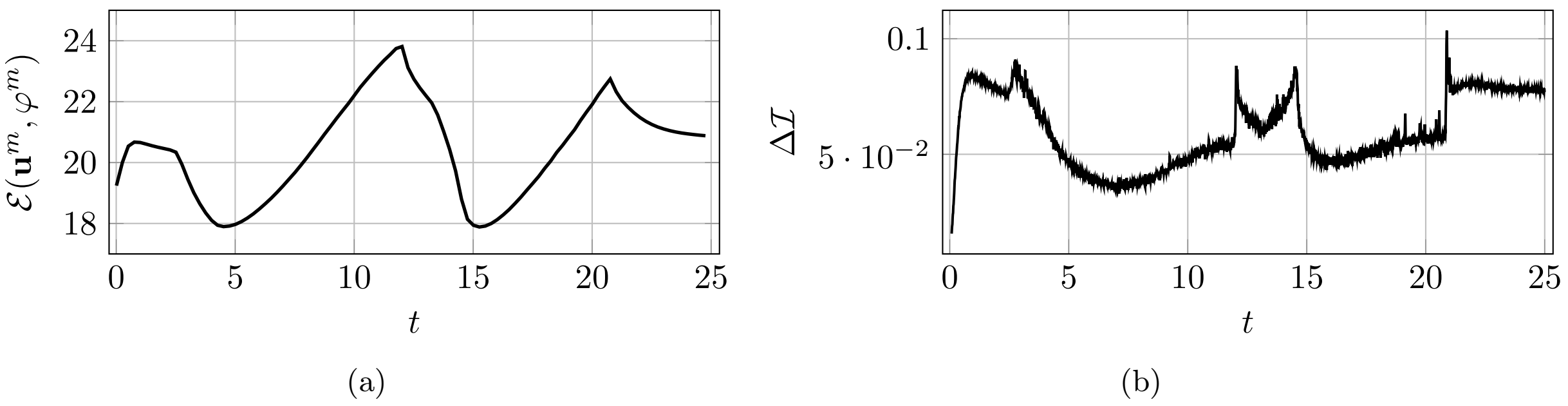
\begin{figure}
\begingroup
\pgfplotsset{    
    every axis/.style={
        width=0.80\textwidth,
        height=2.5cm,
        },
    }
    
    \centering
   
    \begin{subfigure}{0.49\textwidth}
         \centering
         \definecolor{mycolor1}{rgb}{0.0, 0.0, 0.0}%
\definecolor{mycolor2}{rgb}{0.2,0.2,0.2}%
\definecolor{mycolor3}{rgb}{0.4,0.4,0.4}%
\definecolor{mycolor4}{rgb}{0.6,0.6,0.6}%
\definecolor{mycolor5}{rgb}{0.7,0.7,0.7}%
\definecolor{mycolor6}{rgb}{0.8,0.8,0.8}%
\begin{tikzpicture}

\pgfmathsetmacro{\DataSpacing}{10}

\begin{axis}[%
scale only axis,   
xmin=-0.3,
xmax=25.35,
xlabel = {$t$},
ymin=17,
ymax=25,
ylabel = {$\mathcal{E}(\bu^m,\varphi^m)$},
axis background/.style={fill=white},
axis x line*=bottom,
axis y line*=left,
grid = major,
axis lines=box,
legend style={at={(0.75,0.9)}, anchor=center, legend cell align=left, align=left, draw=white!15!black},
]
\addplot [color=mycolor1, line width=1pt]
  table[row sep=crcr]{%
0.01	19.2427\\
0.26	20.0049\\
0.51	20.538\\
0.76	20.6723\\
1.01	20.6583\\
1.26	20.606\\
1.51	20.5511\\
1.76	20.5027\\
2.01	20.4599\\
2.26	20.4164\\
2.51	20.34\\
2.76	19.9413\\
3.01	19.4482\\
3.26	19.0097\\
3.51	18.6504\\
3.76	18.3464\\
4.01	18.1014\\
4.26	17.941\\
4.51	17.8963\\
4.76	17.9158\\
5.01	17.9694\\
5.26	18.06\\
5.51	18.1765\\
5.76	18.3144\\
6.01	18.4731\\
6.26	18.6445\\
6.51	18.8263\\
6.76	19.0244\\
7.01	19.2333\\
7.26	19.4533\\
7.51	19.6724\\
7.76	19.9122\\
8.01	20.1573\\
8.26	20.4166\\
8.51	20.6732\\
8.76	20.9383\\
9.01	21.1917\\
9.26	21.4505\\
9.51	21.6988\\
9.76	21.9434\\
10.01	22.2082\\
10.26	22.475\\
10.51	22.7069\\
10.76	22.9341\\
11.01	23.152\\
11.26	23.3596\\
11.51	23.5447\\
11.76	23.7436\\
12.01	23.8106\\
12.26	23.1148\\
12.51	22.7385\\
12.76	22.4519\\
13.01	22.2017\\
13.26	21.9595\\
13.51	21.5625\\
13.76	21.0211\\
14.01	20.4314\\
14.26	19.729\\
14.51	18.7971\\
14.76	18.1385\\
15.01	17.941\\
15.26	17.8847\\
15.51	17.9148\\
15.76	17.999\\
16.01	18.1251\\
16.26	18.2806\\
16.51	18.4566\\
16.76	18.6576\\
17.01	18.8693\\
17.26	19.099\\
17.51	19.3233\\
17.76	19.5515\\
18.01	19.8133\\
18.26	20.0519\\
18.51	20.3332\\
18.76	20.5692\\
19.01	20.8292\\
19.26	21.0794\\
19.51	21.3808\\
19.76	21.6573\\
20.01	21.9201\\
20.26	22.2226\\
20.51	22.4882\\
20.76	22.7444\\
21.01	22.3263\\
21.26	22.0274\\
21.51	21.8152\\
21.76	21.6329\\
22.01	21.479\\
22.26	21.3529\\
22.51	21.2502\\
22.76	21.1667\\
23.01	21.1017\\
23.26	21.0483\\
23.51	21.0048\\
23.76	20.9703\\
24.01	20.9429\\
24.26	20.9196\\
24.51	20.8985\\
24.76	20.8808\\
};

\end{axis}

\end{tikzpicture}%
         \caption*{(a)}
    \end{subfigure}
    \begin{subfigure}{0.49\textwidth}
         \centering
         \includestandalone{Figures/tikz/energy_diff}
         \caption*{(b)}
    \end{subfigure}
\endgroup
    \caption{ The figure shows (a) the total energy $\mathcal{E}(\bu^m,\varphi^m)$ and (b) the 
    gap $\Delta \mathcal I$ in the energy inequality \eqref{energy:totaldisc} over time. See \eqref{eq:num:energygap} for the definition of $\Delta \mathcal I$.
    Upon cross examination against Figure~\ref{fig:damflow}, Figure~\ref{fig:energy-ineq}~(a) illustrates the relationship between the energy and the behavior of the bubble with respect to its interaction with the obstacle, i.e., the time the bubble initially touches the obstacle, when the bubble touches the tip of the triangle when it starts to leave the obstacle and when the bubble fully leaves the vicinity of the obstacle. Figure~\ref{fig:energy-ineq}~(b) shows that the energy inequality \eqref{energy:totaldisc} remains valid even at the fully discretized setting and despite the mesh refinements done at each time step.
    }
    \label{fig:energy-ineq}
    
\end{figure}

\paragraph{Validation of the energy inequality.} \mbox{}\\
We validate that for our numerical simulation the numerical solution behaves, despite the mesh refinements at each time step,  in accordance to the pointwise-in-time energy inequality \eqref{energy:totaldisc}, see Figure \ref{fig:energy-ineq}.

In Figure \ref{fig:energy-ineq} (a) we show the plot of the total energy $\mathcal{E}(\bu^m,\varphi^m)$ at each iteration $m = 1,\ldots,M$. We observe that the fluctuations in the total energy depict several important epochs at the evolution of the flow. In particular, the fluctuations indicate how the bubble interacts with the obstacle:
\begin{enumerate}[label=\roman*)]
    \item  the decrease around $t = 2.50$ is the point when the bubble first interacts with the obstacle;
    \item  the increase at $t \approx 5.00$ shows the instance when the bubble reaches the tip of the obstacle; 
    \item the change from increasing to decreasing at $t\approx 11.00$ shows the moment that the bubble fully encloses the obstacle;
    \item at around the time $t = 15.00$ the rear of the bubble reaches the tip of the triangular obstacle; 
    \item the energy then increases until $t\approx 20.00$ which is the time when the bubble leaves the obstacle.
\end{enumerate}
We note that gravity is an outer force that can add energy to the system. We thus do not have a dissipative system and especially the energy of the system can increase over time.

On the other hand, Figure~\ref{fig:energy-ineq}(b) plots 
\begin{align}
    {\Delta \mathcal{I} (m) := \mathcal{E}(\bu^{m-1},\varphi^{m-1}) + \tau\mathcal{F}(\bu^{m},\varphi^{m}) - 
    [\mathcal{E}(\bu^m,\varphi^m) + \tau \mathcal{D}(\bu^{m},\varphi^{m})  + \mathcal{P}(\bu^{m},\varphi^{m}) ]},
    \label{eq:num:energygap}
\end{align}
where 
${\mathcal{F}(\bu^{m},\varphi^{m}) = \left( \rho^{m}\bg^m, \bu^m \right)_\Omega}$ accounts for the contribution of the external force to the energy inequality, 
$\mathcal{D}(\bu^{m},\varphi^{m})  = 2\int_\Omega \eta^m|\D \bu^m|^2 \du x + b \int_\Omega \alphi{\psi}|\nabla\mu^m|^2 \du x + \int_\Omega \alpha_{\bu}^{(\delta)}(\psi)|\bu^m|^2\du x$ is the dissipative energy, and 
$\mathcal{P}(\bu^{m},\varphi^{m}) =   \int_\Omega \frac{\rho^{m-1}}{2}|\bu^m-\bu^{m-1}|^2\du x + \frac{\sigma\varepsilon}{2}\int_\Omega \alphi{\psi}|\nabla(\varphi^m-\varphi^{m-1})|^2 \du x $ accounts for the time-discretization errors. 
We see that $\Delta\mathcal{I}(m)$ is the right-hand side minus the left-hand side of the energy inequality \eqref{energy:totaldisc}. According to \eqref{energy:totaldisc}, $\Delta\mathcal{I}(m)$ should be positive at the continuous-in-space setting; and our numerical example illustrates that the same inequality still holds despite the space discretization and the mesh refinement at each time step.

\begin{figure}
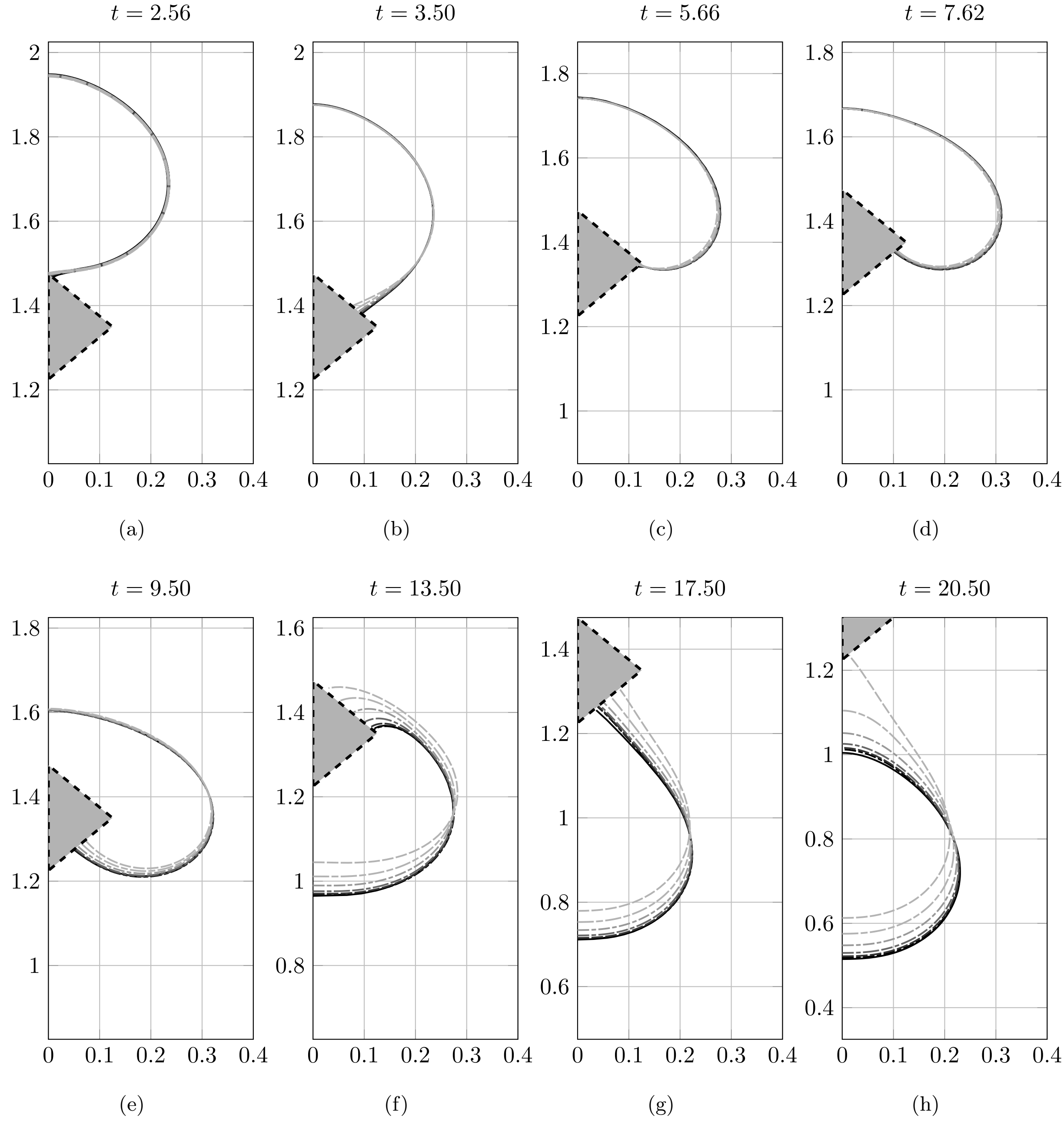

\begingroup
\pgfplotsset{    
    every axis/.style={
        width=4.7cm,
        height=8cm,
        },
    }
    \pgfplotsset{ 
         every axis legend/.append style={
         font=\tiny,
         draw=none,fill=none,text opacity=0, draw opacity=0,
         },
    }
    
    \centering
   
    \begin{subfigure}{0.245\textwidth}
         \centering
         \includestandalone{Figures/tikz/compvdelta/t256}
         \caption*{(a)}
    \end{subfigure}
    \begin{subfigure}{0.245\textwidth}
         \centering
         \includestandalone{Figures/tikz/compvdelta/t350}
         \caption*{(b)}
    \end{subfigure}
    \begin{subfigure}{0.245\textwidth}
         \centering
         \includestandalone{Figures/tikz/compvdelta/t566}
         \caption*{(c)}
    \end{subfigure}
    \begin{subfigure}{0.245\textwidth}
         \centering
         \includestandalone{Figures/tikz/compvdelta/t762}
         \caption*{(d)}
    \end{subfigure}
    
    \bigskip
    
    \begin{subfigure}{0.245\textwidth}
         \centering
         \includestandalone{Figures/tikz/compvdelta/t950}
         \caption*{(e)}
    \end{subfigure}
    \begin{subfigure}{0.245\textwidth}
         \centering
         \includestandalone{Figures/tikz/compvdelta/t1350}
         \caption*{(f)}
    \end{subfigure}
    \begin{subfigure}{0.245\textwidth}
         \centering
         \includestandalone{Figures/tikz/compvdelta/t1750}
         \caption*{(g)}
    \end{subfigure}
    \begin{subfigure}{0.245\textwidth}
         \centering
         \includestandalone{Figures/tikz/compvdelta/t2050}
         \caption*{(h)}
    \end{subfigure}
\endgroup
    \caption{The figure shows the zero-level sets of $\varphi_\delta$, with $\delta = 1.0\times10^{-2}, 5.0\times10^{-3}, 2.5\times10^{-3},1.2\times 10^{-3}, 6.0\times 10^{-4}$ and $3.0\times10^{-4}$, and of $\varphi_S$. The level sets for $\varphi_\delta$ are shown with broken lines which are colored such that the darker the line is the smaller $\delta$ is, while the solid black line indicates the sharp setting.
    We show these data 
    at the indicated times (a) $t = 2.56$, (b) $t = 3.50$, (c) $t = 5.66$, (d) $t = 7.62$, (e) $t = 9.50$, (f) $t = 13.50$, (g) $t = 17.50$ and (h) $t = 20.50$. We observe that the zero-level sets exhibit distinction from each other the more the droplet interacts with the obstacle. We also visually 
    observe the convergence of $\varphi_\delta$ towards $\varphi_S$ as $\delta \to 0$, which is in total agreement with the results in Section \ref{sec:asymptotics}.}
    \label{fig:delta-comp}
    
\end{figure}

\paragraph{Investigation for  $\delta\to 0$.}\mbox{}\\
In Section~\ref{sec:asymptotics} we have analytically  investigated the behavior of our model for $\delta \to 0$.
We repeat the above simulation with varying values of 
$\delta = 1.0\times10^{-2}, 5.0\times10^{-3}, 2.5\times10^{-3},1.2\times 10^{-3}, 6.0\times 10^{-4}$ and $3.0\times10^{-4}$.
To serve as a reference, we also run a simulation where we physically remove the triangular obstacle from the domain and solve \eqref{system:prime} in remaining domain. We shall refer to this simulation as the \textit{sharp setting} and denote by $\varphi_S$ the corresponding phase-field.
Figure \ref{fig:delta-comp} shows the behavior of the zero-level sets of the concentration $\varphi_{\delta}$, with varying values of $\delta$, and of the zero-level set of the concentration $\varphi_S$.

We observe that the zero-level sets of the concentrations are almost identical until they approach the obstacle. 
Figure~\ref{fig:delta-comp}~(a) shows that the zero-level set of $\varphi_S$ is already traversing the obstacle while the zero-level sets of the $\varphi_\delta$'s are just about to touch the vertex of the triangle. 
This can be attributed to the fact that the interface of $\psi$ causes loss to the relative concentration mass as we have discussed in the previous set up.  
As the liquid interacts more with the obstacle, the distinction between the phase field $\varphi_\delta$, for varying $\delta$, and $\varphi_S$ becomes more apparent, see Figure~\ref{fig:delta-comp}~(b). 
We also observe from Figure~\ref{fig:delta-comp}~(c), that the phenomena of \textit{nearly pinning} on the vertex of the obstacle causes the contact between the zero-level sets and the obstacle to catch each other. Nevertheless, we can still see some differences on the zero-level sets that do not touch the obstacle, perhaps due to the loss of relative mass. 
As the droplets overcome the \textit{pinning} phenomena, the distinction between the level sets becomes more obvious, see Figures \ref{fig:delta-comp} (d-h). What is noteworthy about the simulations illustrated in Figure \ref{fig:delta-comp} is the ordering of the level sets. We notice that as we decrease the value of $\delta$, the zero-level sets of $\varphi_\delta$ get closer to that of $\varphi_S$, this result agrees with the asymptotic analysis done in Section \ref{sec:asymptotics}.

\begin{figure} 
    \centering
     {
    \pgfplotsset{
        every axis/.style={
        width=0.98\textwidth,
        height=6cm
        }
        }
    \definecolor{mycolor1}{rgb}{0.0, 0.0, 0.0}%
\definecolor{mycolor2}{rgb}{0.2,0.2,0.2}%
\definecolor{mycolor3}{rgb}{0.4,0.4,0.4}%
\definecolor{mycolor4}{rgb}{0.6,0.6,0.6}%
\definecolor{mycolor5}{rgb}{0.7,0.7,0.7}%
\definecolor{mycolor6}{rgb}{0.8,0.8,0.8}%
\definecolor{mycolor6}{rgb}{0.3882353, 0.3176471, 0.0000000}%
\definecolor{mycolor7}{rgb}{0.0039216, 0.3764706, 0.1764706}
\begin{tikzpicture}

\pgfmathsetmacro{\DataSpacing}{50}

\begin{axis}[%
xmin=-0.1,
xmax=25.1,
xlabel = {$t$},
ymin=-2.3045,
ymax=-2.2855,
ylabel = {$\mathcal{M}$},
ytick={-2.312, -2.308, -2.304, -2.300, -2.296, -2.292, -2.288, -2.284},
yticklabels ={$-2.312$,$-2.308$,$-2.304$,$-2.300$,$-2.296$,$-2.292$,$-2.288$, $-2.284$},
axis background/.style={fill=white},
axis x line*=bottom,
axis y line*=left,
grid = major,
axis lines=box,
legend pos = south west,
]

\addplot [color=mycolor5, line width=1.25pt, dash pattern=on 5pt off 1pt]
  table[row sep=crcr]{%
0.01	-2.28749\\
0.51	-2.28759\\
1.01	-2.28765\\
1.51	-2.28767\\
2.01	-2.28769\\
2.51	-2.28773\\
3.01	-2.28902\\
3.51	-2.2905\\
4.01	-2.2918\\
4.51	-2.29276\\
5.01	-2.29308\\
5.51	-2.29327\\
6.01	-2.29331\\
6.51	-2.29339\\
7.01	-2.29346\\
7.51	-2.29377\\
8.01	-2.29402\\
8.51	-2.29442\\
9.01	-2.29508\\
9.51	-2.29566\\
10.01	-2.29643\\
10.51	-2.29733\\
11.01	-2.2981\\
11.51	-2.29887\\
12.01	-2.29961\\
12.51	-2.30013\\
13.01	-2.30036\\
13.51	-2.30009\\
14.01	-2.29894\\
14.51	-2.29697\\
15.01	-2.29642\\
15.51	-2.2961\\
16.01	-2.29577\\
16.51	-2.29543\\
17.01	-2.29497\\
17.51	-2.29453\\
18.01	-2.29396\\
18.51	-2.29333\\
19.01	-2.29278\\
19.51	-2.29209\\
20.01	-2.29153\\
20.51	-2.29105\\
21.01	-2.29059\\
21.51	-2.2906\\
22.01	-2.2906\\
22.51	-2.2906\\
23.01	-2.29062\\
23.51	-2.29061\\
24.01	-2.29061\\
24.51	-2.29065\\
25.01	-2.29063\\
};
\addlegendentry{$\delta = 1.0\times 10^{-2}$}

\addplot [color=mycolor5, , line width=1.25pt, dash pattern=on 7pt off 1pt on 3pt off 1pt]
  table[row sep=crcr]{%
0.01	-2.28752\\
0.51	-2.28763\\
1.01	-2.28768\\
1.51	-2.28771\\
2.01	-2.28773\\
2.51	-2.28775\\
3.01	-2.28866\\
3.51	-2.28952\\
4.01	-2.29025\\
4.51	-2.29063\\
5.01	-2.29076\\
5.51	-2.29088\\
6.01	-2.291\\
6.51	-2.29114\\
7.01	-2.29131\\
7.51	-2.2915\\
8.01	-2.29178\\
8.51	-2.29217\\
9.01	-2.2926\\
9.51	-2.29307\\
10.01	-2.2936\\
10.51	-2.29413\\
11.01	-2.29465\\
11.51	-2.29513\\
12.01	-2.29542\\
12.51	-2.29564\\
13.01	-2.29574\\
13.51	-2.29496\\
14.01	-2.29382\\
14.51	-2.29338\\
15.01	-2.29327\\
15.51	-2.29314\\
16.01	-2.2929\\
16.51	-2.29261\\
17.01	-2.29226\\
17.51	-2.29184\\
18.01	-2.2914\\
18.51	-2.29097\\
19.01	-2.29059\\
19.51	-2.29026\\
20.01	-2.28999\\
20.51	-2.29\\
21.01	-2.29001\\
21.51	-2.29001\\
22.01	-2.29001\\
22.51	-2.29001\\
23.01	-2.29\\
23.51	-2.29\\
24.01	-2.29\\
24.51	-2.29\\
25.01	-2.28999\\
};
\addlegendentry{$\delta = 5.0\times 10^{-3}$}

\addplot [color=mycolor4, , line width=1.25pt, dash pattern=on 5pt off 1pt on 2pt off 1pt]
  table[row sep=crcr]{%
0.01	-2.28753\\
0.51	-2.28764\\
1.01	-2.28769\\
1.51	-2.28772\\
2.01	-2.28773\\
2.51	-2.28775\\
3.01	-2.28831\\
3.51	-2.28878\\
4.01	-2.28925\\
4.51	-2.28947\\
5.01	-2.28959\\
5.51	-2.28968\\
6.01	-2.28978\\
6.51	-2.28988\\
7.01	-2.29\\
7.51	-2.29017\\
8.01	-2.29039\\
8.51	-2.29063\\
9.01	-2.29089\\
9.51	-2.29118\\
10.01	-2.2915\\
10.51	-2.29182\\
11.01	-2.29211\\
11.51	-2.2923\\
12.01	-2.29243\\
12.51	-2.29252\\
13.01	-2.29247\\
13.51	-2.29194\\
14.01	-2.29137\\
14.51	-2.29135\\
15.01	-2.29129\\
15.51	-2.29116\\
16.01	-2.29098\\
16.51	-2.2908\\
17.01	-2.29062\\
17.51	-2.29044\\
18.01	-2.29022\\
18.51	-2.29002\\
19.01	-2.28982\\
19.51	-2.28969\\
20.01	-2.28969\\
20.51	-2.2897\\
21.01	-2.2897\\
21.51	-2.28969\\
22.01	-2.28969\\
22.51	-2.28969\\
23.01	-2.28969\\
23.51	-2.2897\\
24.01	-2.28969\\
24.51	-2.28968\\
25.01	-2.28968\\
};
\addlegendentry{$\delta = 2.5\times 10^{-3}$}
\addplot [color=mycolor3, , line width=1.25pt, dash pattern=on 5pt off 1pt on 2pt off 1pt on 7pt off 1pt]
  table[row sep=crcr]{%
0.01	-2.28753\\
0.51	-2.28764\\
1.01	-2.28769\\
1.51	-2.28772\\
2.01	-2.28773\\
2.51	-2.28775\\
3.01	-2.28803\\
3.51	-2.28828\\
4.01	-2.2885\\
4.51	-2.28861\\
5.01	-2.28867\\
5.51	-2.28872\\
6.01	-2.28875\\
6.51	-2.28881\\
7.01	-2.28887\\
7.51	-2.28895\\
8.01	-2.28904\\
8.51	-2.28916\\
9.01	-2.28929\\
9.51	-2.28942\\
10.01	-2.28955\\
10.51	-2.28968\\
11.01	-2.28971\\
11.51	-2.2898\\
12.01	-2.28984\\
12.51	-2.28989\\
13.01	-2.28983\\
13.51	-2.28952\\
14.01	-2.28934\\
14.51	-2.28932\\
15.01	-2.28927\\
15.51	-2.28921\\
16.01	-2.28912\\
16.51	-2.28903\\
17.01	-2.28894\\
17.51	-2.28884\\
18.01	-2.28876\\
18.51	-2.28867\\
19.01	-2.28858\\
19.51	-2.28859\\
20.01	-2.2886\\
20.51	-2.2886\\
21.01	-2.2886\\
21.51	-2.2886\\
22.01	-2.2886\\
22.51	-2.28859\\
23.01	-2.28859\\
23.51	-2.28859\\
24.01	-2.28858\\
24.51	-2.28858\\
25.01	-2.28857\\
};
\addlegendentry{$\delta = 1.2\times 10^{-3}$}
\addplot [color=mycolor2, , line width=1.25pt, dash pattern=on 5pt off 1pt]
  table[row sep=crcr]{%
0.01	-2.28753\\
0.51	-2.28764\\
1.01	-2.28769\\
1.51	-2.28772\\
2.01	-2.28773\\
2.51	-2.28774\\
3.01	-2.28789\\
3.51	-2.288\\
4.01	-2.28811\\
4.51	-2.28817\\
5.01	-2.28822\\
5.51	-2.28825\\
6.01	-2.28828\\
6.51	-2.28831\\
7.01	-2.28833\\
7.51	-2.28837\\
8.01	-2.28841\\
8.51	-2.28846\\
9.01	-2.28852\\
9.51	-2.28857\\
10.01	-2.28863\\
10.51	-2.28867\\
11.01	-2.28861\\
11.51	-2.28864\\
12.01	-2.28864\\
12.51	-2.28865\\
13.01	-2.28863\\
13.51	-2.28846\\
14.01	-2.28838\\
14.51	-2.28835\\
15.01	-2.28831\\
15.51	-2.28827\\
16.01	-2.28822\\
16.51	-2.28817\\
17.01	-2.28812\\
17.51	-2.28807\\
18.01	-2.28802\\
18.51	-2.28797\\
19.01	-2.28794\\
19.51	-2.28795\\
20.01	-2.28795\\
20.51	-2.28796\\
21.01	-2.28796\\
21.51	-2.28795\\
22.01	-2.28795\\
22.51	-2.28795\\
23.01	-2.28795\\
23.51	-2.28794\\
24.01	-2.28794\\
24.51	-2.28793\\
25.01	-2.28793\\
};
\addlegendentry{$\delta = 6.0\times 10^{-4}$}
\addplot [color=mycolor1, , line width=1.25pt, dash pattern=on 5pt off 1pt on 2pt off 1pt]
  table[row sep=crcr]{%
0.01	-2.28754\\
0.51	-2.28764\\
1.01	-2.28769\\
1.51	-2.28772\\
2.01	-2.28773\\
2.51	-2.28774\\
3.01	-2.28781\\
3.51	-2.28786\\
4.01	-2.28792\\
4.51	-2.28796\\
5.01	-2.28799\\
5.51	-2.28801\\
6.01	-2.28803\\
6.51	-2.28804\\
7.01	-2.28805\\
7.51	-2.28806\\
8.01	-2.28808\\
8.51	-2.28809\\
9.01	-2.28811\\
9.51	-2.28812\\
10.01	-2.28813\\
10.51	-2.28812\\
11.01	-2.28801\\
11.51	-2.28802\\
12.01	-2.28801\\
12.51	-2.28801\\
13.01	-2.28803\\
13.51	-2.28794\\
14.01	-2.28789\\
14.51	-2.28786\\
15.01	-2.28784\\
15.51	-2.28781\\
16.01	-2.28778\\
16.51	-2.28774\\
17.01	-2.28771\\
17.51	-2.28768\\
18.01	-2.28766\\
18.51	-2.28763\\
19.01	-2.28763\\
19.51	-2.28763\\
20.01	-2.28763\\
20.51	-2.28763\\
21.01	-2.28763\\
21.51	-2.28763\\
22.01	-2.28763\\
22.51	-2.28763\\
23.01	-2.28762\\
23.51	-2.28762\\
24.01	-2.28762\\
24.51	-2.28761\\
25.01	-2.28761\\
};
\addlegendentry{$\delta = 3.0\times 10^{-4}$}

\addplot[only marks, mark=square*, mark options={fill=mycolor1, draw=mycolor1}, mark size=2pt, draw=mycolor2, forget plot] table[row sep=crcr]{%
x	y\\
2.51	-2.28773\\
};


\addplot[only marks, mark=triangle*, mark options={fill=mycolor5, draw=mycolor5}, mark size=3pt, draw=mycolor3, forget plot] table[row sep=crcr]{%
x	y\\
13.01	-2.30036\\
};

\addplot[only marks, mark=triangle*, mark options={fill=mycolor5, draw=mycolor5}, mark size=3pt, draw=mycolor2, forget plot] table[row sep=crcr]{%
x	y\\
13.01	-2.29574\\
};

\addplot[only marks, mark=triangle*, mark options={fill=mycolor4, draw=mycolor4}, mark size=3pt, draw=mycolor2, forget plot] table[row sep=crcr]{%
x	y\\
12.51	-2.29252\\
};
%
\addplot[only marks, mark=triangle*, mark options={fill=mycolor3, draw=mycolor3}, mark size=3pt, draw=mycolor2, forget plot] table[row sep=crcr]{%
x	y\\
12.51	-2.28989\\
};
%
\addplot[only marks, mark=triangle*, mark options={fill=mycolor2, draw=mycolor2}, mark size=3pt, draw=mycolor2, forget plot] table[row sep=crcr]{%
x	y\\
10.51	-2.28867\\
};
%
\addplot[only marks, mark=triangle*, mark options={fill=mycolor1, draw=mycolor1}, mark size=3pt, draw=mycolor2, forget plot] table[row sep=crcr]{%
x	y\\
10.01	-2.28813\\
};


\addplot[only marks, mark=*, mark options={fill=mycolor5, draw=mycolor5}, mark size=2pt, draw=mycolor4, forget plot] table[row sep=crcr]{%
x	y\\
21.01	-2.29059\\
};

\addplot[only marks, mark=*, mark options={fill=mycolor5, draw=mycolor5}, mark size=2pt, draw=mycolor4, forget plot] table[row sep=crcr]{%
x	y\\
20.01	-2.28999\\
};

delta 2.5e-3 
\addplot[only marks, mark=*, mark options={fill=mycolor4, draw=mycolor4}, mark size=2pt, draw=mycolor4, forget plot] table[row sep=crcr]{%
x	y\\
19.51	-2.28969\\
};
%
\addplot[only marks, mark=*, mark options={fill=mycolor3, draw=mycolor3}, mark size=2pt, draw=mycolor4, forget plot] table[row sep=crcr]{%
x	y\\
19.01	-2.28858\\
};
%
\addplot[only marks, mark=*, mark options={fill=mycolor2, draw=mycolor2}, mark size=2pt, draw=mycolor4, forget plot] table[row sep=crcr]{%
x	y\\
19.01	-2.28794\\
};
%
\addplot[only marks, mark=*, mark options={fill=mycolor1, draw=mycolor1}, mark size=2pt, draw=mycolor4, forget plot] table[row sep=crcr]{%
x	y\\
18.51	-2.28763\\
};

\draw[dashed, mycolor1] (0, -2.28753) -- (25, -2.28753);

\end{axis}

\end{tikzpicture}%
    }
    \caption{The figure shows the the evolution of $\mathcal{M}$ over time for varying values of $\delta$.
    We observe that the variation in $\mathcal{M}$ is less pronounced for smaller values of $\delta$ which corroborates the contention that the quantity $\mathcal{M}$ is conserved as $\delta\to 0$.
    As in Figure~\ref{fig:mass-trend} the square and the circle mark the first and last contact of the droplet with the obstacle and the triangle marks the time point with largest spill-over, when the obstacle is (nearly) fully covered by the droplet.
    }
    \label{fig:masses}
\end{figure}
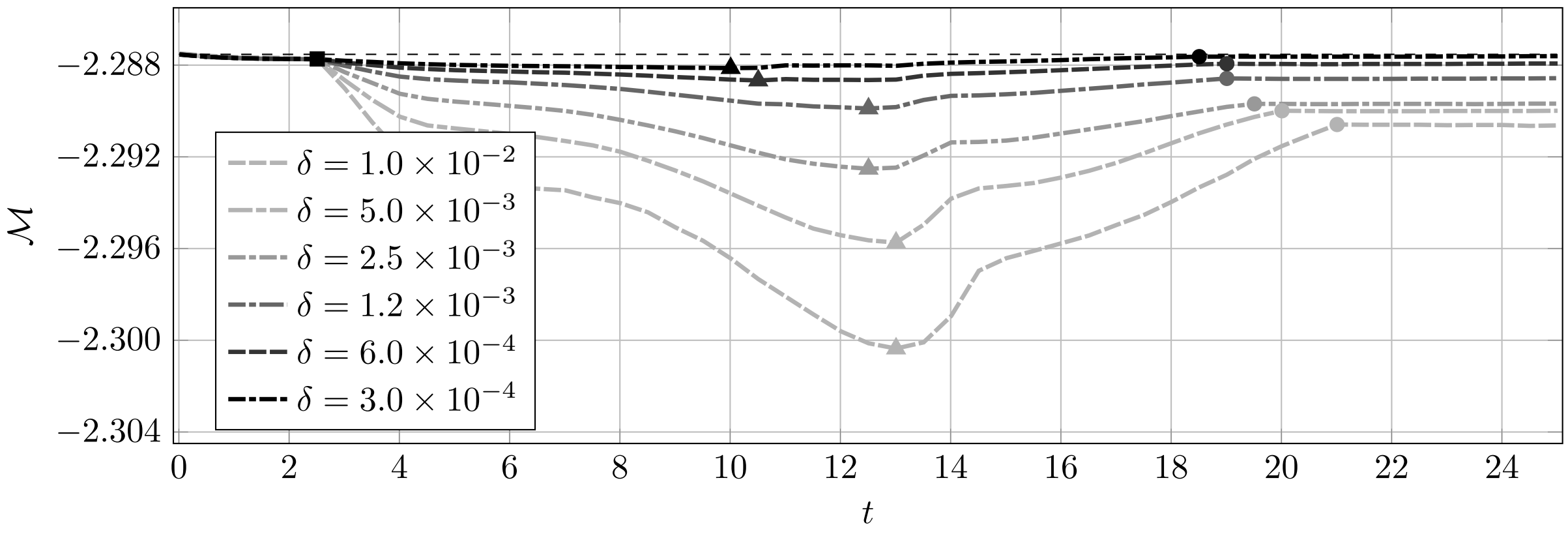

In Figure~\ref{fig:masses} we show the evolution of $\mathcal M$ over time for these simulations.
We observe a reduction in the variation of the quantity $\mathcal{M}$ as $\delta>0$ becomes smaller. This supports the observation that as $\delta\to 0$, we formally recover the conservation of the average of the phase-field $\varphi$. Although not explicitly shown, Figure~\ref{fig:masses} illustrates how the \emph{spill-over} is also reduced as we make $\delta>0$ smaller. This can be observed on the values of $\mathcal{M}$ at the beginning and end of the time interval $[0,25]$, 
i.e., as $\delta$ becomes sufficiently small, 
the amount of the droplet that stays inside the obstacle becomes nearly negligible, confirming improved conservation.

\section*{Conclusion}

In this work, we proposed and analyzed a diffuse-interface model for two-phase flow in the presence of obstacles. By incorporating a porous-medium approximation of the obstacle into a thermodynamically consistent Cahn--Hilliard/Navier--Stokes framework, we obtained a unified formulation that is well-suited for applications in shape and topology optimization.

We performed formal matched asymptotics, which demonstrated that the proposed model recovers the 
sharp-obstacle model
in the limit as the interface parameter associated with the obstacle tends to zero. 
This provides theoretical justification for the diffuse representation of obstacles and confirms the consistency of the approach with the underlying physical model.

We introduced a time-discrete scheme designed to preserve key structural properties of the continuous system, including stability and energy dissipation. We established the existence of solutions to the semi-discrete system and derived global energy estimates, ensuring the robustness of the proposed discretization. The numerical experiments corroborate the theoretical findings, illustrating convergence of the diffuse obstacle approximation toward the sharp-interface limit as the regularization parameter decreases. Moreover, the simulations highlight qualitative features such as interface dynamics and interaction with obstacles, as well as the influence of the diffuse interface thickness on mass conservation and flow behavior.

Overall, the proposed framework provides a mathematically consistent and computationally feasible approach for modeling two-phase flows in implicitly defined domains. It serves as a foundation for future developments in PDE-constrained optimization involving binary fluids. Possible directions for further research include the rigorous analysis of the fully discrete scheme, the incorporation of more complex physical effects such as contact-line dynamics, and the extension to optimization problems with additional constraints and objective functionals.
    
	\section*{Statements and Declarations}
	
	\subsubsection*{Conflict of Interest.}
	The authors declare that they have no conflict of interest.
	
	\subsubsection*{Funding.}
	The authors acknowledge funding of the project \textit{Fluiddynamische Formoptimierung mit Phasenfeldern und Lipschitz-Methoden} by the German Research Foundation under the project number \href{https://gepris.dfg.de/gepris/projekt/543959359?language=en}{543959359}.
	
	\subsubsection*{Author Contribution.}
	(CRediT taxonomy)
	\begin{description}
		\item[M.H.]  Conceptualization, Funding Acquisition,  Methodology, Project Administration, Resources, Supervision,  Writing - Original Draft, Writing - Review \& Editing
		\item[C.K.]  Conceptualization, Funding Acquisition, Methodology, Project Administration, Resources, Supervision, Writing - Original Draft, Writing - Review \& Editing
		\item[J.S.H.S.]  Conceptualization, Data Curation, Methodology, Project Administration, Software, Visualization, Writing - Original Draft, Writing - Review \& Editing
	\end{description}

	\subsubsection*{Acknowledgement.}
	The authors acknowledge funding by the German Research Foundation within the project \textit{Fluiddynamische Formoptimierung mit Phasenfeldern und Lipschitz-Methoden} under the project number \href{https://gepris.dfg.de/gepris/projekt/543959359?language=en}{543959359}.

	\subsubsection*{Data availability.}
 Our implementation in \texttt{FreeFEM++} is documented in the repository
 \mbox{\href{https://github.com/jhsimon/ThermoConsistentCHNS}{jhsimon/ThermoConsistentCHNS}.}

 \subsubsection*{Disclosure of AI  Usage.}
 AI tools were used exclusively for linguistic editing and optimization to improve clarity, consistency and readability during the final stage of the writing process. After using any AI tool, the authors reviewed and edited the content as needed and take full responsibility for the content of the publishes article.

	
	
	\begingroup
	\printbibliography
	\endgroup

\end{document}